\documentclass{amsart}
\usepackage{amssymb, amsmath, latexsym}

\renewcommand{\baselinestretch}{\baselinestretch}
\renewcommand{\baselinestretch}{1.1}
\author{Constantin-Nicolae Beli}
\title{Analogues of the $p^n$th Hilbert symbol in characteristic $p$ (updated)} 
\date{}

   \def\m{\lim}
   
\def\p{\partial}     
    \def\te{\theta}
   
 \def\({\overline}

\def\){\underline} \def\<{\cdot} \def\go{\mathfrak}
\def\>{~~~~~~~} \def\#{{\bf
Definition}} \def\*{\section} \def\be{\begin{equation}}
\def\ee{\end{equation}}

\def\sb{\subset} \def\sp{\supset} \def\sbq{\subseteq} \def\spq{\supseteq} 
\def\ti{\times}  \def\oo{{\cal O}} 
   \def\FF{{\mathbb F}}
 \def\ff{\dot{F}} \def\ooo{{\oo^\ti}} 
 \def\mo{{\rm mod}~}  
  \def\fs{\ff^2}  
\def\p{\go p}    
\def\*{\sharp}  \def\0{} 
 \def\1{^{-1}}  
 \def\[{\prec} \def\]{\succ} 
\def\bmat{\left(\begin{array}} \def\emat{\end{array}\right)} 
   
\def\N{{\rm N}}  
  
 \def\m2{~(\mo 2)} \def\no{\noindent}
 \def\btm{\begin{thm}}
\def\etm{\end{tm}}
 \def\blem{\begin{lem}}
\def\elem{\end{lem}}
\def\kk{K^\times}

\newtheorem{theorem}{Theorem}[section]
\newtheorem{proposition}[theorem]{Proposition}
\newtheorem{lemma}[theorem]{Lemma}
\newtheorem{definition}{Definition}
\newtheorem{corollary}[theorem]{Corollary}

\newtheorem{bof}[theorem]{}
\newtheorem{teorema}{Theorem}

\def\qed{\mbox{$\Box$}\vspace{\baselineskip}}
\def\pf{$Proof.\,\,$} 
\def\bco{\begin{corollary}} \def\eco{\end{corollary}} 
\def\bdf{\begin{definition}} \def\edf{\end{definition}} 
\def\btm{\begin{theorem}} \def\etm{\end{theorem}} 
\def\bpr{\begin{proposition}} \def\epr{\end{proposition}}  
\def\blm{\begin{lemma}} \def\elm{\end{lemma}} 
\def\bff{\begin{bof}\rm} \def\eff{\end{bof}}
\def\btr{\begin{teorema}} \def\etr{\end{teorema}}

\def\de{\newcommand} \de\tm[1]{{\no\bf Theorem~#1}} 
 \def\mb{\mathbb} 
 \def\QQ{{\mb Q}}  \def\ZZ{{\mb Z}}
\def\NN{{\mb N}} \def\FF{{\mb F}}
\def\II{{\mathcal I}}
\de\lm[1]{{\no\bf Lemma~#1}}
\de\df[1]{{\no\bf Definition~#1}} \de\co[1]{{\no\bf Corollary~#1}}
\de\tp[1]{\te (#1 )} \de\ts[1]{\te (O^-(#1 ))} \de\ty[1]{\te
(O(#1 ))} \de\tx[1]{\te (#1 )} \de\up[1]{(1+\p^{#1} )\fs}
 \de\ups[1]{((1+\p^{#1})\fs )^*} \de\upo[1]{(1+\p^{#1} )\ooo^2}
\de\upon[2]{(1+\p^{#1})\ooo^2\cap\N (#2 )}
\de\lr[1]{\longrightarrow^{\!\!\!\!\!\!\!\! #1}}
\de\lf[1]{\longleftarrow^{\!\!\!\!\!\!\!\! #1}}
\de\si[1]{\sim^{\!\!\!\!\! #1}} \de\apr[1]{\approx^{\!\!\!\!\! #1}}
\de\leg[2]{\left(\frac {#1}{#2}\right)}
\DeclareMathOperator\Br{Br}

\de\Brr[1]{{}_{#1}\Br}
\de\kkk[1]{K^{\times{#1}}}
\newcommand{\diff}{\mathop{}\!\mathrm{d}}
\DeclareMathOperator\car{char}
\DeclareMathOperator\ad{ad}
\DeclareMathOperator\dlog{dlog}

\begin{document}

\maketitle

\begin{abstract}

The $p$th degree Hilbert symbol $(\cdot,\cdot
)_p:\kk/\kkk p\times\kk/\kkk p\to\Brr p(K)$ from characteristic $\neq
p$ has two analogues in characteristic $p$, 
$$[\cdot,\cdot )_p:K/\wp (K)\times\kk/\kkk p\to\Brr p(K),$$
where $\wp$ is the Artin-Schreier map $x\mapsto x^p-x$, and 
$$((\cdot,\cdot ))_p:K/K^p\times K/K^p\to\Brr p(K).$$

The symbol $[\cdot,\cdot )_p$ generalizes to an analogue of
$(\cdot,\cdot )_{p^n}$ via the Witt vectors,
$$[\cdot,\cdot )_{p^n}:W_n(K)/\wp
(W_n(K))\times\kk/\kkk{p^n}\to\Brr{p^n}(K).$$

Here $W_n(K)$ is the truncation of length $n$ of the ring of
$p$-typical Witt wectors, i.e. $W_{\{1,p,\ldots,p^{n-1}\}}(K)$.

In this paper we construct similar generalizations for $((\cdot,\cdot
))_p$. Our construction involves Witt vectors and Weyl algebras. In
the process we obtain a new kind of Weyl algebras in characteristic
$p$, with many interesting properties.

The symbols we introduce, $((\cdot,\cdot ))_{p^n}$ and, more
generally, $((\cdot,\cdot ))_{p^m,p^n}$, which here are defined in
terms of central simple algebras, coincide with the homonymous symbols
we introduced in [B] in terms of the symbols $[\cdot,\cdot
)_{p^n}$. This will be proved in a future paper. In the present paper
we only introduce the symbols and we prove that they have the same
properties with the symbols from [B]. These properties are enough to
obtain the representation theorem for $\Brr{p^n}(K)$ from [B, Theorem
4.10]. 

\end{abstract}



\section{Introduction}

If $A$ is a central division algebra (c.s.a.) over a field $K$ then we
denote by $[A]$ its class in the Brauer group $(\Br (K),+)$. We have
$[A]+[B]=[A\otimes B]$, $0=[K]=\{ M_n(K)\mid n\geq 1\}$ and
$-[A]=[A^{op}]$. We denote by $\Brr n(K)$ the $n$-torsion of $\Br
(K)$. 

From now on we assume that $\car K=p$. We denote by $F$ the
Frobenius map, $x\mapsto x^p$ and by $\wp =F-1$ the
Artin-Schreier map, $x\mapsto x^p-x$. 

Reall that if $\car K\neq p$ and $\mu_p\sb K$ then we have the
bilinear and skew-symmetric Hilbert symbol $(\cdot,\cdot )_p:\kk/\kkk
p\times\kk/\kkk p\to\Brr p(K)$. In characteristic $p$, besides
$(\kk/\kkk p,\cdot)$, we have two more groups, $(K/\wp (K),+)$ and
$(K/K^p,+)$. These three greoups are involved in two bilinear symbols
with values in $\Brr p(K)$, which are analogues of $(\cdot,\cdot
)_p$. 

The symbol $[\cdot,\cdot )_p=[\cdot,\cdot )_{K,p}:K/\wp
(K)\ti\kk/\kkk p\to\Brr p(K)$ is given by $[a,b)_p=[A_{[a,b)_p}]$, where
$A_{[a,b)_p}$ is a c.s.a. of degree $p$ over $K$ generated by $x,y$, with
the relations $x^p-x=a$, $y^p=b$ and $yxy^{-1}=x+1$,
i.e. $yx=xy+y$. The symbol $[\cdot,\cdot )_p$ is bilinear. 

The symbol $((\cdot,\cdot ))_p=((\cdot,\cdot ))_{K,p}:K/K^P\ti
K/K^p\to\Brr p(K)$ is given by $((a,b))_p=[A_{((a,b))_p}]$, where
$A_{((a,b))_p}$ is a c.s.a. of degree $p$ over $K$ generated by $x,y$,
with the relations $x^p=a$, $y^p=b$ and $[y,x]:=yx-xy=1$.

The symbol $((\cdot,\cdot ))_p$ is bilinear and skew-symmetric. It
also has the property $((ab,c))_p+((bc,a))_p+((ca,b))_p=0$ $\forall
a,b,c\in K$. This enables us to define linear map
$\alpha_p:\Omega^1(K)/\diff K\to\Brr p(K)$ by $a\diff b\mapsto
((a,b))_p$. Here $\Omega^1(K)$ is the $K$-module generated by $\diff a$
with $a\in K$, subject to $\diff (a+b)=\diff a+\diff b$ $\diff
(ab)=a\diff b+b\diff a$ for $a,b\in K$. 

Unlike $[\cdot,\cdot )_p$, the symbol $((\cdot,\cdot ))_p$ is not
widely used. We found this notation in [KMRT, page 25], but only when
$p=2$. The properties of $((\cdot,\cdot ))_p$ listed above appear in
[BK1, 8.1.1], where ${\mathcal A}_{f,g}$ is used to denote
$A_{((f,g))_p}^{op}$. 

The symbols $[\cdot,\cdot )_p$ and $((\cdot,\cdot ))_p$ are related by the
relations $((a,b))_p=[ab,b)_p$ if $b\neq 0$ and
$((a,0))_p=0$. Therefore the symbol $((\cdot,\cdot ))_p$ defined here
is the same with the symbol from [B]. (See [B, Remark 3.1(2)].)

We didn't find the formula $((a,b))_p=[ab,b)_p$ in the literature so
we prove it here. We will produce an isomorphism $f$ between
$A_{[a,b)_p}=K\langle z,u\mid z^p-z=a,\, u^p=b,\, uz=zu+u\rangle$ and
$A_{((a,b))_p}=K\langle x,y\mid x^p=a,\, y^p=b,\, [y,x]=1\rangle$. We
take $f$ with $f(z)=xy$ and $f(u)=y$. To prove that there is a
morphism with these properties we must show that $f$ preserves the
relations between generators, i.e. that $(xy)^p-xy=a$, $y(xy)=(xy)y+y$
and $y^p=b$. We already have the third relation and for the second we
just note that $yxy-xyy=[y,x]y=y$. For the first relation we note that
$[\cdot,x]$ is a derivation so $[y,x]=1$ implies $[y^n,x]=ny^{n-1}$,
i.e. $y^nx=xy^n+ny^{n-1}$ for $n\geq 1$. It follows that
$x^{n+1}y^{n+1}=x^n(xy^n)y=x^n(y^nx-ny^{n-1})y=x^ny^nxy-nx^ny^n=x^ny^n(xy-n)$. Then,
by induction on $n$, we get $x^ny^n=(xy)(xy-1)\cdots (xy-n+1)$. In 
particular, $ab=x^py^p=(xy)(xy-1)\cdots (xy-p+1)=(xy)^p-(xy)$. So $f$
is defined and it is obviously surjective ($f(zu^{-1})=x$ and
$f(u)=y$). Since $A_{[ab,b)_p}$ is a c.s.a. $f$ is an isomorphism.


The symbol $[\cdot,\cdot )_p$ generalizes to a symbol with values
in the $p^n$-torsion of the Brauer group via Witt vectors. Namely, we
have a symbol $[\cdot,\cdot )_{p^n}=[\cdot,\cdot )_{K,p^n}:W_n(K)/\wp
(W_n(K))\times\kk/\kkk{p^n}\to\Brr{p^n}(K)$, where the
Artin-Schreier map $\wp$ is defined on $p$-typical Witt vectors by
$\wp =F-1$, i.e. if $x=(x_0,\ldots,x_{n-1})\in W_n(K)$ then
$\wp (x)=Fx-x=(x_0^p,\ldots,x_{n-1}^p)-(x_0,\ldots,x_{n-1})$. If
$a=(a_0,\ldots,a_{n-1})\in W_n(K)$, $b\in\kk$ then
$[a,b)_p=[A_{[a,b)_p}]$, where $A_{[a,b)_p}$ is a c.s.a. of degree $p^n$
generated by $x=(x_0,...,x_{n-1})$ and $y$, where $x_0,\ldots,x_{n-1}$
commute with each other, with the relations $\wp (x)=a$, $y^{p^n}=b$
and $yxy^{-1}=x+1$. Here $x$ is regarded as a Witt vector and in the
last relation $yxy^{-1}:=(yx_0y^{-1},\ldots,yx_{n-1}y^{-1})$ and $x+1$
is a sum of Witt vectors, $x+1=(x_0,\ldots,x_{n-1})+(1,0,\ldots,0)$. 



In this paper we will produce similar generalizations for
$((\cdot,\cdot ))_p$. Our construction involves Weyl algebras and Witt
vectors. As a by-product, we construct a new class of Weyl algebras in
characteristic $p$. In a future paper we will prove that the symbols
$((\cdot,\cdot ))_{p^m,p^n}$ we introduce here are the same with the
ones from [B]. For now, we only prove they have the same properties.

\section{Universal $B$ algebra}

Throughout this paper $\NN$ denotes $\NN_0=\ZZ_{\geq 0}$ and $\NN^*$
denotes $\NN_1=\ZZ_{\geq 1}$.

Unless otherwise specified, all rings are assumed to be commutative,
with unity.

By $[\cdot,\cdot ]$ we denote the commutator, $[a,b]=ab-ba$. For every
$a$ the map $[a,\cdot]$ is a derivation, i.e. $[a,b_1\cdots
b_n]=\sum_{i=1}^nb_1\cdots b_{i-1}[a,b_i]b_{i+1}\cdots b_n$. Similarly
for $[\cdot,b]$. In particular, if $[a,b]=1$ then for any $n\geq 1$ we
have $[a,b^n]=nb^{n-1}$ and $[a^n,b]=na^{n-1}$. 
\medskip

For every $S\sbq\NN^*$ we denote by $S^{-1}=\{ n^{-1}\mid n\in
S\}$. Note that $\ZZ [S^{-1}]=\ZZ [p^{-1}\mid p\text{ prime, }\exists
n\in S, p\mid n]$. A ring $R$ has a structure of $\ZZ [S^{-1}]$-ring
iff $S\sbq R^\times$ or, equivalently, iff $p\in R^\times$ for all
primes $p$ dividing elements from $S$. 

If $S\sbq\NN^*$ then we denote
$$\II_S=\{ i=(i_n)_{n\in S}\in\NN^S\mid
i_n=0\text{ for almost all }n\in S\}.$$
If $x=(x_n)_{n\in S}$, where all $x_n$ commute with each other, and
$i=(i_n)_{n\in S}\in\II_S$ then we define $x^i=\prod_{n\in
  S}x_n^{i_n}$. Since $i_n=0$ for almost all $n$ this is a finite
product. When $S=\NN^*$ we denote $\II =\II_{\NN^*}$. 

On $\II_S$ we define the lexicographic order $\leq$ as follows. If
$i=(i_n)_{n\in S},\, j=(j_n)_{n\in S}\in\II_S$ we say that $i<j$ if
there is $n\in S$ such that $i_n<j_n$ and $i_k=j_k$ for $k<n$. Note
that for any $S\sbq\NN^*$ we can regard $(\II_S,\leq )$ as a subset of
$(\II,\leq )$ by identifying $(i_n)_{n\in S}\in\II_S$ with 
$(i_n)_{n\geq 1}\in\II$, where $i_n:=0$ for $n\in\NN^*\setminus S$. If
$S=\emptyset$ we put $I_\emptyset =\{ 0\}$. 

If $x=(x_n)_{n\in S}$ has commuting entries and $i=(i_n)_{n\in
S}\in\II_S$ then we denote by $x^i=\prod_{n\in S}x_n^{i_n}$. Since
$i_n=0$ for almost all $n$ this is a finite product. If $S=\emptyset$
then by $x=(x_n)_{n\in\emptyset}$ we mean an empty sequence of length
zero and the set $\{ x^i\mid i\in\II_\emptyset\}$ is just $\{ 1\}$.

If $n\in\NN\cup\{\infty\}$ we put $\II_n=\II_{\{ 1,\ldots,n\}}$. When
$n=0$ by $\{ 1,\ldots,n\}$ we mean $\emptyset$ so $\II_0=\II_\emptyset
=\{ 0\}$. If $n=\infty$ then $\{ 1,\ldots,n\}$ means $\NN^*$ so
$\II_\infty=\II$.

Note that any $S\sbq\NN^*$ can be written as $S=\{ s_1<s_2<\cdots
<s_n\}$ for some $n\in\NN\cup\{ \infty\}$ so $(S,\leq )\cong\{
1,\ldots,n\}$. Also if $x=(x_k)_{k\in S}$ and $i=(i_k)_{k\in
  S}\in\II_S$ then the product $x^i$ writes as $x^i=y^j$, where
$y=(y_1,\ldots,y_n)$, with $y_k=x_{s_k}$, and
$j=(j_1,\ldots,j_n)\in\II_n$, where $j_k=i_{s_k}$. So we can restrict
ourselves to subsets $S$ of $\NN^*$ of the type $\{1,\ldots,n\}$ for
some $n\in\NN\cup\{\infty\}$. However, as we will see later, we need
sequences $x$ indexed by subsets of $\NN^*$ that are not of this
form. 
\medskip

If $R$ is a ring, $n\in\NN\cup\{\infty\}$ and $X=(X_1,\ldots,X_n)$
then $X^i$ with $i\in\II_n$ are a basis for
$R[X]=R[X_1,\ldots,X_n]$. (If $n=0$ then $R[X_1,\ldots,X_n]:=R$.)
A polynomial in $R[X]$ has the form $P=\sum_{i\in\II_n}a_iX^i$, where
$a_i\in R$ are almost all zero and we define $\deg_XP=\max\{
i\in\II_n\mid a_i\neq 0\}$ and $\deg_X0=-\infty$. Then
$\deg_X:R[X]\to\II_n\cup\{ -\infty\}$ has the usual properties of the
degree: $\deg_X(P+Q)\leq\max\{\deg_XP,\deg_XQ\}$ and
$\deg_XPQ\leq\deg_XP+\deg_XQ$, with equality when $R$ is an integral
domain.
\medskip

If $R$ is a ring, $m,n\in\NN\cup\{\infty\}$, $X=(X_1,\ldots,X_m)$ and
$Y=(Y_1,\ldots,Y_n)$ then we define $R\langle [X],[Y]\rangle=R\langle
[X_1,\ldots,X_m],[Y_1,\ldots,Y_n]\rangle$ as
$$R\langle [X],[Y]\rangle =R\langle X,Y\mid [X_i,X_j]=0~\forall i,j,\,
[Y_i,Y_j]=0~\forall i,j\rangle.$$
Let $\{ X\} =\{ X_1,\ldots,X_m\}$ and $\{ Y\} =\{
Y_1,\ldots,Y_n\}$. We call a word in $R\langle [X],[Y]\rangle$ a
product $Z=Z_1\cdots Z_k$ with $Z_h\in\{ X\}\cup\{ Y\}$. We denote
$Z_X=\prod_{Z_h\in\{ X\}}Z_h$ and $Z_Y=\prod_{Z_h\in\{ Y\}}Z_h$ and we
define $\deg_XZ:=\deg_XZ_X$ and $\deg_YZ:=\deg_YZ_Y$. Note that every
word writes uniquely as $Z=X^{i_1}Y^{j_1}\cdots X^{i_s}Y^{j_s}$ for
some $s\geq 1$, $i_h\in\II_m$, $j_h\in\II_n$, with $i_h\neq 0$ if
$h>1$, $j_h\neq 0$ if $h<s$. Then $Z_X=X^{i_1+\cdots +i_s}$,
$Z_Y=Y^{j_1+\cdots +j_s}$ so $\deg_XZ=i_1+\cdots +i_s$,
$\deg_YZ=j_1+\cdots +j_s$. The words are a basis for $R\langle
[X],[Y]\rangle$. If $P\in R\langle [X],[Y]\rangle$, $P=\sum_Za_ZZ$,
where $Z$ covers all words and $a_Z\in R$ are almost all zero, then we
define $\deg_XP=\max\{ \deg_XZ\mid a_Z\neq 0\}$ and $\deg_YP=\max\{
\deg_YZ\mid a_Z\neq 0\}$. 
\medskip

We denote by $R[X][Y]$ the submodule of $R\langle [X],[Y]\rangle$
generated by the words $X^iY^j$ with $i\in\II_m$, $j\in\II_n$. An
element in $R[X][Y]$ has the form
$\sum_{i\in\II_m,j\in\II_n}a_{i,j}X^iY^j$, where $a_{i,j}\in R$ and
$a_{i,j}=0$ for almost all $i,j$. We have $\deg_XP=\max\{
i\in\II_m\mid\exists j\in\II_n,\, a_{i,j}\neq 0\}$ and $\deg_YP=\max\{
j\in\II_n\mid\exists i\in\II_m,\, a_{i,j}\neq 0\}$. 

Note that $R[X][Y]$ is the image of the injective linear map $\mu
:R[X]\otimes_RR[Y]\to R\langle [X][Y]\rangle$, given by $P(X)\otimes
Q(Y)\mapsto P(X)Q(Y)$. 
\medskip

For convenience, if $C$ is an $R$-algebra and $x_1,\ldots,x_n\in C$ we
say that $C=R[x_1,\ldots,x_n]$ {\em strictly} if $C\cong
R[X_1,\ldots,X_n]$ relative to the generators
$x_1,\ldots,x_n$. Similarly, we say that $C=R\langle
x_1,\ldots,x_n\rangle$ {\em strictly} if $C$ is freely generated by
$x_1,\ldots,x_n$. 
\medskip

For every $n\in\NN^*$ we put
$$D(n)=\{d\in\NN^*\, :\, d\mid n\}\text{ and }D^*(n)=D(n)\setminus\{
n\}.$$

Recall that the elements of in $W$, the ring of universal Witt
vectors, write as $x=(x_1,x_2,\ldots)$. The ghost functions $w_n$ are
defined as $w_n(x)=\sum_{d\mid n}dx_d^{n/d}$. Over $\QQ$ the sum and
the product of the Witt vectors $x=(x_1,x_2,\ldots )$,
$y=(y_1,y_2,\ldots )$ are given by $x+y=z=(z_1,z_2,\ldots )$ and
$xy=t=(t_1,t_2,\ldots )$, where $z$ and $t$ satisfy
$w_n(z)=w_n(x)+w_n(y)$ and $w_n(t)=w_n(x)w_n(y)$. One proves easily
that $z_n=s_n(x,y)$, $t_n=p_n(x,y)$ for some $s_n,p_n\in\QQ
[X_d,Y_d\mid d\in D(n)]$. But it turns out that in fact $s_n$ and
$p_n$ have coefficients in $\ZZ$. This allows the definition of the
ring of Witt vectors to be extended over arbitrary rings by defining
$x+y=(s_1(x,y),s_2(x,y),\ldots )$, $xy=(p_1(x,y),p_2(x,y),\ldots )$.

We also consider truncation sets, i.e. subsets $P$ of $\NN^*$ with the
property that if $n\in P$ then $D(n)\sbq P$. They appear in the
definiton of the truncated Witt vectors $W_P$, whose elements have the
form $(x_n)_{n\in P}$. The operations on $W_p$ are defined the same
way as for the universal Witt vectors. If $x=(x_n)_{n\in P}$,
$y=(y_n)_{n\in P}$ then $x+y=(s_n(x,y))_{n\in P}$, $xy=(p_n(x,y))_{n\in
P}$. 

If $P,Q$ are truncation sets with $Q\sbq P$ and $x=(x_n)_{n\in P}$ is
a Witt vector from $W_P$ then we denote by $x_Q$ its truncation in
$W_Q$, $x_Q=(x_n)_{n\in Q}$.

\bdf If $R$ is a $\QQ$-ring then we define $B(R)$, the universal $B$
algebra over $R$, as the $R$-algebra generated by $x=(x_1,x_2,\ldots
)$ and $y=(y_1,y_2,\ldots )$, with the relations
$[w_m(x),w_n(x)]=[w_m(y),w_n(y)]=0$ and
$[w_n(y),w_m(x)]=\delta_{m,n}m$ $\forall m,n\in\NN^*$.

More generally, if $P,Q$ are truncation sets and $R$ is a $\ZZ
[P^{-1},Q^{-1}]$-ring, i.e. with $P,Q\sbq R^\times$, we define the
algebra $B_{P,Q}(R)$ generated by $x=(x_m)_{m\in P}$ and
$y=(y_n)_{n\in Q}$, with the relations $[w_m(x),w_n(x)]=0$ $\forall
m,n\in P$, $[w_m(y),w_n(y)]=0$ $\forall m,n\in Q$ and
$[w_n(y),w_m(x)]=\delta_{m,n}m$, $\forall m\in P,\, n\in Q$.

If $P=Q$ we denote $B_P(R)=B_{P,P}(R)$. In particular,
$B(R)=B_{\NN^*}(R)$. 
\edf
{\bf Remark.} If $P,Q,P',Q'$ are truncation sets with $P'\sbq P$ and
$Q'\sbq Q$ then any $\ZZ [P^{-1},Q^{-1}]$-ring $R$ is also a $\ZZ
[P'^{-1},Q'^{-1}]$-ring. Also the generators of $B_{P',Q'}(R)$ are
amongst the generators of $B_{P,Q}(R)$ and the relations among
generators in $B_{P',Q'}(R)$ also hold in $B_{P,Q}(R)$. So we have a
cannonical morphism $B_{P',Q'}(R)\to B_{P,Q}(R)$.

\blm Let $C$ be an algebra over a ring $R$ and let
$n\in\NN\cup\{\infty\}$. Assume that $x=(x_1,\ldots,x_n)$ and
$y=(y_1,\ldots,y_n)$ are two sequences with entries in $C$ such
that for $1\leq i\leq n$ we have $y_i\in a_ix_i+\langle
x_1,\ldots,x_{i-1}\rangle$ for some $a_i\in R^\times$. 

(i) For every $1\leq i\leq n$ we have $x_i\in a_i^{-1}y_i+\langle
y_1,\ldots,y_{i-1}\rangle$. As a consequence, $\langle
x_1,\ldots,x_n\rangle =\langle y_1,\ldots,y_n\rangle$. In particular,
$x_1,\ldots,x_n$ commute with each other iff $y_1,\ldots,y_n$ do so.

(ii) We have $C=R\langle x_1,\ldots,x_n\rangle$ strictly iff
$C=R\langle y_1,\ldots,y_n\rangle$ strictly.

(iii) We have $C=R[x_1,\ldots,x_n]$ strictly iff $C=R[y_1,\ldots,y_n]$
strictly.
\elm
\pf (i) Note that the second statement follows from the first by
double inclusion. We use induction on $i$. If $i=1$ by hypothesis
$y_1=a_1x_1+b$ for some $b\in R$. It follows that
$x_1=a_1^{-1}y_1-a_1^{-1}b$ and we are done. Assume now that (i) holds
for indices $<i$. Then $y_1,\ldots,y_{i-1}$ can be written in terms of
$x_1,\ldots,x_{i-1}$ and vice versa so $\langle
x_1,\ldots,x_{i-1}\rangle =\langle y_1,\ldots,y_{i-1}\rangle$. Then
$y_i\in a_ix_i+\langle x_1,\ldots,x_{i-1}\rangle =a_ix_i+\langle
y_1,\ldots,y_{i-1}\rangle$ implies that $x_i\in a_i^{-1}y_i+\langle
y_1,\ldots,y_{i-1}\rangle$.

(ii) In the view (i) our statement is symmetric in $x$ and $y$ so we
will prove only that if $C=R\langle x_1,\ldots,x_n\rangle$ strictly
then $C=R\langle y_1,\ldots,y_n\rangle$ strictly. We consider a free
algebra $C'=R\langle z_1,\ldots,z_n\rangle$. Then there is a unique
morphism of algebras $f:C'\to C$ given by $f(z_i)=y_i$. We have
$x_i\in a_i^{-1}y_i+\langle y_1,\ldots,y_{i-1}\rangle
=a_i^{-1}f(z_i)+\langle f(z_1),\ldots,f(z_{i-1})\rangle
=f(a_i^{-1}z_i+\langle z_1,\ldots,z_{i-1}\rangle )$. Hence there is
$t_i\in a_i^{-1}z_i+\langle z_1,\ldots,z_{i-1}\rangle$ with
$f(t_i)=x_i$. We denote $z=(z_1,\ldots,z_n)$,
$t=(t_1,\ldots,t_n)$. Since $C=R\langle x_1,\ldots,x_n\rangle$
strictly there is a unique morphism of algebras $g:C\to C'$ with
$g(x_i)=t_i$. Then for any $i$ we have $f(g(x_i))=f(t_i)=x_i$ so
$f\circ g=1_C$. Since $C'=R\langle z_1,\ldots,z_n\rangle$ strictly,
$t_i\in a_i^{-1}z_i+\langle z_1,\ldots,z_{i-1}\rangle$ and
$g(x_i)=t_i$ we have that $C'$, $z$, $t$, $C$, $x$ and $g$ are in the
same situation as $C$, $x$, $y$, $C'$, $z$ and $f$. Then, by the same
reasoning as for $f$, we get that $g$ too admits an inverse to the
right, i.e. there is $h:C'\to C$ such that $g\circ h=1_{C'}$. Together
with $f\circ g=1_C$, this implies that $f=h$ so $f$ and $g$ are
inverse to each other isomorphisms. Since $C'$ is freely generated by
$z_1,\ldots,z_n$ this implies that $C$ is freely generated by
$f(z_1),\ldots,f(z_n)$, i.e. by $y_1,\ldots,y_n$.

(iii) is similar to (ii) but this time we define
$C'=R[z_1,\ldots,z_n]$ and for the construction of $f$ and $g$ we use
the universal property for polynomial algebras instead of free
algebras. \qed

\bco If $R$ is a $\ZZ [P^{-1},Q^{-1}]$-ring then the condtitons
$[w_m(x),w_n(x)]=0$ $\forall m,n\in P$ and $[w_m(y),w_n(y)]=0$
$\forall m,n\in Q$ from the definition of $B_{P,Q}(R)$ are equivalent
to $[x_m,x_n]=0$ $\forall m,n\in P$ and $[y_m,y_n]=0$ $\forall m,n\in
Q$, respectively. 

Also $\langle w_m(x)\mid m\in P\rangle =\langle x_m\mid m\in P\rangle$
and  $\langle w_n(y)\mid n\in Q\rangle =\langle y_n\mid n\in
Q\rangle$.
\eco
\pf Since $R$ is a $\ZZ [P^{-1},Q^{-1}]$-ring we have $P,Q\sbq
R^\times$. For every $m\in P$ we have $w_m(x)=\sum_{d\mid
m}dx_d^{m/d}\in mx_m+\langle x_k\mid k\in P,\, k<m\rangle$ and $m\in
R^\times$. Hence $[w_m(x),w_n(x)]=0$ $\forall m,n\in P$ is equivalent
to $[x_m,x_n]=0$ $\forall m,n\in P$ and $\langle w_m(x)\mid m\in
P\rangle =\langle x_m\mid m\in P\rangle$ by Lemma 2.1(i). Similarly
for the equivalence between $[w_m(y),w_n(y)]=0$ $\forall m,n\in Q$ and
$[y_m,y_n]=0$ $\forall m,n\in Q$ and for $\langle w_n(y)\mid n\in
Q\rangle =\langle y_n\mid n\in Q\rangle$. \qed 
\bigskip

\blm Let $C$ be an algebra over a ring $R$ and let
$m,n\in\NN\cup\{\infty\}$. Let $x=(x_1,\ldots,x_m)$,
$y=(y_1,\ldots,y_n)$ be sequences in $C$ such that $x_i$'s commute
with each other and $y_j$'s commute with each other. The following are
equivalent:

(i) $x^iy^j$ with $i\in\II_m$, $j\in\II_n$ are a basis of $C$.

(ii) Every element $\alpha\in C$ writes uniquely as $\alpha
=P(x_1,\ldots,x_m,y_1,\ldots,y_n)$ for some $P\in
R[X_1,\ldots,X_m][Y_1,\ldots,Y_n]$. 

(iii) If $C'$ and $C''$ are the subalgebras of $C$ generated by
$x_1,\ldots,x_m$ and $y_1,\ldots,y_n$, respectively, then
$C'=R[x_1,\ldots,x_m]$ and $C''=R[y_1,\ldots,y_n]$ strictly and the
linear map $\mu :C'\otimes_RC''\to C$, given by
$\alpha\otimes\beta\mapsto\alpha\beta$, is bijective.
\elm
\pf The equivalence between (i) and (ii) is trivial.

The condition that $C'=R[x_1,\ldots,x_m]$ strictly from (iii) is
equivalent to the fact that $x^i$ with $i\in\II_m$ are linearly
independent, i.e. that they are a basis for $C'$. But this is a
consequence of (i). Similarly, $C''=R[y_1,\ldots,y_n]$ strictly means
that $y^j$ with $j\in\II_n$ are a basis of $C''$ and this is a
consequence of (i). Assuming that the two conditions are fulfilled,
$x^i\otimes y^j$ with $i\in\II_m$, $j\in\II_n$ are a basis for
$C'\otimes_RC''$. Then the condition that $\mu$ is a bijection is
equivalent to the fact that $\mu (x^i\otimes y^j)=x^iy^j$, with
$i\in\II_m$, $j\in\II_n$, are a basis for $C$, i.e. it is equivalent to
(i). \qed

Let $P,Q\sbq\NN^*$ be truncation sets and let $R$ be a ring with
$P,Q\sbq R^\times$. Let $x=(x_m)_{m\in P}$, $y=(y_n)_{n\in Q}$. If
$m\in P$, $n\in Q$ then $m,n\in R^\times$. We denote $z_m=w_m(x)$ and
$t_n=n^{-1}w_n(y)$, which is defined since $n\in R^\times$. We have
$z_m\in mx_m+\langle x_k\mid k\in P,\, k<m\rangle$, with $m\in
R^\times$ and $t_n\in y_n+\langle y_l\mid l\in Q,\, l<n\rangle$, with
$1\in R^\times$. Hence if $z=(z_m)_{m\in P}$ and $t=(t_n)_{n\in Q}$
then Lemma 2.1 applies both to $x$ and $z$ and to $y$ and $t$.

We now prove that $B_{P,Q}(R)$ writes in terms of the Weyl
algebras. Recall that the $N$-th Weyl algebra $A_N(R)$ is the
$R$-algebra generated by $x=(x_1,\ldots,x_N)$ and $y=(y_1,\ldots,y_N)$
with the relations $[x_m,x_n]=[y_m,y_n]=0$ and
$[y_n,x_m]=\delta_{m,n}$ for $1\leq m,n\leq N$. It has the property
that $x^iy^j$ with $i,j\in\II_N$ form a basis for $A_N(R)$.

We write $A_N(R)(x,y)$ if we want to specify the generators $x$ and
$y$. 

\blm (i) Let $z'=(z_m)_{m\in P\cap Q}$, $z''=(z_m)_{m\in P\setminus
Q}$, $t'=(t_n)_{n\in P\cap Q}$, $t''=(t_n)_{n\in Q\setminus
P}$. Then $B_{P,Q}(R)=A_{|P\cap Q|}(z',t')[z'',t'']$ strictly.

In particular, if $P=Q=\NN^*$ then $B(R)=A_\infty (z,t)$. 

(ii) $z^it^j$ with $i\in\II_P$, $j\in\II_Q$ are a basis of
$B_{P,Q}(R)$ over $R$.
\elm
\pf (i) By Lemma 2.1(ii) the free algebras $R\langle x\rangle$ and
$R\langle y\rangle$ are freely generated by $z$ and $t$. Hence
$R\langle x,y\rangle$ is freely generated by $z$ and $t$. Since
$z_m=w_m(x)$, $t_n=n^{-1}w_n(y)$ the relations $[w_m(x),w_n(x)]=0$
$[w_m(y),w_n(y)]=0$ and $[w_n(y),w_m(x)]=\delta_{m,n}m$ write as
$[z_m,z_n]=0$, $[t_m,t_n]=0$ and $[t_n,z_m]=\delta_{m,n}$. Hence
$B_{P,Q}(R)=R\langle z,t\mid [z_m,z_n]=0,\, [t_m,t_n]=0,\,
[z_m,t_n]=\delta_{m,n}\rangle$. Note that the only pairs of generators
that do not commute are $z_n,t_n$ with $n\in P\cap Q$, when we have
$[t_n,z_n]=1$, so they involve only the entries of $z'$ and $t'$. The
relations among generators involving only the entries of $z'$ and
$t'$ are $[z_m,z_n]=[t_m,t_n]=0$ and $[z_m,t_n]=\delta_{m,n}$ $\forall
m,n\in P\cap Q$ and they define the algebra $A_{|P\cap
Q|}(z',t')$. The relations involving $z_m$ with $m\in P\setminus Q$
and $t_n$ with $n\in Q\setminus P$, i.e. the entries of $z''$ and
$t''$, are the commutativity relations with all the other generators. It
follows that $B_{P,Q}(R)=C[z'',t'']$ strictly, where $C=R\langle
z',t'\mid [z_m,z_n]=[t_m,t_n]=0,\,
[t_n,z_m]=\delta_{m,n}\rangle=A_{|P\cap Q|}(z',t')$.

(ii) By the theory of Weyl algebras $z'^at'^b$, $a,b\in\II_{P\cap Q}$,
are a basis of $C=A_{|P\cap Q|}(z',t')$ over $R$. Since
$B_{P,Q}(R)=C[z'',t'']$ strictly $z''^ct''^d$, $c\in\II_{P\setminus Q}$,
$d\in\II_{Q\setminus P}$, are a basis of $B_{P,Q}(R)$ over $C$. Hence
$z'^at'^bz''^ct''^d=(z'^az''^c)(t'^bt''^d)$, with $a,b\in\II_{P\cap
Q}$, $c\in\II_{P\setminus Q}$, $d\in\II_{Q\setminus P}$, are a basis
of $B_{P,Q}(R)$ over $R$. But $\{ z'^az''^c\mid a\in\II_{P\cap Q},\,
c\in\II_{P\setminus Q}\} =\{ z^i\mid i\in\II_P\}$ and $\{
t'^bt''^d\mid b\in\II_{P\cap Q},\, d\in\II_{Q\setminus P}\} =\{
t^j\mid j\in\II_Q\}$. Hence the conclusion. \qed

\blm $x^iy^j$ with $i\in\II_P$, $j\in\II_Q$ are a basis of
$B_{P,Q}(R)$ over $R$. 
\elm
\pf By Lemma 2.4(ii) and Lemma 2.3 if $C'=\langle z_m\mid m\in
P\rangle$, $C''=\langle t_n\mid n\in Q\rangle$ then $C'=R[z_m\mid m\in
P]$ and $C''=R[t_n\mid n\in Q]$ strictly and the multiplication map $\mu
:C'\otimes C''\to B_{P,Q}(R)$ is a bijection.

Recall that Lemma 2.1 applies to $x$ and $z$ and to $y$ and $t$. By
Lemma 2.1(i) we get $C'=\langle x_m\mid m\in P\rangle$ and
$C''=\langle y_n\mid n\in Q\rangle$ and by Lemma 2.1(iii)
$C'=R[x_m\mid m\in P]$ and $C''=R[y_n\mid n\in Q]$ strictly. Together
with the bijectivity of $\mu$, by Lemma 2.3 this implies that $x^iy^j$, with
$i\in\II_P$, $j\in\II_Q$, are a basis of $B_{P,Q}(R)$. \qed

\blm Let $R$ be a ring and let $m,n\in\NN\cup\{\infty\}$. Let
$\overline C=R\langle [X],[Y]\rangle$, where $X=(X_1,\ldots,X_m)$ and
$Y=(Y_1,\ldots,Y_n)$. Let $C$ be an $R$-algebra generated by
$x=(x_1,\ldots,x_m)$ and $y=(y_1,\ldots,y_n)$ such that in $C$ we have
$[x_i,x_j]=0$, $[y_i,y_j]=0$ and $[y_j,x_i]\in\langle
x_1,\ldots,x_{i-1},y_1,\ldots,y_{j-1}\rangle$.

We denote by $f:\overline C\to C$ the surjective morphism of algebras
given by $X_i\mapsto x_i$ and $Y_j\mapsto y_j$, i.e. $f(P)=P(x,y)$
$\forall P\in\overline C$. 

(i) $C=f(R[X_1,\ldots,X_m][Y_1,\ldots,Y_n])=
R[x_1,\ldots,x_m][y_1,\ldots,y_n]$. Equivalently, $C$ is spanned by
$x^iy^j$ with $i\in\II_m,\, j\in\II_n$.

(ii) For every word $Z$ of $\overline C$ we have $f(Z-Z_XZ_Y)=f(P)$
for some $P\in R[X][Y]$ with $\deg_XP<\deg_XZ$, $\deg_YP<\deg_YZ$.
\elm
\pf Note that the condition $[x_i,y_j]\in\langle
x_1,\ldots,x_{i-1},y_1,\ldots,y_{j-1}\rangle$ means that
$[x_i,y_j]=f(P_{i,j})$ for some $P_{i,j}\in R\langle
X_1,\ldots,X_{i-1},Y_1,\ldots,Y_{j-1}\rangle$.

If $Z$ is a word of $\overline C$ satisfying (ii) then let $P\in
R[X][Y]$ with $\deg_XP<\deg_XZ$, $\deg_YP<\deg_YZ$ such that
$f(Z-Z_XZ_Y)=f(P)$. It follows that $f(Z)=f(Q)$, where $Q=Z_XZ_Y+P\in
R[X][Y]$. Moreover, since $\deg_XP<\deg_XZ=\deg_XZ_XZ_Y$ and
$\deg_YP<\deg_YZ=\deg_YZ_XZ_Y$ we have $\deg_XQ=\deg_XZ$ and
$\deg_YQ=\deg_YZ$. In particular, $f(Z)=f(Q)\in f(R[X][Y])$.

Hence if (ii) holds then $f(Z)\in f(R[X][Y])$ for all words $Z$. Since
$\overline C$ is spanned by words this implies that $C=f(\overline
C)=f(R[X][Y])$, i.e. we have (i).

Now we prove (ii) by induction on $\deg_XZ\in\II_m$. (Recall that
$\II_m$ is a well-ordered set.) If $\deg_XZ=0$ then $Z_X=1$ and
$Z=Z_Y$ so $Z-Z_XZ_Y=0$ and we may take $P=0$. Suppose now that
$\deg_XZ>0$. Let $Z=Z_1\cdots Z_k$ with $Z_h\in\{ X\}\cup\{ Y\}$. We
denote by $M$ the set of all elements of $P\in \( C$
such that $f(P)=f(Q)$ for some $Q\in R[X][Y]$ such that
$\deg_XQ<\deg_XZ$, $\deg_YQ<\deg_YZ$. Then $M$ is an $R$-submodule of
$\overline C$. Note that if $f(P)=f(P')$ then $P\in M$ iff $P'\in
M$. We must prove that $Z-Z_XZ_Y\in M$.

If $Z'$ is a word in $\overline C$ with $\deg_XZ'<\deg_XZ$ and
$\deg_YZ'<\deg_YZ$ then by the induction hypothesis we have that (ii)
holds for $Z'$ so, by the reasoning above, $f(Z')=f(Q)$ for some $Q\in
R[X][Y]$ with $\deg_XQ=\deg_XZ'<\deg_XZ$ and
$\deg_YQ=\deg_YZ'<\deg_YZ$. It follows that $Z'\in M$. 

For any permutation $\sigma\in S_k$ we denote by $Z_\sigma =Z_{\sigma
(1)}\cdots Z_{\sigma (k)}$. Note that $(Z_\sigma )_X=Z_X$ and
$(Z_\sigma )_Y=Z_Y$ $\forall\sigma\in S_k$. We prove that all
$Z_\sigma$ are congruent modulo $M$. In particular, since $Z_XZ_Y$
writes as $Z_\sigma$ for some $\sigma$, we get $Z\equiv Z_XZ_Y\mod M$,
which proves our claim. Since $S_k$ is generated by transpositions, it
is enough to consider the case $\sigma =(r,r+1)$ for some $1\leq r\leq
k-1$. Then $Z-Z_\sigma =Z_1\cdots Z_k-Z_1\cdots
Z_{r-1}Z_rZ_{r+1}Z_{r+2}\cdots Z_k=Z_1\cdots
Z_{r-1}[Z_r,Z_{r+1}]Z_{r+2}\cdots Z_k$. If $Z_r,Z_{r+1}\in\{ X\}$ or
$Z_r,Z_{r+1}\in\{ Y\}$ then $[Z_r,Z_{r+1}]=0$ so $Z-Z_\sigma
=0$. Suppose that $Z_r=Y_j$, $Z_{r+1}=X_i$. Then 
$f([Z_r,Z_{r+1}])=[f(Z_r),f(Z_{r+1})]=[y_j,x_i]=f(P_{i,j})$. It
follows that $f(Z-Z_\sigma)=f(Z_1\cdots Z_{r-1}P_{i,j}Z_{r+2}\cdots
Z_k)$ so we must prove that $\alpha :=Z_1\cdots
Z_{r-1}P_{i,j}Z_{r+2}\cdots Z_k\in M$. Since $P_{i,j}\in\langle
X_1,\ldots,X_{i-1},Y_1,\ldots,Y_{j-1}\rangle$ it can be written as
a linear combination of words $T=T_1\cdots T_l$ with $T_h\in\{
X_1,\ldots,X_{i-1},Y_1,\ldots,Y_{j-1}\}$. Hence $\alpha$ writes as a
linear combination of words $Z'=Z_1\cdots Z_{r-1}TZ_{r+2}\cdots Z_k$,
with $T$ of this type. Then it suffices to prove that each such $Z'$
belongs to $M$. To do this we prove that $\deg_XZ'<\deg_XZ$ and
$\deg_YZ'<\deg_YZ$. But $Z'$ is obtained from $Z_X$ by removing
the factors $Z_rZ_{r+1}=Y_jX_i$ and replacing them by $T$. But $T_X$
is a product of factors from $\{ X_1,\ldots,X_{i-1}\}$ so
$\deg_XT=\deg_XT_X<\deg_XX_i=\deg_XY_jX_i$. Hence $\deg_X
Z'<\deg_XZ$. Similarly $\deg_YZ'<\deg_YZ$.

The case $Z_r=X_i$, $Z_{r+1}=Y_j$ is similar. (Here we have
$[Z_r,Z_{r+1}]=-[Y_j,X_i]$.)  \qed

\bco With the hypothesis of Lemma 2.6, for any $1\leq i\leq m$, $1\leq
j\leq n$ there is $c_{i,j}\in
R[X_1,\ldots,X_{i-1}][Y_1,\ldots,Y_{j-1}]$ such that
$[y_j,x_i]=c_{i,j}(x,y)$.
\eco
\pf We use Lemma 2.6(i) for $Z=Y_jX_i$. Then
$[y_j,x_i]=y_jx_i-x_iy_j=f(Y_jX_i-X_iY_j)=f(Z-Z_XZ_Y)=f(P)=P(x,y)$
for some $P\in R[X][Y]$ with $\deg_XP<\deg_XZ=\deg_XX_i$ and
$\deg_YP<\deg_YZ=\deg_YY_j$. But this simply means that $P\in
R[X_1,\ldots,X_{i-1}][Y_1,\ldots,Y_{j-1}]$. (If $T=X^aY^b$, with
$a\in\II_m$, $b\in\II_n$, is a word that appears with a nonzero
coefficient in $P$ then $\deg_XX^a=\deg_XT<\deg_XX_i$, which means
that $X^a$ is a product of factors from $\{X_1,\ldots,X_{i-1}\}$
only. Similarly, $Y^b$ is a product of factors from
$\{Y_1,\ldots,Y_{j-1}\}$ only.)

Hence we may take $c_{i,j}=P$. \qed

\blm Let $C$ be an $R$-algebra generated by $x=(x_1,\ldots,x_m)$ and
$y=(y_1,\ldots,y_n)$ such that $[x_i,x_j]=0$ $\forall i,j$,
$[y_i,y_j]=0$ $\forall i,j$ and for any $1\leq i\leq m$, $1\leq j\leq
n$ there are $\alpha,\beta,\gamma\in\langle
x_1,\ldots,x_{i-1},y_1,\ldots y_{j-1}\rangle$ such that
$[y_j,x_i]+[\beta,x_i]+[y_j,\alpha]+\gamma =0$.

Then there are $c_{i,j}\in R[X_1,\ldots,X_{i-1}][Y_1,\ldots,Y_{j-1}]$
such that for every $1\leq i\leq m$, $1\leq j\leq n$ we have
$[y_j,x_i]=c_{i,j}(x,y)$. Also $C$ is spanned by $x^iy^j$ with
$i\in\II_m$, $j\in\II_n$.
\elm
\pf By Corollary 2.7 it is enough to prove that $[y_j,x_i]\in\langle
x_1,\ldots,x_{i-1},y_1,\ldots y_{j-1}\rangle$ $\forall i,j$. We use
induction on $i+j$. If $i+j=2$, i.e. if $i=j=1$, then by hyopthesis
there are $\alpha,\beta,\gamma\in R$ with
$[y_1,x_1]+[\beta,x_1]+[y_1,\alpha ]+\gamma=0$, i.e. $[y_1,x_1]+\gamma
=0$. Hence $[y_1,x_1]=c_{1,1}$, where $c_{1,1}=-\gamma\in R$.

Suppose now that our statement is true when $i+j<N$. Let $i,j$ with
$i+j=N$. Let $\alpha,\beta,\gamma\in\langle
x_1,\ldots,x_{i-1},y_1,\ldots y_{j-1}\rangle$ such that
$[y_j,x_i]+[\beta,x_i]+[y_j,\alpha ]+\gamma =0$. To prove that
$[y_j,x_i]\in\langle x_1,\ldots,x_{i-1},y_1,\ldots y_{j-1}\rangle$
$\forall i,j$ it is enough to prove that $[\beta,x_i],[y_j,\alpha
]\in\langle x_1,\ldots,x_{i-1},y_1,\ldots y_{j-1}\rangle$. For
$[\beta,x_i]$, since $\beta$ is a linear combination of words
$z=z_1\cdots z_k$ with $z_h\in\{
x_1,\ldots,x_{i-1},y_1,\ldots,y_{j-1}\}$, it is enough to prove that
$[z,x_i]\in\langle x_1,\ldots,x_{i-1},y_1,\ldots y_{j-1}\rangle$ for
$z$ of this form. But we have $[z,x_i]=\sum_{h=1}^kz_1\cdots
z_{h-1}[z_h,x_i]z_{h+1}\cdots z_k$ so it is enough to prove that the
terms of this sum belong to $\langle x_1,\ldots,x_{i-1},y_1,\ldots
y_{j-1}\rangle$. But $z_1\cdots z_{h-1},z_{h+1}\cdots z_k\in\{
x_1,\ldots,x_{i-1},y_1,\ldots y_{j-1}\}$ so it suffices to prove
that $[z_h,x_i]\in\langle x_1,\ldots,x_{i-1},y_1,\ldots
y_{j-1}\rangle$. If $z_h=x_l$ for some $1\leq l\leq i-1$ then
$[z_h,x_i]=0$. If $z_h=y_l$ for some $1\leq l\leq j-1$ then
$i+l<i+j=N$ so, by the induction hypothesis, $[z_h,x_i]\in\langle
x_1,\ldots,x_{i-1},y_1,\ldots y_{l-1}\rangle\sbq\langle
x_1,\ldots,x_{i-1},y_1,\ldots y_{j-1}\rangle$. So we are done. The
relation $[y_j,\alpha ]\in\langle x_1,\ldots,x_{i-1},y_1,\ldots
y_{j-1}\rangle$ proves similarly. \qed 

\blm For every $m,n\in\NN^*$ there is $c_{m,n}\in\ZZ
[m^{-1},n^{-1}][X_d\mid d\in D^*(m)][Y_e\mid e\in D^*(n)]$ such that
for every truncation sets $P,Q$ with $m\in P$, $n\in Q$ and every $\ZZ
[P^{-1},Q^{-1}]$-ring $R$ in $B_{P,Q}(R)$ we have
$[y_n,x_m]=c_{m,n}(x,y)$.

In particular, $c_{1,1}=1$. 
\elm
\pf Take first the case $P=D(m)$, $Q=D(n)$ and $R=\ZZ
[D(m)^{-1},D(n)^{-1}]=\ZZ [m^{-1},n^{-1}]$. In $B_{D(m),D(n)}(\ZZ
[m^{-1},n^{-1}])$ for every $d\in D(m)$, $e\in D(n)$ we have
$w_d(x)=dx_d+a$ and $w_e(y)=ey_e+b$ for some $a\in\langle x_k\mid k\in
D(m),\, k<d\rangle$, $b\in\langle y_l\mid l\in D(n),\,
l<e\rangle$. Then $\delta_{d,e}d=[w_e(y),w_d(x)]=[ey_e+b,dx_d+a]$ so
$[y_e,x_d]+[\beta,x_d]+[y_e,\alpha ]+\gamma=0$, where $\alpha
=d^{-1}a$, $\beta =e^{-1}b$ and $\gamma
=d^{-1}e^{-1}[b,a]-\delta_{d,e}e^{-1}$. Obviously
$\alpha,\beta,\gamma\in\langle x_k\mid k\in D(m),\, k<d,\, y_l\mid
l\in D(n),\, l<e\rangle$. Since also, by Corollary 2.2, $x$ and $y$
have commuting entries, by Corollary 2.7 we get that for every $d\in
D(m)$, $e\in D(n)$ we have $[y_e,x_d]\in\ZZ[m^{-1},n^{-1}] [x_k\mid
k\in D(m),\, k<d][y_l\mid l\in D(n),\, l<e]$. In particular, when
$d=m$, $e=n$ there is $c_{m,n}\in\ZZ[m^{-1},n^{-1}][x_k\mid k\in
D(m),\, k<m][y_l\mid l\in D(n),\, l<n]=\ZZ[m^{-1},n^{-1}] [x_d\mid
d\in D(m)^*][y_e\mid e\in D(n)^*]$ such that $[y_n,x_m]=c_{m,n}(x,y)$.

Let now $P,Q$ be arbitrary truncation sets with $m\in P$, $n\in Q$ and
let $R$ be a $\ZZ [P^{-1},Q^{-1}]$-ring. Then $R$ is also a $\ZZ
[m^{-1},n^{-1}]$-ring. So the relation $[y_n,x_m]=c_{m,n}(x,y)$, which
holds in $B_{D(m),D(n)}(\ZZ [m^{-1},n^{-1}])$, will also hold in
$B_{D(m),D(n)}(R)=B_{D(m),D(n)}(\ZZ [m^{-1},n^{-1}])\otimes_{\ZZ
[m^{-1},n^{-1}]}R$. But $D(m)\sbq P$ and $D(n)\sbq Q$ so this
relation also holds in $B_{P,Q}(R)$. (See the Remark following
Definition 1.)

We have $[y_1,x_1]=[w_1(y),w_1(x)]=\delta_{1,1}\cdot 1=1$ so
$c_{1,1}=1$. \qed 

\blm Let $R$ be a ring and let $C=R\langle X\rangle/{\mathcal R}$,
$C'=R\langle X'\rangle/{\mathcal R}'$ with $X'\sbq X$ and ${\mathcal
R}'\sbq{\mathcal R}$, where $X,X'$ are sets of generators and
${\mathcal R}\sbq R\langle X\rangle$, ${\mathcal R}'\sbq R\langle
X'\rangle$ are the ideals of relations.

If there are $(v_i)_{i\in I}$ in $R\langle X'\rangle$ such that $v_i$
span $C'$ and they are linearly independent in $C$ then ${\mathcal
R}'={\mathcal R}\cap R\langle X'\rangle$ so $C'\sbq C$. 

If moreover $X'=X$ or $(v_i)_{i\in I}$ is a basis for $C$ then
$C'=C$.
\elm
\pf Since ${\mathcal R}'\sbq{\mathcal R}$ the map $f:C'\to C$,
$f(x)=x$ is well defined. (More precisely, $f$ is given by
$x+{\mathcal R}'\mapsto x+{\mathcal R}$ $\forall x\in R\langle
X'\rangle$.) 

Then $f$ is injective, so $C'\sbq C$, iff ${\mathcal R}'={\mathcal
R}\cap R\langle X'\rangle$. We must prove that $\ker f=0$. Let
$\alpha\in\ker f$. Then $\alpha$ writes as a linear combination
$\alpha =\sum_{h=1}^sa_hv_{i_h}$ with $a_h\in R$ and $i_h\in I$
mutually distinct. Then in $C$ we have $0=f(\alpha
)=\sum_{h=1}^sa_hv_{i_h}$. Since $v_i$ are linearly independent in $C$
we get $a_h=0$ $\forall h$ and so $\alpha =0$. 

If $X'=X$ or $(v_i)_{i\in I}$ is a basis for $C$ then $f$ is also
surjective. Thus it is a bijection, i.e. $C'=C$. \qed

\blm For any truncation sets $P,Q$ and any $\ZZ [P^{-1},Q^{-1}]$-ring
$R$ the algebra $B_{P,Q}(R)$ is generated by $x=(x_m)_{m\in P}$,
$y=(y_n)_{n\in Q}$, with the relations $[x_m,x_n]=0$ $\forall m,n\in
P$, $[y_m,y_n]=0$ $\forall m,n\in Q$ and $[y_n,x_m]=c_{m,n}(x,y)$
$\forall, m\in P,\, n\in Q$. (Here $c_{m,n}$ are those from Lemma
2.9.) 
\elm
\pf Let $C$ be the $R$-algebra generated by $x$ and $y$ with the
relations $[x_m,x_n]=0$, $[y_m,y_n]=0$ and
$[y_n,x_m]=c_{m,n}(x,y)$. The algebras $C$ and $B_{P,Q}(R)$ have the
same generators and the relations among generators in $C$ also hold
in $B_{P,Q}(R)$. In $C$ we have $[x_m,x_n]=0$, $[y_m,y_n]=0$ and
$[y_n,x_m]=c_{m,n}(x,y)\in\langle x_k\mid k\in P,\, k<m,\, y_l\mid
l\in Q,\, l<n\rangle$. By Lemma 2.6(i) $x^iy^j$ with $i\in\II_P$,
$j\in\II_Q$ span $C$ and by Lemma 2.5 they are a basis for
$B_{P,Q}(R)$. Hence $B_{P,Q}(R)=C$ by Lemma 2.10. \qed

By Lemma 2.11 the relations among generators in $B_{P,Q}$ are written
in terms of $c_{m,n}$, which so far have coefficients in $\QQ$. We
prove that in fact their coefiicients are integers so the definition
of $B_{P,Q}$ can be extended over arbitrary rings. It is a situation
similar to that from the theory of Witt vectors, where the polynomials
$p_n$ and $s_n$, which give the sum and the product in $W$, have a
priori rational coefficients but it turns out their coefficients are
integers. We will use the same series $\Lambda$ that is used in the
theory of Witt vectors to prove that $s_n$ and $p_n$ have integral
coefficients. 

If $x=(x_1,x_2,\ldots )$, where $x_1,x_2,\ldots$ commute with each
other, then we denote by $\Lambda (x;t)\in\ZZ [x][[t]]$ the formal
series $\Lambda (x;t)=\prod_n(1-x_nt^n)^{-1}$.

We have $\log\Lambda (x;t)=\sum_{d\geq 1}\log
(1-x_dt^d)^{-1}=\sum_{d,e\geq 1}\frac 1ex_d^et^{de}$. The coefficient
of $t^n$ in this series is $\sum_{de=n}\frac 1ex_d^e=\sum_{d\mid
n}\frac dnx_d^{n/d}=\frac 1nw_n(x)$. Thus $\log\Lambda
(x;t)=\sum_{n\geq 1}w_n(x)\frac{t^n}n$ so $\Lambda (x;t)=\exp
(\sum_{n\geq 1}w_n(x)\frac{t^n}n)$. (In the theory of Witt vectors
this formula appears in the equivalent form $t\frac d{dt}\log\Lambda
(x;t)=\sum_{n\geq 1}w_n(x)t^n$.)

\bdf For every ring $R$ we define the algebra $B'(R)$ generated by
$x=(x_1,x_2,\ldots )$, $y=(y_1,y_2,\ldots )$ with the relations
$[x_m,x_n]=[y_m,y_n]=0$ $\forall m,n$ and $\Lambda (y;t)\Lambda
(x;s)=\Lambda (x;s)\Lambda (y;t)(1-st)^{-1}$.
\edf

\blm For every $m,n\in\NN^*$ there is $c'_{m,n}\in\ZZ
[X_1,\ldots,X_{m-1}][Y_1,\ldots,Y_{n-1}]$ such that for every ring $R$
in $B'(R)$ we have $[y_n,x_m]=c'_{m,n}(x,y)$.

Also $B'(R)$ is spanned by $x^iy^j$ with $i,j\in\II$.
\elm
\pf (i) Since $B'(R)=B'(\ZZ )\otimes_\ZZ R$ it is enough to take the
case $R=\ZZ$.

Let $\Lambda (x;s)=a_0+a_1s+\cdots$ and $\Lambda
(y;t)=b_0+b_1t+\cdots$. For every $m\in\NN^*$ we have $\prod_{i\leq
m-1}(1-x_is^i)^{-1}\equiv 1\mod (s)$ and $\prod_{i\geq
m}(1-x_is^i)^{-1}\equiv (1-x_ms^m)^{-1}\equiv 1+x_ms^m\mod
(s^{m+1})$. Hence
\begin{multline*}
\sum_{k=0}^\infty a_ks^k=\Lambda (x;s)\equiv
\left(\prod_{i\leq m-1}(1-x_is^i)^{-1}\right) (1+x_ms^m)\\
\equiv \prod_{i\leq m-1}(1-x_is^i)^{-1}+x_ms^m\mod (s^{m+1}).
\end{multline*}
Therefore if $k<m$ then $a_k$ equals the coefficient of $s^k$ in
$\prod_{i\leq m-1}(1-x_is^i)^{-1}$ so $a_k\in\langle
x_1,\ldots,x_{m-1}\rangle$. Also $a_m=x_m+\alpha$, where $\alpha$ is
the the coefficient of $s^m$ in $\prod_{i\leq m-1}(1-x_is^i)^{-1}$ so
$\alpha\in\langle x_1,\ldots,x_{m-1}\rangle$. Similarly, if
$n\in\NN^*$ then for $l<n$ we have $b_l\in\langle
y_1,\ldots,y_{n-1}\rangle$ and also $b_n=y_n+\beta$ for some
$b_l\in\langle y_1,\ldots,y_{n-1}\rangle$.

The relation $\Lambda (x;s)\Lambda (y;t)=\Lambda (y;t)\Lambda
(x;s)(1-st)^{-1}$ writes
as 
$$\left(\sum_lb_lt^l\right)\left(\sum_ka_ks^k\right)
=\left(\sum_ka_ks^k\right)\left(\sum_lb_lt^l\right)
\left(\sum_rs^rt^r\right).$$ 
We identify the coefficients of $s^mt^n$ and we get
$b_na_m=\sum_{r=0}^{\min\{ m,n\}}a_{m-r}b_{n-r}$,
i.e. $[y_n+\beta,x_m+\alpha]=[b_n,a_m]=c$, where $c=\sum_{r=1}^{\min\{
m,n\}}a_{m-r}b_{n-r}$. Thus $[y_n,x_m]+[\beta,x_m]+[y_n,\alpha]+\gamma
=0$, where $\gamma =[\beta,\alpha ]-c$. But $\alpha\in\langle
x_1,\ldots,x_{m-1}\rangle$, $\beta\in\langle
y_1,\ldots,y_{n-1}\rangle$ and for $r\geq 1$ we have
$a_{m-r}\in\langle x_1,\ldots,x_{m-1}\rangle$, $b_{n-r}\in\langle
y_1,\ldots,y_{n-1}\rangle$ so $c\in\langle
x_1,\ldots,x_{m-1},y_1,\ldots,y_{n-1}\rangle$. Hence
$\alpha,\beta,\gamma\in\langle
x_1,\ldots,x_{m-1},y_1,\ldots,y_{n-1}\rangle$ and our result follows
by Lemma 2.8. \qed

\blm If $X,Y$ belong to a $\QQ$-algebra $C$ and $a:=[Y,X]$ commutes
with $X$ and $Y$ (in partcular, if $a\in\QQ$) then we have the
equality of formal series $\exp (tY)\exp (sX)=\exp (sX)\exp (tY)\exp
(ast)$.
\elm
\pf Since $[\cdot,X]$ is a derivation and $a=[Y,X]$ commutes with $Y$
we have $[X,Y^n]=\sum_{h=1}^nY^{h-1}[Y,X]Y^{n-h}=nY^{n-1}a$. By the
linearity of $[\cdot,X]$, if $f(Y)$ is a polynomial in the variable
$Y$ or, more generally, a series with coefficients in $\QQ [Y]$ then
$[f(Y),X]=\frac d{dY}f(Y)a$. When we take $f(Y)=\exp (tY)$, so $\frac
d{dY}f(Y)=t\exp (tY)$, we get $\exp (tY)X-X\exp (tY)=[\exp
(tY),X]=t\exp (tY)a$. It follows that $\exp (tY)X=(X+ta)\exp (tY)$,
i.e. $\exp (tY)X\exp (-tY)=X+ta$. Since $\alpha\mapsto\exp
(tY)\alpha\exp (-tY)$ is an automorphism of $C[[s,t]]$ we get $\exp
(tY)g(X)\exp (-tY)=g(X+ta)$ so $\exp (tY)g(X)=g(X+ta)\exp (tY)$ for
every polynomial $g(X)$ in the variable $X$ or, more generally, for
any series with coefficients in $\QQ [X]$. When we take $g(X)=\exp
(sX)$ we get $\exp (tY)\exp (sX)=\exp (sX+sta)\exp (tY)$. But $sta$
commutes with $sX$ and $tY$ so $\exp (sX+sta)\exp (tY)=\exp (sX)\exp
(tY)\exp (ast)$ and we get our result.




Alternatively, one can use the weaker Baker-Campbell-Hausdorff formula
$e^Ae^B=e^{A+B+\frac 12[A,B]}=e^{A+B}e^{\frac 12[A,B]}$, which holds
when $[A,B]$ commutes with $A$ and $B$. Together with
$e^Be^A=e^{B+A}e^{\frac 12[B,A]}=e^{A+B}e^{-\frac 12[A,B]}$, this
implies $e^Ae^B=e^Be^Ae^{[A,B]}$. Then our result follows by taking
$A=tY$ $B=sX$. Indeed, we have $[B,A]=ast$, which commutes with $A$
and $B$. \qed

\btm We have $B(\QQ )=B'(\QQ )$ and for every $m,n\in\NN^*$ we have\\
$c_{m,n}=c'_{m,n}\in\ZZ [X_d\mid d\in D^*(m)][Y_e\mid e\in D^*(n)]$.
\etm
\pf The algebras $B(\QQ )$ and $B'(\QQ )$ have the same generators. We
prove that the relations from $B'(\QQ )$ also hold in $B(\QQ )$. We
use the following obvious result. If
$\alpha_1,\ldots,\alpha_N,\beta_1,\ldots,\beta_N$ belong to some
algebra $C$ such that $\beta_n\alpha_m=\alpha_m\beta_n\gamma_{m,n}$
for some $\gamma_{m,n}\in Z(C)$ then
$\beta_1\cdots\beta_N\alpha_1\cdots\alpha_N=
\alpha_1\cdots\alpha_N\beta_1\cdots\beta_N=\prod_{m,n}\gamma_{m,n}$. We
take $C=B(\QQ )[[s,t]]$, $\alpha_m=\exp (w_m(x)\frac{s^m}m)$ and
$\beta_n=\exp (w_n(y)\frac{t^n}n)$. If $m\neq n$ then $w_m(x)$ and
$w_n(y)$ commute so $\alpha_m$ and $\beta_n$ commute. Hence we may
take $\gamma_{m,n}=1$. When $m=n$ we have $[w_n(y),w_n(x)]=n$ so
$[\frac 1nw_n(y),\frac 1nw_n(x)]=\frac 1n\in\QQ$. By Lemma 2.13 we get
$\beta_n\alpha_n=\alpha_n\beta_n\gamma_{n,n}$, where
$\gamma_{n,n}=\exp (\frac 1ns^nt^n)$. In both cases $\gamma_{m,n}\in
Z(C)$. We have
$$\alpha_1\cdots\alpha_N=\exp\left(\sum_{n=1}^Nw_n(x)\frac{s^n}n\right)
\equiv\exp\left(\sum_{n=1}^\infty w_n(x)\frac{s^n}n\right) =\Lambda
(x;s)\mod (s^{N+1}).$$
Similarly, $\beta_1\cdots\beta_N\equiv\Lambda (y;t)\mod (t^{N+1})$. We
also have
$$\prod_{m,n}\gamma_{m,n}=\exp\left(\sum_{n=1}^N\frac
1ns^nt^n\right)\equiv\exp\left(\sum_{n=1}^\infty\frac 1ns^nt^n\right)
=(1-st)^{-1}\mod (s^{N+1}t^{N+1}).$$
Therefore $\beta_1\cdots\beta_N\alpha_1\cdots\alpha_N=
\alpha_1\cdots\alpha_N\beta_1\cdots\beta_N\prod_{m,n}\gamma_{m,n}$
implies $\Lambda (y;t)\Lambda (x;s)\equiv\Lambda (x;s)\Lambda
(y;t)(1-st)^{-1}\mod (s^{N+1},t^{N+1})$. Since this holds for every
$N$ we have $\Lambda (y;t)\Lambda (x;s)=\Lambda (x;s)\Lambda
(y;t)(1-st)^{-1}$. By Corrolary 2.2 in $B(\QQ )$ we also have
$[x_m,x_n]=[y_m,y_n]=0$ so all the relations from $B'(\QQ )$ also hold
in $B(\QQ )$. But by Lemma 2.12 $x^iy^j$ with $i,j\in\II$ span
$B'(\QQ )$ and by Lemma 2.5 they are a basis for $B(\QQ )$. Then
$B(\QQ )=B'(\QQ )$ by Lemma 2.10.

In $B(\QQ )=B'(\QQ )$ we have
$[y_n,x_m]=c_{m,n}(x,y)=c'_{m,n}(x,y)$. But by Lemma 2.5 every element
in $B(\QQ )$ writes uniqely as $P(x,y)$ for some $P\in\QQ
[X][Y]$. Hence $c_{m,n}=c'_{m,n}$. Since $c'_{m,n}$ has integral
coefficients so does $c_{m,n}$ so $c_{m,n}\in\ZZ [X_d\mid d\in
D^*(m)][Y_e\mid e\in D^*(n)]$. \qed

Since $c_{m,n}$ have integral coefficents the alternative definition
of $B_{P,Q}(R)$ from Lemma 2.11 extends to arbitrary rings as follows.

\bdf For every ring $R$ and every truncation sets $P,Q$ we define
$B_{P,Q}(R)$ as the $R$ alegebra generated by $x=(x_m)_{m\in P}$ and
$y=(y_n)_{n\in Q}$, with the relations $[x_m,x_n]=0$, $[y_m,y_n]=0$ and
$[y_n,x_m]=c_{m,n}(x,y)$.
\edf

The following result generalizes Lemma 2.5, which is only for $\ZZ
[P^{-1},Q^{-1}]$-rings.

\bpr For every ring $R$ and every truncation sets $P,Q$ the products
$x^iy^j$ with $i\in\II_P$, $j\in\II_Q$ are a basis for $B_{P,Q}(R)$.
\epr
\pf Since $B_{P,Q}(R)=B_{P,Q}(\ZZ )\otimes_\ZZ R$ it is enough to
consider the case $R=\ZZ$.

As $\ZZ$-algebras, $B_{P,Q}(\ZZ )$ and $B_{P,Q}(\QQ )$ are generated by
$x$ and $y$ and by $x$, $y$ and $\QQ$, respectively. The relations
among generators in $B_{P,Q}(\ZZ )$, $[x_m,x_n]=0$, $[y_m,y_n]=0$
and $[y_n,x_m]=c_{m,n}(x,y)$, also hold in $B_{P,Q}(\QQ )$. Now for
every $m\in P$, $n\in Q$ in $B_{P,Q}(\ZZ )$ we have
$[y_n,x_m]=c_{m,n}(x,y)\in\langle x_k\mid k\in P,\, k<m,\, y_l\mid
l\in Q,\, l<n\rangle$ so, by Lemma 2.6(i), $x^iy^j$ with $i\in\II_P$,
$j\in\II_Q$ span $B_{P,Q}(\ZZ )$. By Lemma 2.5, in $B_{P,Q}(\QQ )$ they
are linearly independent over $\QQ$ and so over $\ZZ$. Then we have
$B_{P,Q}(\ZZ )\sbq B_{P,Q}(\QQ )$ by Lemma 2.10. Since $x^iy^j$ are
linearly independent (over $\ZZ$) in $B_{P,Q}(\QQ )$, they are also
linearly independent in $B_{P,Q}(\ZZ )$ so they are a basis of
$B_{P,Q}(\ZZ )$. \qed

Since $B_{P,Q}(R)$ is a free $R$-module we have:

\bco If $P,Q$ are truncation sets and $R\sbq S$ then $B_{P,Q}(R)\sbq
B_{P,Q}(S)$.
\eco

\bpr For any truncation sets $P,Q,P',Q'$ with $P'\sbq P$, $Q'\sbq Q$
and any ring $R$ we have $B_{P',Q'}(R)\sbq B_{P,Q}(R)$. Equivalently,
if $x,y$ are the generators of $B_{P,Q}(R)$ then $B_{P',Q'}(R)$ is the
subalgebra of $B_{P,Q}(R)$  generated by $x_{P'}$ an $y_{Q'}$.
\epr
\pf $B_{P,Q}(R)$ is generated by $x=(x_m)_{m\in P}$, $y=(y_n)_{n\in
Q}$ and $B_{P',Q'}(R)$ by $x':=x_{P'}=(x_m)_{m\in P'}$,
$y':=y_{Q'}=(y_n)_{n\in Q'}$. The relations among generators in
$B_{P',Q'}(R)$ also hold in $B_{P,Q}(R)$. By Lemma 2.15 $x'^iy'^j$ with
$i\in\II_{P'}$, $j\in\II_{Q'}$ are a basis in $B_{P',Q'}(R)$. They
are also linearly independent in $B_{P,Q}(R)$, where they are a part
of the basis $x^iy^j$ with $i\in\II_P$, $j\in\II_Q$. Then
$B_{P',Q'}(R)\sbq B_{P,Q}(R)$ by Lemma 2.10. \qed

\bpr $B(R)=B'(R)$ holds for every ring $R$.
\epr
\pf Both $B(R)$ and $B'(R)$ are generated by $x=(x_1,x_2,\ldots )$ and
$y=(y_1,y_2,\ldots )$. The relations among generators from $B'(R)$,
$[x_m,x_n]=[y_m,y_n]=0$ and $[y_n,x_m]=c'_{m,n}(x,y)=c_{m,n}(x,y)$,
also hold in $B(R)$. By Lemma 2.12 $x^iy^j$ with $i,j\in\II$ span
$B'(R)$ and by Lemma 2.15 they are a basis of $B(R)$. Then $B(R)=B'(R)$
by Lemma 2.10. \qed

\blm For every truncation sets $P,Q$ and any ring $R$ we have
$$B_{P,Q}(R)(x,y)^{op}=B_{Q,P}(R)(y,x)=B_{P,Q}(R)(x,-y)=B_{P,Q}(R)(-x,y).$$
Here $-x,-y$ are the opposites of $x,y$ as Witt vectors. 
\elm
\pf We first consider the case $R=\QQ$. Then $B_{P,Q}(\QQ )(x,y)$ is
the algebra generated by $x$ and $y$, where each of $x$ and $y$ has
mutually commuting entires and we have the extra relations
$[w_n(y),w_m(x)]=\delta_{m,n}p^m$. In the opposite algebra the
condition $[w_n(y),w_m(x)]=\delta_{m,n}m$ is replaced by
$[w_m(x),w_n(y)]=\delta_{m,n}m=\delta_{n,m}n$. But this is simply
the definition of $B_{Q,P}(\QQ )(y,x)$. For the second equality note
that $\QQ [x]=\QQ [-x]$ strictly, as $x\mapsto -x$ gives a
self-inverse isomorphism of $\QQ [x]$. The relations
$[w_m(x),w_n(y)]=\delta_{m,n}m$ can also be written as
$[w_n(y),w_m(-x)]=[w_n(y),-w_m(x)]=\delta_{m,n}m$. Hence we have
$B_{P,Q}(\QQ )(x,y)^{op}=B_{P,Q}(\QQ )(-x,y)$. Similarly, $B_{P,Q}(\QQ
)(x,y)^{op}=B_{P,Q}(\QQ )(x,-y)$.

Since $\ZZ\sbq\QQ$, by considering the $\ZZ$-subalgebra generated by
$x$ and $y$ in $B_{P,Q}(\QQ )(x,y)^{op}=B_{Q,P}(\QQ )(y,x)=B_{P,Q}(\QQ
)(x,-y)=B_{P,Q}(\QQ )(-x,y)$, we get $B_{P,Q}(\ZZ
)(x,y)^{op}=B_{Q,P}(\ZZ )(y,x)=B_{P,Q}(\ZZ )(x,-y)=B_{P,Q}(\ZZ
)(-x,y)$. From here, by taking the tensor product with $R$, we get
our lemma for arbitrary $R$. \qed

\blm If $R$ is a ring, $P,Q$ are truncation sets, $a=(a_m)_{m\in P}\in
W_P(R)$ and $b=(b_n)_{n\in Q}\in W_Q(R)$ then
$B_{P,Q}(R)(x,y)=B_{P,Q}(R)(x+a,y+b)$.

Here $x+a$ and $y+b$ are sums of Witt vectors.
\elm
\pf Let $x'=(x'_m)_{m\in P}$ and $y'=(y'_n)_{n\in Q}$ be
multi-variables. We take first $R$ to be a $\QQ$-ring, so we can use
Definition 1 for $B_{P,Q}(R)$. We prove that there exists an
isomorphism $f:B_{P,Q}(R)(x',y')\to B_{P,Q}(R)(x,y)$ which on
generators is given by $x'\mapsto x+a$, $y'\mapsto y+b$. To prove that
there is a morphism $f$ defined this way on generators we must show
that $f$ preserves the relations among generators. Since each of $x$
and $y$ has commmuting entries, so will $x+a$ and $y+b$. For any
$m\in P$, $n\in Q$ we have $w_m(a),w_n(b)\in R\sbq Z(B_{P,Q}(R)(x,y))$
so
$[w_n(y+b),w_m(x+a)]=[w_m(y)+w_n(b),w_m(x)+w_m(a)]
=[w_n(y),w_m(x)]=\delta_{m,n}m$. So $f$ is a morphism. By a similar
reasoning, there is a morphism $g:B_{P,Q}(R)(x,y)\to
B_{P,Q}(R)(x',y')$ given by $x\mapsto x'-a$, $y\mapsto
y'-b$. Obviously $f$ and $g$ are inverse to each other. So $f$ is an
isomorphism. Since $f$ is given by $x'\mapsto x+a$, $y'\mapsto y+b$ we
have $B_{P,Q}(R)(x,y)=B_{P,Q}(R)(x+a,y+b)$.

We now consider the multivariables $z=(z_m)_{m\in P}$,
$t=(t_n)_{n\in Q}$ and we take $R=\QQ [z,t]$ and $a=z$, $b=t$. Then we
have $B_{m,n}(\QQ [z,t])(x,y)=B_{m,n}(\QQ [z,t])(x+z,y+t)$. Since $\ZZ
[z,t]\sbq\QQ [z,t]$, by considering the $\ZZ [z,t]$-subalgebra
generated by $x,y$ we get $B_{m,n}(\ZZ [z,t])(x,y)=B_{m,n}(\ZZ
[z,t])(x+z,y+t)$. (Note that the $\ZZ [z,t]$-subalgebra generated by
$x,y$ is the same with the $\ZZ [z,t]$-subalgebra generated by
$x+z,y+t$.) 

Take now an arbitrary ring $R$ and let $a=(a_m)_{m\in P}\in W_P(R)$
and $b=(b_n)_{n\in Q}\in W_Q(R)$. On $R$ we consider the $\ZZ
[z,t]$-module structure given by the morphism $h:\ZZ [z,t]\to R$ given
by $z\mapsto a$, $t\mapsto b$. Then we have $B_{P,Q}(\ZZ
[z,t])(x,y)\otimes_{\ZZ [z,t]]}R=B_{P,Q}(\ZZ
[z,t])(x+z,y+t)\otimes_{\ZZ [z,t]}R$,
i.e. $B_{P,Q}(R)(x,y)=B_{P,Q}(R)(x+a,y+b)$. \qed

\blm If $P,Q$ are truncation sets, $R$ is a ring and $x=(x_m)_{m\in
P}$, $y=(y_n)_{n\in Q}$, $z=(z_m)_{m\in P}$ and $t=(t_n)_{n\in Q}$ are
multivariables then
$$B_{P,Q}(R)(x,y)\otimes_RB_{P,Q}(R)(z,t)=
B_{P,Q}(R)(x+z,y)\otimes_RB_{P,Q}(R)(z,t-y).$$

Here if $A,B$ are $R$-algebras we identify every $a\in A$ and $b\in B$
as the elements $a\otimes 1$ and $1\otimes b$ of $A\otimes_RB$. Also
$x+z$ and $t-y$ are sums of Witt vectors.
\elm
\pf Note that $B_{P,Q}(R)(x,y)\otimes_RB_{P,Q}(R)(z,t)$ is a free
$R$-module with the basis $x^iy^jz^kt^l$, with $i,k\in\II_P$,
$j,l\in\II_Q$. Hence, same as for $B_{P,Q}$, if $R\sbq S$ then
$B_{P,Q}(R)(x,y)\otimes_RB_{P,Q}(R)(z,t)\sbq
B_{P,Q}(S)(x,y)\otimes_SB_{P,Q}(S)(z,t)$. We will use this property
for $R=\ZZ$, $S=\QQ$.

We consider first the case when $R=\QQ$, so we can use the original
definition of $B_{P,Q}$. Then $B_{P,Q}(\QQ )(x,y)\otimes_\QQ
B_{P,Q}(\QQ )(z,t)$ is the algebra generated by $x,y,z,t$, where each
of $x,y,z,t$ has commuting entries, the entries of $x$ and $y$ commute
with those of $z$ and $t$, $[w_n(y),w_m(x)]=\delta_{m,n}m$ and
$[w_n(t),w_m(z)]=\delta_{m,n}m$. Let now $x',y',z',t'$ be
multivariables similar to $x,y,z,t$. We prove that there is a morphism
$f:B_{m,n}(\QQ )(x',y')\otimes_\QQ B_{m,n}(R)(z',t')\to B_{m,n}(\QQ
)(x,y)\otimes_\QQ B_{m,n}(\QQ )(z,t)$ with $f(x')=x+z$, $f(y')=y$,
$f(z')=z$ and $f(t')=t-y$. We have to prove that $x+z,y,z,t-y$ satisfy
the same relations as $x',y',z',t'$. The commutativity relations
required for $x+z,y,z,t-y$ follow directly from the similar
commutativity relations involving $x,y,z,t$, with the exception of the
commutativity between the entries of $x+z$ and those of $t-y$. If
$C=\langle (x+z)_m\mid m\in P\rangle$ and $C'=\langle (t-y)_n\mid n\in
Q\rangle$ then by Corollary 2.2 we have $C=\langle w_m(x+z)\mid m\in
P\rangle$ and $C'=\langle w_n(t-y)\mid n\in Q\rangle$. It follows that
the conditions $[(t-y)_n,(x+z)_m]=0$ $\forall m\in P,\, n\in Q$ are
equivalent to $[w_n(t-y),w_m(x+z)]=0$ $\forall m\in P,\, n\in
Q$. (Both are equivalent to $[\beta,\alpha]=0$ $\forall\alpha\in
C,\,\beta\in C'$.) Since the entries of $x,y$ commute with those of
$z,t$ we have $[w_n(t-y),w_m(x+z)]=
[w_n(t)-w_n(y),w_m(x)+w_m(z)]=[w_n(y),w_m(x)]-[w_n(t),w_m(z)]=
\delta_{m,n}m-\delta_{m,n}m=0$ so we are done. The remaining
relations are $[w_n(y),w_m(x+z)]=[w_n(y),w_m(x)+w_n(z)]=
[w_n(y),w_m(x)]=\delta_{m,n}m$ and
$[w_n(t-y),w_m(z)]=[w_n(t)-w_n(y),w_m(z)]=
[w_n(t),w_m(z)]=\delta_{m,n}m$. Similarly we prove that there is a
morphism $g:B_{P,Q}(\QQ )(x,y)\otimes_\QQ B_{P,Q}(\QQ )(z,t)\to
B_{P,Q}(\QQ )(x',y')\otimes_\QQ B_{P,Q}(\QQ )(z',t')$ with
$g(x)=x'-z'$, $g(y)=y'$, $g(z)=z'$ and $g(t)=t'+y'$. Obviously $f$ and
$g$ are inverse to each other.

We obviously have 
$$f(B_{P,Q}(\ZZ )(x',y')\otimes_\ZZ B_{P,Q}(\ZZ )(z',t'))\sbq
B_{P,Q}(\ZZ )(x,y)\otimes_\ZZ B_{P,Q}(\ZZ )(z,t)$$
$$g(B_{P,Q}(\ZZ )(x,y)\otimes_\ZZ B_{P,Q}(\ZZ )(z,t))\sbq B_{P,Q}(\ZZ
)(x',y')\otimes_\ZZ B_{P,Q}(\ZZ )(z',t').$$
Hence by resticting $f$ we obtain an isomorphism
$$f_\ZZ :B_{P,Q}(\ZZ )(x',y')\otimes_\ZZ B_{P,Q}(\ZZ
)(z',t')\to B_{P,Q}(\ZZ )(x,y)\otimes_\ZZ B_{P,Q}(\ZZ )(z,t).$$
If $R$ is arbitrary we take the tensor product $\otimes R$ and we
obtain an isomorphism
$$f_R:B_{P,Q}(R)(x',y')\otimes_RB_{P,Q}(R)(z',t')\to
B_{P,Q}(R)(x,y)\otimes_R B_{P,Q}(R)(z,t),$$
given by $x'\mapsto x+z$, $y'\mapsto y$, $z'\mapsto z$, $t'\mapsto
t-y$. Hence the conclusion. \qed

\blm If $P$ is a truncation set, $R$ is a ring and for $\alpha =1,2,3$
we have the multivariables $x_\alpha =(x_{\alpha,m})_{m\in P}$ and $y_\alpha
=(y_{\alpha,n})_{n\in P}$ then
\begin{multline*}
B_P(R)(x_1,y_1)\otimes B_P(R)(x_2,y_2)\otimes B_P(R)(x_3,y_3)\\
=B_P(R)(x_1,y_1-x_2x_3)\otimes_R B_P(R)(x_2,y_2-x_1x_3)\otimes_R
B_P(R)(x_3,y_3-x_1x_2).
\end{multline*}
Here we make the same conventions as in Lemma 2.21. Also $y_1-x_2x_3$,
$y_2-x_1x_3$ and $y_3-x_1x_2$ are formulas with Witt vectors.
\elm
\pf For concenience we denote
$C(R)(x_1,y_1,x_2,y_2,x_3,y_3)=B_P(R)(x_1,y_1)\otimes_R
B_P(R)(x_2,y_2)\otimes_R B_P(R)(x_3,y_3)$. Note that
$C(R)(x_1,y_1,x_2,y_2,x_3,y_3)$ is a free $R$-module with the basis
$x_1^{i_1}y_1^{j_1}x_2^{i_2}y_2^{j_2}x_3^{i_3}y_3^{j_3}$. So if $R\sbq
S$ then $C(R)(x_1,y_1,x_2,y_2,x_3,y_3)\sbq
C(S)(x_1,y_1,x_2,y_2,x_3,y_3)$. We will use this with $R=\ZZ$,
$S=\QQ$.

Take first $R=Q$ so we can use the original definition for $B_P$. Then
$C(\QQ )(x_1,y_1,x_2,y_2,x_3,y_3)$ is the $\QQ$-algebra generated by
$x_1,y_1,x_2,y_2,x_3,y_3$, where each $x_\alpha$ and each $y_\alpha$
has commuting entries, for $\alpha\neq\beta$ the entries of $x_\alpha$
and $y_\alpha$ commute with those of $x_\beta$ and $y_\beta$ and for
every $\alpha$ and every $m,n\in P$ we have $[w_n(y_\alpha
),w_m(x_\alpha )]=\delta_{m,n}m$. Let $x'_1,y'_1,x'_2,y'_2,x'_3,y'_3$
be multivariables similar to $x_1,y_1,x_2,y_2,x_3,y_3$. We prove that
there is an morphism of $\QQ$-algebras $f:C(\QQ
)(x'_1,y'_1,x'_2,y'_2,x'_3,y'_3)\to C(\QQ )(x_1,y_1,x_2,y_2,x_3,y_3)$
given by $x'_1\mapsto x_1$, $x'_2\mapsto x_2$, $x'_3\mapsto x_3$,
$y'_1\mapsto y_1-x_2x_3$, $y'_2\mapsto y_2-x_1x_3$ and $y'_3\mapsto
y_3-x_1x_2$. We must prove that
$x_1,y_1-x_2x_3,x_2,y_2-x_1x_3,x_3,y_3-x_1x_2$ satisfy the same
relations as $x'_1,y'_1,x'_2,y'_2,x'_3,y'_3$. The commutativity
conditions follow directly from the similar conditions involving
$x_1,y_1,x_2,y_2,x_3,y_3$ with exception of the commutativity amongst
the entries of $y_1-x_2x_3$, $y_2-x_1x_3$ and $y_3-x_1x_2$. It
suffices to prove for $y_1-x_2x_3$ and $y_2-x_1x_3$. If $D=\langle
(y_1-x_2x_3)_m\mid m\in P\rangle$ and $D'=\langle (y_2-x_1x_3)_n\mid
n\in P\rangle$ then by Corollary 2.2 we have $D=\langle
w_m(y_1-x_2x_3)\mid m\in P\rangle$ and $D'=\langle w_n(y_2-x_1x_3)\mid
n\in P\rangle$. So the conditions $[(y_1-x_2x_3)_m,(y_2-x_1x_3)_n]=0$
$\forall m,n\in P$ are equivalent to
$[w_m(y_1-x_2x_3),w_n(y_2-x_1x_3)]=0$ $\forall m,n\in P$. Since the
entries of $x_1,x_2,x_3$ commute with each other, the entries of $y_1$
commute with those of $y_2$ and the entries of $x_3$ commute with
those of $y_1$ and $y_2$ we have $[w_m(y_1-x_2x_3),w_n(y_2-x_1x_3)]=
[w_m(y_1)-w_m(x_2)w_m(x_3),w_n(y_2)-w_n(x_1)w_n(x_3)]=
-[w_m(y_1),w_n(x_1)w_n(x_3)]-[w_m(x_2)w_m(x_3),w_n(y_2)]=
-[w_m(y_1),w_n(x_1)]w_n(x_3)-[w_m(x_2),w_n(y_2)]w_m(x_3)=
-\delta_{n,m}w_m(y_3)+\delta_{m,n}w_m(y_3)=0$. The remaining relations
are $[w_n(y_1-x_2x_3),w_m(x_1)]=
[w_n(y_1)-w_m(x_2)w_m(x_3),w_m(x_1)]=[w_n(y_1),w_m(x_1)]=\delta_{m,n}m$
and the similar ones for $y_2-x_1x_3$ and  $x_2$ and for $y_3-x_1x_2$
and $x_3$. 

Similarly, we have a morphism $g:C(\QQ
)(x_1,y_1,x_2,y_2,x_3,y_3)\to C(\QQ )(x'_1,y'_1,x'_2,y'_2,x'_3,y'_3)$
given by $x_1\mapsto x'_1$, $x_2\mapsto x'_2$, $x_3\mapsto x'_3$,
$y_1\mapsto y'_1+x'_2x'_3$, $y_2\mapsto y'_2+x'_1x'_3$ and $y_3\mapsto
y'_3+x'_1x'_2$. Obviously $f$ and $g$ are inverse to each other.

From here on we continue like in the proof of Lemma 2.21. $f$ and $g$
will send $C(\ZZ )(x'_1,y'_1,x'_2,y'_2,x'_3,y'_3)$ and $C(\ZZ
)(x_1,y_1,x_2,y_2,x_3,y_3)$ to each other so by restricting $f$ we get
an isomorphism $f_\ZZ :C(\ZZ )(x'_1,y'_1,x'_2,y'_2,x'_3,y'_3)\to C(\ZZ
)(x_1,y_1,x_2,y_2,x_3,y_3)$. Then for $R$ arbitrary, by taking the
tensor product $\otimes_\ZZ R$, we get an isomorphism
$f_R:C(R)(x'_1,y'_1,x'_2,y'_2,x'_3,y'_3)\to
C(R)(x_1,y_1,x_2,y_2,x_3,y_3)$ given by $x'_1\mapsto x_1$,
$x'_2\mapsto x_2$, $x'_3\mapsto x_3$, $y'_1\mapsto y_1-x_2x_3$,
$y'_2\mapsto y_2-x_1x_3$ and $y'_3\mapsto y_3-x_1x_2$. Hence the
conclusion. \qed

For any $k\in\NN^*$ we denote by $F_k,V_k:W\to W$ the Frobenius and
Verschiebung maps of order $k$.

Recall that if $x=(x_n)_{n\geq 1}$ then $(V_kx)_n=x_{n/k}$ if $k\mid
n$ and $(V_kx)_n=0$ otherwise. In terms of ghost functions,
$w_n(V_kx)=kw_{n/k}(x)$ if $k\mid n$ and $w_n(V_kx)=0$ otherwise.

For $F_kx$, in terms of ghost functions we have
$w_n(F_kx)=w_{kn}(x)$. Also $(F_kx)_n\in\ZZ [x_d\mid d\in
D(kn)]$. In fact we have a more precise result,
$$(F_kx)_n\in kx_{kn}+\ZZ [x_d\mid d\in D^*(kn)].$$
Indeed, for any $e\in D^*(n)$ we have $(F_kx)_e\in\ZZ [x_d\mid d\in
D(ke)]\sbq\ZZ [x_d\mid d\in D^*(kn)]$ so $w_n(F_kx)=\sum_{e\mid
n}e(F_kx)_e^{n/e}\in n(F_kx)_n+\ZZ [x_d\mid d\in D^*(kn)]$. Also
$w_{kn}(x)=\sum_{d\mid kn}dx_d^{kn/d}\in knx_{kn}+\ZZ [x_d\mid d\in
D^*(kn)]$. Therefore $w_n(F_kx)=w_{kn}(x)$ implies that
$n(F_kx)_n$ and $knx_{kn}$ differ from each other by a polynomial
in $x_d$ with $d\in D^*(kn)$ and same happens for $(F_kx)_n$ and
$kx_{kn}$.

In the particular case when $k$ is a prime number $p$ and the base
ring has characteristic $p$ we have $(F_px)_n=x_n^p$.

\blm Let $X=(X_n)_{n\in\NN^*}$ be a multivariable regarded as a Witt
vector and let $k\geq 1$. Then for any ring $R$ the algebra morphism
$f:R[X]\to R[X]$ given on generators by $X\mapsto F_kX$ is
injective. 

Equivalently, $(F_kX)^i$ with $i\in\II$ are linearly independent.
\elm
\pf If $k_1,k_2\geq 1$ and $f_1,f_2:R[X]\to R[X]$ are the algebra
morphisms given by $X\mapsto F_{k_1}X$ and $X\mapsto F_{k_2}X$,
respectively then $f_1f_2:R[X]\to R[X]$ is given by
$X\mapsto F_{k_2}F_{k_1}X=F_{k_1k_2}X$. If $f_1,f_2$ are
injective then so is $f_1f_2$. Hence if our statement holds for
$k_1,k_2$ then it also holds for $k_1k_2$. Thus it suffices to
consider the case when $k$ is a prime number $p$.

We first prove two particular cases.

{\em Case 1.} If $\car R=p$ then $F_pX=(X_n^p)_{n\in\NN^*}$ so $f$
is given by $P(X_1,X_2,\ldots )\mapsto P(X_1^p,X_2^p,\ldots )$ for
every polynomial $P$. Obviously in this case $f$ is injective.

{\em Case 2.} If $p$ is not a divisor of zero in $R$ then let
$F_pX=Y=(Y_n)_{n\geq 1}$. Then $f(X_n)=Y_n\in pX_{pn}+\ZZ [X_d\mid
d\in D^*(pn)]$ so $Y_n=pX_{pn}+$ a sum of monomials of smaller degrees
(in the lexicographic order). We denote $Z=(Z_n)_{n\geq 1}$, with
$Z_n=X_{pn}$, and for $i=(i_1,i_2,\ldots )\in\II$ we denote by
$|i|=i_1+i_2+\cdots$. Then for any monomial $aX^i$ with $0\neq a\in R$
and $i\in\II$ we have $f(aX^i)=ap^{|i|}Z^i+$ a sum of monomials of
smaller degrees. But $p$ is not a zero divisor so $ap^{|i|}\neq
0$. Thus $\deg_Xf(aX^i)=\deg Z^i$. Also note that if $i,j\in\II$ with
$i<j$ then $\deg_XZ^i<\deg_XZ^j$. It follows that if $0\neq P\in
R[X]$, $P=\sum_{i\in\II}a_iX^i$ has $\deg_XP=i_0$ then
$\deg_Xf(P)=\deg_XZ^{i_0}$. (We have $f(P)=a_{i_0}p^{|i_0|}Z^{i_0}+$ a
sum of monomials of smaller degrees.) In particular $f(P)\neq 0$.

For the general case we denote by $f_R:R[X]\to R[X]$ the algebra
morphism given by $X\mapsto F_pX$. Then $f_R=f_\ZZ\otimes 1_R$. We
prove a more general result, namely that $f_\ZZ$ is universally
injective, i.e. that for every $\ZZ$-module $M$ the morphism of
$\ZZ$-modules $f_\ZZ\otimes 1_M:\ZZ [X]\otimes_\ZZ M\to\ZZ
[X]\otimes_\ZZ M$ is injective. It suffices to consider the case when
$M$ is finitely generated. Since a finitely generated $\ZZ$-module is
a direct sum of modules of the form $\ZZ$ or $\ZZ/q^s\ZZ$ for some
prime $q$ and $s\geq 1$, it suffices to take $M$ of this form. Since
$\ZZ$ and $\ZZ/q^s\ZZ$ are rings, we have reduced our problem to rings
of this type.

If $R=\ZZ$ or $\ZZ/q^s\ZZ$ for some prime $q\neq p$ then $p$ is not a
zero divisor in $R$ so we are in the case 2, proved above.

If $R=\ZZ/p^s\ZZ$ then $R[X]=\ZZ [X]/p^s\ZZ [X]$ and $f_R$ is
injective iff $f_\ZZ^{-1}(p^s\ZZ [X])=p^s\ZZ [X]$. If $s=1$ then
$R=\FF_p$ is of characteristic $p$ so our result holds by case
1. Hence $f_\ZZ^{-1}(p\ZZ [X])=p\ZZ [X]$. Suppose now that $s\geq 1$
is arbitrary. Assume that $f_R$ is not injective so there is $P\in
f_\ZZ^{-1}(p^s\ZZ [X])\setminus p^s\ZZ [X]$, i.e. $f_\ZZ (P)\in p^s\ZZ
[X]$, but $P\notin p^s\ZZ [X]$. Let then $t<s$ be maximal with $P\in
p^t\ZZ [X]$. Then $P=p^tQ$ with $Q\in\ZZ [X]\setminus p\ZZ [X]$. Let
also $f_\ZZ (P)=p^sT$ for some $T\in\ZZ [X]$. Hence $p^tf_\ZZ
(Q)=f_\ZZ (P)=p^sT$. Since $\ZZ [X]$ is torsion-free we get $f_\ZZ
(Q)=p^{s-t}T\in p\ZZ [X]$. Since $Q\notin p\ZZ [X]$, this contradicts
$f_\ZZ^{-1}(p\ZZ [X])=p\ZZ [X]$. Hence the conclusion.

Since $X^i$ with $i\in\II$ are a basis for $R[X]$ the injectivity of
$f$ is equivalent to the linear independence of $f(X^i)=(F_kX)^i$
with $i\in\II$. \qed

\bdf For any $k\geq 1$ and any Witt vector $x$ we denote by
$V_{k^{-1}}x=(x_{kn})_{n\geq 1}$.
\edf
The notation is justified by the fact that $V_{k^{-1}}:W\to W$ is an
inverse to the left for $V_k$.

Note that $\NN^*\setminus k\NN^*$ is a truncation set and we may write
$x=(x_{\NN^*\setminus k\NN^*},V_{k^{-1}}x)$, in the sense that
$x_{\NN^*\setminus k\NN^*}$ contains the entries $x_n$ of $x$ with
$k\nmid n$ and $V_{k^{-1}}x$ contains those with $k\mid n$. 

\blm We have $w_{kn}(x)\in kw_n(V_{k^{-1}}x)+\ZZ [x_{\NN^*\setminus k\NN^*}]$. 
\elm
\pf Since $V_{k^{-1}}x=(x_{kn})_{n\geq 1}$ we have
$kw_n(V_{k^{-1}}x)=k\sum_{e\mid n}ex_{ke}^{n/e}=\sum_{e\mid
n}kex_{ke}^{kn/ke}$. So $kw_n(V_{k^{-1}}x)$ is the sum of all terms
$dx_d^{kn/d}$ from $w_{kn}(x)=\sum_{d\mid kn}dx_d^{kn/d}$ with
$k\mid d$. It follows that $w_{kn}x=kw_n(V_{k^{-1}}x)+\sum_{d\mid
kn,k\nmid d}dx_d^{kn/d}$. But $\sum_{d\mid
kn,k\nmid d}dx_d^{kn/d}\in\ZZ [x_{\NN^*\setminus k\NN^*}]$, so we get our
claim. \qed

\bpr Let $R$ be a ring and let $k,l\geq 1$. Then in $B(R)$ we have:

(i) $\langle x,F_ky\rangle =R[x_{\NN^*\setminus
 k\NN^*}]\otimes_RB(R)(V_{k^{-1}}x,F_ky)$.

(ii) $\langle F_lx,y\rangle
=B(R)(F_lx,V_{l^{-1}}y)\otimes_RR[y_{\NN^*\setminus l\NN^*}]$.
\epr
\pf (i) Let $z=(z_m)_{m\geq 1}$, $t=(t_n)_{n\geq 1}$ be
multivariables. We claim that there is a morphism of algebras
$f=f_R:R[x_{\NN^*\setminus k\NN^*}]\otimes_RB(R)(z,t)\to B(R)$ given by
$x_{\NN^*\setminus k\NN^*}\mapsto x_{\NN^*\setminus k\NN^*}$,
$z\mapsto V_{k^{-1}}x$ and $t\mapsto F_k y$. 

We first take the case $R=\QQ$. We must
prove that the relations among the generators $x_{\NN^*\setminus
k\NN^*}$, $z$ and $t$ of $\QQ [x_{\NN^*\setminus k\NN^*}]\otimes_\QQ
B(\QQ )(z,t)$ are preserved by $x_{\NN^*\setminus
k\NN^*}$, $V_{k^{-1}}x$ and $F_ky$ in $B(\QQ )$. But these
relations are the mutual commutativity of the entries of each of
$x_{\NN^*\setminus k\NN^*}$, $z$ and $t$, the commutativity between the
entries of $x_{\NN^*\setminus k\NN^*}$ and those of $z$ and $t$, and
$[w_n(t),w_m(z)]=\delta_{m,n}m$ for $m,n\geq 1$. The corresponding
conditions for $x_{\NN^*\setminus k\NN^*}$,  $V_{k^{-1}}x$ and $F_ky$ in
$B(\QQ )$ are trivial, with the exception of the commutativity between
the entries of $x_{\NN^*\setminus k\NN^*}$ and those of $F_ky$ and
$[w_n(F_ky),w_m(V_{k^{-1}}x)]=\delta_{m,n}m$ for $m,n\geq 1$. 

If $m\in\NN^*\setminus k\NN^*$ and $n\in\NN^*$ then $m\neq kn$ so
$[w_n(F_ky),w_m(x)]=[w_{kn}(y),w_m(x)]=\delta_{m,kn}m=0$. Hence
every element of $C:=\langle w_m(x)\mid m\in\NN^*\setminus
k\NN^*\rangle$ will commute with every element of $C':=\langle
w_n(F_ky)\mid n\in\NN^*\rangle$. But $\NN^*\setminus k\NN^*$ and
$\NN^*$ are truncation sets so by Corollary 2.2 we have $C:=\langle
x_m\mid m\in\NN^*\setminus k\NN^*\rangle =\langle x_{\NN^*\setminus
k\NN*}\rangle$ and $C':=\langle (F_ky)_n\mid n\in\NN^*\rangle
=\langle F_ky\rangle$. So we have the commutativity between the
entries of $x_{\NN^*\setminus k\NN^*}$ and those of $F_ky$. 

We have $[w_n(F_ky),w_{km}(x)]=[w_{kn}(y),w_{km}(x)]=
\delta_{km,kn}km=\delta_{m,n}km$. By Lemma 2.24 we also have
$w_{km}(x)=kw_m(V_{k^{-1}}x)+\alpha$ for some $\alpha\in\langle
x_{\NN^*\setminus k\NN^*}\rangle$. But, as we have just proved, the
entries of $x_{\NN^*\setminus k\NN^*}$ commute with those of
$F_ky$. It follows that $\alpha$ commutes with
$w_n(F_ky)$. Therefore $\delta_{m,n}km=[w_n(F_ky),w_{km}(x)]=
[w_n(F_ky),kw_m(V_{k^{-1}}x)+\alpha ]=k[w_n(F_ky),w_m(V_{k^{-1}}x)]$
so $[w_n(F_ky),w_m(V_{k^{-1}}x)]=\delta_{m,n}m$. 

So we proved the existence of $f_\QQ$. Since $f_\QQ$ sends the
generators $x_{\NN^*\setminus k\NN^*}$, $z$, $t$ to elements from $B(\ZZ
)$ we have $f_\QQ (\ZZ [x_{\NN^*\setminus k\NN^*}]\otimes_\ZZ
B(\ZZ)(z,t))\sbq B(\ZZ )$. Therefore $f_\ZZ$ is simply defined as the
restriction of $f_\QQ$. Then for an arbitrary ring $R$ we obtain $f_R$
from $f_\ZZ$ by taking the tensor product $\otimes_\ZZ R$.

Next we prove that $f$ is injective. To do this we prove that the
basis $x_{\NN^*\setminus k\NN^*}^{i_1}\otimes z^{i_2}t^j$ of
$R[x_{\NN^*\setminus k\NN^*}]\otimes_RB(R)(z,t)$, with
$i_1\in\II_{\NN^*\setminus k\NN^*}$ and $i_2,j\in\II$, is mapped by $f$
to a linearly independent set. We have $f(x_{\NN^*\setminus
k\NN^*}^{i_1}\otimes z^{i_2}t^j)=x_{\NN^*\setminus
k\NN^*}^{i_1}(V_{k^{-1}}x)^{i_2}(Fy)^j=\mu (x_{\NN^*\setminus
k\NN^*}^{i_1}(V_{k^{-1}}x)^{i_2}\otimes (Fy)^j)$, where $\mu
:R[x]\otimes_RR[y]\to B(R)$ is the multiplication map,
$\alpha\otimes\beta\mapsto\alpha\beta$. Since $\mu$ is bijective we
must prove that $x_{\NN^*\setminus
k\NN^*}^{i_1}(V_{k^{-1}}x)^{i_2}\otimes (Fy)^j$ are linerly
independent in $R[x]\otimes_RR[y]$. But this will follow from the fact
that $x_{\NN^*\setminus k\NN^*}^{i_1}(V_{k^{-1}}x)^{i_2}$ with
$i_1\in\II_{\NN^*\setminus k\NN^*}$, $i_2\in\II$ are linearly
independent in $R[X]$ and $(F_ky)^j$ with $j\in\II$ are linearly
independent in $R[y]$. Indeed, we have $(x_{\NN^*\setminus
k\NN^*},V_{k^{-1}}x)=x$ so $\{ x_{\NN^*\setminus
k\NN^*}^{i_1}(V_{k^{-1}}x)^{i_2}\mid i_1\in\II_{\NN^*\setminus
k\NN^*},\, i_2\in\II\} =\{x^i\mid i\in\II\}$, which is a basis for
$R[x]$. And by Lemma 2.23 $(FX)^j$ with $j\in\II$ are linearly
independent in $R[X]$. Since $R[X]\cong R[y]$ this implies that
$(Fy)^j$ with $j\in\II$ are linearly independent in $R[y]$.

Since $f$ is an injective morphism of algebras we have ${\rm Im}\,
f=R[x_{\NN^*\setminus k\NN^*}]\otimes_RB(R)(V_{k^{-1}}x,F_ky)$. But
${\rm Im}\, f=\langle x_{\NN^*\setminus
k\NN^*},V_{k^{-1}}x,F_ky\rangle =\langle x,F_ky\rangle$, which
concludes our proof. 

(ii) proves similarly. \qed

\bco For any $k,l,m,n\in\NN^*$ and any ring $R$ in $B(R)$ we have
$$(i)\quad [(F_ky)_n,x_m]=\begin{cases}c_{m/k,n}(x_{kd}\mid d\in
D^*(m/k),\,(F_ky)_e\mid e\in D^*(n))&\text{if }k\mid m\\
0&\text{if }k\nmid m\end{cases}.$$

$$(ii)\quad [y_n,(F_lx)_m]=\begin{cases}c_{m,n/l}((F_lx)_d\mid d\in
D^*(m),\,y_{le}\mid e\in D^*(n/l))&\text{if }l\mid n\\
0&\text{if }l\nmid n\end{cases}.$$
\eco
\pf Statement (i) in the case $k\mid m$ can also be written as
$[(F_ky)_n,x_{km}]=c_{m,n}(x_{kd}\mid d\in D^*(m),\,(F_ky)_e\mid
e\in D^*(n))$. Since $V_{k^{-1}}x=(x_{km})_{m\geq 1}$ this also writes
as $[(F_ky)_n,(V_{k^{-1}}x)_m]=c_{m,n}(V_{k^{-1}}x,F_ky)$. But this
is a relation among generators from $B(R)(V_{k^{-1}}x,F_ky)$, which
exists by Proposition 2.25(i). 

Also by Proposition 2.25(i) the entries of $x_{\NN^*\setminus k\NN^*}$
commute with those of $F_ky$, i.e. every $x_m$ with $k\nmid m$
commutes with every $(F_ky)_n$. But this is just the statement (i)
in the case $k\nmid m$.

The proof of (ii) is similar. \qed

\bpr Let $R$ be a ring of characteristic $p$ and let $k,l\geq 0$. We
write the generators $x,y$ of $B(R)$ as $x=(x',x'')$ and $y=(y',y'')$,
where $x'=x_{\NN^*\setminus p^k\NN^*}$ and
$x''=V_{p^{-k}}x=(x_{p^km})_{m\geq 1}$, $y'=y_{\NN^*\setminus p^l\NN^*}$
and $y''=V_{p^{-l}}y=(y_{p^ln})_{n\geq 1}$. Then in $B(R)$ we have
$$\langle x',y',F_{p^l}x'',F_{p^k}y''\rangle =B_{\NN^*\setminus
p^k\NN^*,\NN^*\setminus
p^l\NN^*}(R)(x',y')\otimes_RB(R)(F_{p^l}x'',F_{p^k}y'').$$
\epr
\pf First note that the subalgebra of $B(R)$ generated by
$x'=x_{\NN^*\setminus p^k\NN^*}$ and $y'=y_{\NN^*\setminus p^l\NN^*}$ is
$B_{\NN^*\setminus p^k\NN^*,\NN^*\setminus p^l\NN^*}(R)$. Let
$z=(z_m)_{m\geq 1}$ and $t=(t_n)_{n\geq 1}$ be multivariables. We must
prove that there is an isomorphism of algebras
$$f:B_{\NN^*\setminus p^k\NN^*,\NN^*\setminus
p^l\NN^*}(R)\otimes_RB(R)(z,t)\to C:=\langle
x',y',F_{p^l}x'',F_{p^k}y''\rangle$$
given by $x'\mapsto x'$, $y'\mapsto y'$, $z\mapsto F_{p^l}x''$,
$t\mapsto F_{p^k}y''$. 

First we note that we have a morphism of algebras
$f_1:B_{\NN^*\setminus p^k\NN^*,\NN^*\setminus p^l\NN^*}(R)\to C$,
which is simply the inclusion map, $x'\mapsto x'$, $y'\mapsto y'$.

Since $\car R=p$ the Frobenius map $F_{p^l}$ is given by
$(x_1,x_2,\ldots )\mapsto (x_1^{p^l},x_2^{p^l},\ldots )$. Then we get
$F_{p^l}x''=F_{p^l}V_{p^{-k}}x=V_{p^{-k}}F_{p^l}x=
(x_{p^km}^{p^l})_{m\geq 1}$. Similarly,
$F_{p^k}y''=F_{p^k}V_{p^{-l}}y=V_{p^{-l}}F_{p^k}y$. By
Proposition 2.25 we have the existence of the algebras
$B(R)(V_{p^{-k}}x,F_{p^k}y)$ and
$B(R)(F_{p^l}x,V_{p^{-l}}y)$. It follows that we have the algebra
morphisms $g_1,g_2:B(R)\to B(R)$, with $g_1$ given by $x\mapsto
V_{p^{-k}}x$, $y\mapsto F_{p^k}y$ and $g_2$ by
$x\mapsto F_{p^l}x$, $y\mapsto V_{p^{-l}}y$. Then $g_1g_2:B(R)\to
B(R)$ is given by $x\mapsto V_{p^{-k}}F_{p^l}x=F_{p^l}x''$,
$y\mapsto F_{p^k}V_{p^{-l}}y=F_{p^k}y''$. By changing the
variables for $B(R)$ to $z,t$ we get a morphism of
algebras $f_2:B(R)(z,t)\to C\sbq B(R)$ given by
$z\mapsto F_{p^l}x''$, $t\mapsto F_{p^k}y''$.

We have ${\rm Im}\, f_1=\langle x',y'\rangle$ and ${\rm Im}\,
f_2=\langle F_{p^l}x'',F_{p^k}y''\rangle$. By Proposition 2.25(i)
the entries of $x'=x_{\NN^*\setminus p^k\NN^*}$ commute with those of
$F_{p^k}y$. Therefore they commute also with the entries of
$V_{p^{-l}}F_{p^k}y=F_{p^k}y''$. Similarly the entries of $y'$
commute with those of $F_{p^l}x''$. Hence the images of $f_1$ and
$f_2$ commute with each other. So we have a morphism $f=f_1\otimes
f_2:B_{\NN^*\setminus p^k\NN^*,\NN^*\setminus p^l\NN^*}(R)\otimes
B(R)(z,t)\to C$ given by $x'\mapsto x'$, $y'\mapsto y'$, $z\mapsto
F_{p^l}x''$, $t\mapsto F_{p^k}y''$. Obviously, $f$ is surjective. For
injectivity we prove that the basis $\{x'^ay'^bz^ct^d\mid
a\in\II_{\NN^*\setminus p^k\NN^*},\, b\in\II_{\NN^*\setminus p^l\NN^*},\,
c,d\in\II\}$ of $B_{\NN^*\setminus p^k\NN^*,\NN^*\setminus p^l\NN^*}(R)\otimes
B(R)(z,t)$ is mapped by $f$ to a linearly independent set. Since
$x=(x',x'')$, $y=(y',y'')$ and the entries of $y'$ commute with those
of $F_{p^l}x''$ we have
$f(x'^ay'^bz^ct^d)=x'^ay'^b(F_{p^l}x'')^c(F_{p^k}y'')^d=
x'^a(F_{p^l}x'')^cy'^b(F_{p^k}y'')^d=x'^ax''^{p^lc}y'^by''^{p^kd}=x^iy^j$,
where $i,j\in\II$, $i=(a,p^lc)$, $j=(b,p^kd)$. But $x^iy^j$ with
$i,j\in\II$ are a basis of $B(R)$ so they are linearly independent and
the map $(a,b,c,d)\mapsto ((a,p^lc),(b,p^kd))$ is injective. Hence our
claim. \qed

We will need a truncated version of Proposition 2.27. If $P$
is a truncation set, $x=(x_n)_{n\in P}$ is a truncated Witt vector of
type $W_P$ and $k\in\NN^*$ then $V_kx$ is a Witt vector of type
$W_{D(k)P}$, while $F_kx$ is of type $W_{P/k}$. (Here $D(k)P:=\{
dn\mid d\in D(k), n\in P\}$ and $P/k:=\{ n\in\NN^*\mid kn\in P\}$.)

Same as for $F_k$, $V_{k^{-1}}x$ is defined as a Witt vector of
type $W_{P/k}$. More precisely, we have $V_{k^{-1}}x=(x_{kn})_{kn\in
P}=(x_{kn})_{n\in P/k}$. Since $V_{k^{-1}}x$ contains the entries of
$x$ with indices multiples of $k$ we have $x=(x_{P\setminus
k\NN^*},V_{k^{-1}}x)$. Note that $P\setminus k\NN^*$ is a truncation
set.

Note that, while in general $F_kx$ is of type $W_{P/k}$, if the
characteristic is $p$ and $k$ is a $p$-power then we can define
$F_kx$ as a vector of type $W_P$. Namely, in characteristic $p$ we
have $F_{p^k}x=(x_n^{p^k})_{n\in P}\in W_P$.

\bco Let $R$ be a ring of characteristic $p$, let $P,Q$ be
truncation sets and let $k,l\geq 0$. We write the generators
$x=(x_m)_{m\in P}$ and $y=(y_n)_{n\in Q}$ of $B_{P,Q}(R)$ as
$x=(x',x'')$ and $y=(y',y'')$, where $x'=x_{P\setminus p^k\NN^*}$,
$x''=V_{p^{-k}}x=(x_{p^km})_{m\in P/p^k}$, $y'=y_{Q\setminus p^l\NN^*}$,
$y''=V_{p^{-l}}y=(x_{p^ln})_{n\in Q/p^l}$. Then in $B_{P,Q}(R)$ we
have $$\langle x',y',F_{p^l}x'',F_{p^k}y''\rangle =B_{P\setminus
p^k\NN^*,Q\setminus
p^l\NN^*}(R)(x',y')\otimes_RB_{P/p^k,Q/p^l}(R)(F_{p^l}x'',F_{p^k}y'').$$
\eco
\pf Let $\overline x=(x_m)_{m\geq 1}$, $\overline y=(y_n)_{n\geq 1}$
be the generators of $B(R)$. Then $x=\overline x_P$ and $y=\overline
y_Q$.

We write $\overline x=(\overline x',\overline x'')$ and $\overline
y=(\overline y',\overline y'')$, where $\overline x'=\overline
x_{\NN^*\setminus p^k\NN^*}$, $\overline x''=V_{p^{-k}}\overline x$,
$\overline y'=\overline y_{Q\setminus p^l\NN^*}$, $\overline
y''=V_{p^{-l}}\overline y$. Then by Proposition 2.27 in $B(R)$ we have
$\langle\overline x',\overline y',F_{p^l}\overline
x'',F_{p^k}\overline y''\rangle =B_{\NN^*\setminus
p^k\NN^*,\NN^*\setminus p^l\NN^*}(R)(\overline x',\overline
y')\otimes_RB(R)(F_{p^l}\overline x'',F_{p^k}\overline y'')$. 

We have $x'=(x_m)_{m\in P\setminus p^k\NN^*}$ and $\overline
x'=(x_m)_{\NN^*\setminus p^k\NN^*}$ so $x'=\overline x'_{P\setminus
  p^k\NN^*}$. We have $F_{p^l}x''=(x^{p^l}_{p^km})_{m\in P/p^k}$ and
$F_{p^l}\overline x''=(x^{p^l}_{p^km})_{m\geq 1}$ so
$F_{p^l}x''=(F_{p^l}\overline x'')_{P/p^k}$. Similarly
$y'=\overline y'_{Q\setminus p^l\NN^*}$ and
$F_{p^k}y''=(F_{p^k}\overline y'')_{Q/p^l}$. By Proposition 2.17
the subalgebra of $B_{\NN^*\setminus p^k\NN^*,\NN^*\setminus
p^l\NN^*}(R)(\overline x',\overline y')$ generated by $x'=\overline
x'_{P\setminus p^k\NN^*}$ and $y'=\overline y'_{Q\setminus p^l\NN^*}$ is
$\langle x',y'\rangle =B_{P\setminus p^k\NN^*,Q\setminus
p^l\NN^*}(R)(x',y')$ and the subalgebra of $B(R)(F_{p^l}\overline
x'',F_{p^k}\overline y'')$ generated by
$F_{p^l}x''=(F_{p^l}\overline x'')_{P/p^k}$ and
$F_{p^k}y''=(F_{p^k}\overline y'')_{Q/p^l}$ is
$\langle F_{p^l}x'',F_{p^k}y''\rangle
=B_{P/p^k,Q/p^l}(R)(F_{p^l}x'',F_{p^k}y'')$. It follows that the
subalgebra $\langle x',y',F_{p^l}x'',F_{p^k}y''\rangle$ of 
$B_{\NN^*\setminus p^k\NN^*,\NN^*\setminus p^l\NN^*}(R)(\overline x',\overline
y')\otimes_RB(R)(F_{p^l}\overline x'',F_{p^k}\overline y'')$
writes as $B_{P\setminus p^k\NN^*,Q\setminus
p^l\NN^*}(R)(x',y')\otimes_RB_{P/p^k,Q/p^l}(R)(F_{p^l}x'',F_{p^k}y'')$. \qed 

\bco If $R$ is a ring of characteristic $p$, $m,n\geq 1$ and $k,l\geq
0$ then in $B(R)$ we
have $$[y_n^{p^l},x_m^{p^l}]=\begin{cases}c_{m/p^k,n/p^l}(x_{p^kd}^{p^l}\mid
d\in D^*(m/p^k),y_{p^le}^{p^k}\mid e\in D^*(n/p^l))&\text{if }p^k\mid
m,\, p^l\mid n\\
0&\text{otherwise}\end{cases}.$$
\eco
\pf Recall that in $B(R)$ we have $[y_n,x_m]=c_{m,n}(x_d\mid d\in
D^*(m),y_e\mid e\in D^*(n))$ $\forall m,n\geq 1$. But by Proposition
2.27 $F_{p^l}x''=(x_{p^km}^{p^l})_{m\geq 1}$ and
$F_{p^k}y''=(y_{p^ln}^{p^k})_{n\geq 1}$ generate
$B(R)(F_{p^l}x'',F_{p^k}y'')$ so $\forall m,n\geq 1$ we
have
$$[y_{p^ln}^{p^k},x_{p^km}^{p^l}]=c_{m,n}(x_{p^kd}^{p^l}\mid d\in
D^*(m),y_{p^le}^{p^k}\mid e\in D^*(n)).$$
This gives the formula for $[y_n^{p^l},x_m^{p^l}]$ when $p^k\mid m$,
$p^l\mid n$.

By Proposition 2.25(i) the entries of $x_{\NN^*\setminus
p^k\NN^*}=(x_m)_{m\in\NN^*\setminus p^k\NN^*}$ commute with those
of $F_{p^k}y=(y_n^{p^k})_{n\geq 1}$. Hence for any $m,n\geq 1$ with
$p^k\nmid m$ we have $[y_n^{p^k},x_m]=0$, which implies
$[y_n^{p^k},x_m^{p^l}]=0$. Similarly, by using Proposition 2.25(ii),
if $p^l\nmid n$ then $[y_n,x_m^{p^l}]=0$, so
$[y_n^{p^k},x_m^{p^l}]=0$, . \qed

\section{The $p$-typical $B$ algebra in characteristic $p$}

From now on we fix a prime $p$ and we only consider truncation sets of
the type $P=\{ 1,p,\ldots,p^{n-1}\}$ with $n\in\NN\cup\{\infty\}$. (If
$n=0$ then $P=\emptyset$ and if $n=\infty$ then $P=\{
1,p,p^2,\ldots\}$.) We denote by $W_n(R)$ the ring of truncated
$p$-typical Witt vectors of length $n$,
$W_n(R)=W_{\{1,p\ldots,p^{n-1}\}}(R)$. It's elements will be written
as $(x_0,x_1,\ldots,x_{n-1})$ instead of
$(x_1,x_p,\ldots,x_{p^{n-1}})$. In particular, $W(R):=W_\infty (R)$ is
the ring of $p$-typical Witt vectors. If $m\geq n$ and
$x=(x_0,\ldots,x_{m-1})\in W_m(R)$ then we denote by $x_{(n)}$ its
truncation in $W_n(R)$, $x_{(n)}=(x_0,\ldots,x_{n-1})$.

The ghost function $w_{p^i}$ will be renamed $w_i$. In the new
notation if $X=(X_0,X_1,\ldots )$ then
$w_i(X)=\sum_{k=0}^ip^kX_k^{p^{i-k}}$. 

More generally, every multivariable indexed by the set $P=\{
1,p,\ldots,p^{n-1}\}$ will now be indexed by the set $\{
0,1,\ldots,n-1\}$, i.e. instead of $(X_1,X_p,\ldots,X_{p^{n-1}})$ we
write $(X_0,X_1,\ldots,X_{n-1})$. Moreover, $\II_{\{
1,p,\ldots,p^{n-1}\}}$ will be renamed $\II_n$ and it's elements will
be denoted by $(i_0,i_1,\ldots,i_{n-1})$ instead of
$(i_1,i_p,\ldots,i_{p^{n-1}})$. So if $X=(X_0,\ldots,X_{n-1})$ and
$i=(i_0,\ldots,i_{n-1})\in\II_n$ then $X^i:=X_0^{i_0}\cdots
X_{n-1}^{i_{n-1}}$. When $n=\infty$ we put $\II =\II_\infty$.

If $m,n\in\NN\cup\{\infty\}$ we denote by $B_{m,n}(R)$ the algebra
$B_{P,Q}(R)$ with $P=\{ 1,p,\ldots,p^{m-1}\}$, $Q=\{
1,p,\ldots,p^{n-1}\}$. When $m=n$ we denote $B_n(R)=B_{n,n}(R)$. We
say that $B(R):=B_\infty (R)$ is the $p$-typical $B$ algebra over $R$
and $B_{m,n}(R)$ are it's truncations. The generators of $B_{m,n}(R)$
will be renamed as $x=(x_0,\ldots,x_{m-1})$ and
$y=(y_0,\ldots,y_{n-1})$ instead of $x=(x_1,x_p,\ldots,x_{p^{m-1}})$
and $y=(y_1,y_p,\ldots,y_{p^{n-1}})$. Note that $\ZZ
[P^{-1},Q^{-1}]=\ZZ [p^{-1}]$ (or $\ZZ$, if $m,n\in\{ 0,1\}$). So if
$R$ is a $\ZZ [p^{-1}]$-ring, i.e. with $p\in R^\times$, then
$B_{m,n}(R)$ is generated by $x=(x_0,\ldots,x_{m-1})$ and
$y=(y_0,\ldots,y_{n-1})$ with the relations $[w_i(x),w_j(x)]=0$,
$[w_i(y),w_j(y)]=0$ and $[w_j(y),w_i(x)]=\delta_{i,j}p^i$.

The polynomials $c_{p^i,p^j}\in\ZZ
[X_1,\ldots,X_{p^{i-1}}][Y_1,\ldots,Y_{p^{j-1}}]$ will be renamed
$c_{i,j}$, with $c_{i,j}\in\ZZ
[X_0,\ldots,X_{i-1}][Y_0,\ldots,Y_{j-1}]$. In particular, $c_{0,0}$ is
the old $c_{1,1}$, i.e. $c_{0,0}=1$.

Therefore, for an arbitrary ring $R$, $B_{m,n}(R)$ is the algebra
generated by $x=(x_0,\ldots,x_{m-1})$ and $y=(y_0,\ldots,y_{n-1})$,
with the relations $[x_i,x_j]=0$, $[y_i,y_j]=0$ and
$[y_j,x_i]=c_{i,j}(x,y)$. If $m'\leq m$, $n'\leq n$ then Proposition
2.17 states that $B_{m',n'}(R)$ is the subalgebra of $B_{m,n}(R)$
generated by $x_{(m')}$ and $y_{(n')}$. 

If $m$ or $n=0$ then the sequence $x$ or $y$, respectively,
is empty. Hence $B_{m,0}(R)$ is generated by $x=(x_0,\ldots,x_{m-1})$,
with the relations $[x_i,x_j]=0$ for $0\leq i,j\leq m-1$,
i.e. $B_{m,-1}(R)=R[x]$ strictly. Similarly, $B_{0,n}(R)=R[y]$
strictly, where $y=(y_0,\ldots,y_{n-1})$ and $B_0(R)=R$.
\medskip

From now on all $p$-typical Witt vectors will be over rings of
characteristic $p$.
\medskip

Most of the results from \S 2, such as Lemmas 2.19 - 2.22, can be
easily translated in the new notation for $p$-typical algebras by
simply replacing $B_{P,Q}$ for some truncation sets $P,Q$ with
$B_{m,n}$ for some $m,n\in\NN\cup\{\infty\}$. Corollaries
2.28 and 2.29 need a little more attention.

The Frobenius and Verschiebung maps $F_p$ and $V_p$ will be renamed
$F$ and $V$. In this notation $F_{p^k}$ and $V_{p^k}$ write as
$F^k$ and $V^k$. Recall that $V$ is given by $(x_0,x_1,\ldots )\mapsto
(0,x_0,x_1,\ldots )$ and, since we are in characterictic $p$, $F$ is
given by $(x_0,x_1,\ldots )\mapsto (x_0^p,x_1^p,\ldots )$. 

The map $V_{p^{-1}}$ will be renamed $V^{-1}$. It is given by
$(x_0,x_1,x_2,\ldots )\mapsto (x_1,x_2,\ldots )$. More generally, the
map $V_{p^{-k}}=V_{p^{-1}}^k$ will be written as
$(V^{-1})^k=V^{-k}$ and it is given by $(x_0,x_1,\ldots )\mapsto
(x_k,x_{k+1},\ldots )$. Note that $V^{-1}$ is an inverse only to the
left for $V$. More generally, if $k,l\in\ZZ$ then $V^kV^l=V^{k+l}$
holds in all cases except when $k>0$, $l<0$. Recall that on truncated
Witt vectors $V_{k^{-1}}$ is defined as $V_{k^{-1}}:W_P\to
W_{P/k}$. If we take $P=\{ 1,p,\cdots,p^{n-1}\}$ then $P/p^k=\{
1,p,\cdots,p^{n-k-1}\}$. Hence in the new notation $V^{-k}$ is defined
as $V^{-k}:W_n\to W_{n-k}$ and is given by
$(x_0,\ldots,x_{n-1})\mapsto (x_k,\ldots,x_{n-1})$.


We are now ready to state the $p$-typical version of Corollary
2.28. We take $P=\{ 1,p,\ldots,p^{m-1}\}$ and $Q=\{
1,p,\ldots,p^{n-1}\}$ and we take with $k\leq m$, $l\leq n$. We write
the formula $x=(x',x'')$, where $x'=x_{P\setminus p^k\NN^*}$ and
$x''=V_{p^{-k}}x$. Since $P\setminus p^k\NN^*=\{ 1,p,\ldots,p^{k-1}\}$
in the new notation we have $x'=x_{(k)}$ and $x''=V^{-k}x$. Similarly,
$y=(y',y'')$, where $y'=y_{(l)}$ and $y''=V^{-l}y$. We get:

\bpr Let $m,n\in\NN\cup\{\infty\}$ and $k,l\in\NN$ with $k\leq n$,
$l\leq n$. We write the generators $x$ and $y$ of $B_{m,n}(R)$ as
$x=(x',x'')$ and $y=(y',y'')$, where $x'=x_{(k)}=(x_0,\ldots,x_{k-1})$,
$y'=y_{(l)}=(y_0,\ldots,y_{l-1})$, $x''=V^{-k}x=(x_k,\ldots,x_{m-1})$ and
$y''=V^{-l}y=(y_l,\ldots,y_{n-1})$. Then in $B_{m,n}(R)$ we have
$$\langle x',y',F^lx'',F^ky''\rangle
=B_{k,l}(R)(x',y')\otimes_RB_{m-k,n-l}(R)(F^lx'',F^ky'').$$ 
\epr

When we take in Corollary 2.28 $m=p^i$ and $n=p^j$ we get in the new
notation: 

\bpr If $i,j,k,l\in\NN$ then in $B(R)$ we have
$$[y_j^{p^k},x_i^{p^l}]=\begin{cases}
c_{i-k,j-l}(x_k^{p^l},\ldots,x_{i-1}^{p^l},y_l^{p^k},\ldots,y_{j-1}^{p^k})
&\text{if }i\geq k,\, j\geq l\\ 
0&\text{otherwise}
\end{cases}.$$
In particular, $[y_j^{p^i},x_i^{p^j}]=c_{0,0}=1$. Also if $k>i$ then
$[y_j^{p^k},x_i]=0$ and if $l>j$ then $[y_j,x_i^{p^l}]=0$.
\epr

\bco Let $C$ be a subalgebra of $B(R)$ and let $I\sb C$ be an
ideal. If $x_i^{p^l},y_j^{p^k}\in C$ then in $C/I$ we have
$[y_j^{p^k},x_i^{p^l}]=0$ iff $i<k$ or $j<l$.
\eco
\pf The "if" part follows directly from Proposition 3.2. For the "only
if" part assume that $i\geq k$, $j\geq l$ and
$[y_j^{p^k},x_i^{p^l}]=0$. Since $x_i^{p^l}$ commutes with
$y_j^{p^k}$, $x_i^{p^j}=(x_i^{p^l})^{p^{j-l}}$ will commute with
$y_j^{p^i}=(y_j^{p^k})^{p^{i-k}}$. But by Proposition 3.2 we have
$[y_j^{p^i},x_i^{p^j}]=1\neq 0$. Contradiction. \qed

\bco If $m,n\in\NN\cup\{\infty\}$, $0\leq l_0\leq\cdots\leq
l_{m-1}$ and $0\leq k_0\leq\cdots\leq k_{n-1}$ then in $B(R)$ we have 
$$\langle x_0^{p^{l_0}},\cdots,x_{m-1}^{p^{l_{m-1}}},
y_0^{p^{k_0}},\cdots,y_{n-1}^{p^{k_{n-1}}}\rangle
=R[x_0^{p^{l_0}},\cdots,x_{m-1}^{p^{l_{m-1}}}]
[y_0^{p^{k_0}},\cdots,y_{n-1}^{p^{k_{n-1}}}].$$
\eco
\pf If $i\leq m-1$, $j\leq n-1$ then by Proposition 3.2 we have
$[y_j^{p^{k_j}},x_i^{p^{l_i}}]=0$ or
$c_{i-k,j-l}(x_{k_j}^{p^{l_i}},\ldots,x_{i-1}^{p^{l_i}},
y_{l_i}^{p^{k_j}},\ldots,y_{j-1}^{p^{k_j}})$. Since if $a\leq i-1$
then $l_a\leq l_i$ and if $b\leq j-1$ then $k_b\leq k_j$, in both
cases we have $[y_j^{p^{k_j}},x_i^{p^{l_i}}]\in\langle
x_0^{p^{l_0}},\ldots,x_{i-1}^{p^{l_{i-1}}},
y_0^{p^{k_0}},\ldots,y_{j-1}^{p^{k_{j-1}}}\rangle$. Then we get our
result from Lemma 2.6(i). \qed

From now on we focus on the finitely generated case, of algebras
$B_{m,n}(R)$, where $m,n\in\NN$. Since $B_{m,n}(R)\sbq
B(R)$ the result from Proposition 3.2 will also hold in
$B_{m,n}(R)$ when $i<m$, $j<n$. 

Suppose that $C=B_{m,n}(R)/I$, where $I\sb B_{m,n}(R)$ is an
ideal. Let $0\leq m'\leq m$, $0\leq n'\leq n$ and let $0\leq
l_0\cdots\leq l_{m'-1}$, $0\leq k_0\cdots\leq k_{n'-1}$. Let $D$ be
the subalgebra of $C$ generated by
$x_0^{p^{l_0}},\cdots,x_{m'-1}^{p^{l_{m'-1}}},
y_0^{p^{k_0}},\cdots,y_{n'-1}^{p^{k_{n'-1}}}$. By Corollary 3.4 we
have $D=R[x_0^{p^{l_0}},\cdots,x_{m'-1}^{p^{l_{m'-1}}}]
[y_0^{p^{k_0}},\cdots,y_{n'-1}^{p^{k_{n'-1}}}]$.

We are interested in $C(D)$, the centralizer of $D$ in $C$. First we
determine the powers $x_i^{p^l}$ and $y_j^{p^k}$ that belong to
$C(D)$. Since $x_i^{p^l}$ commutes with
$x_0^{p^{l_0}},\cdots,x_{m'-1}^{p^{l_{m'-1}}}$, we have $x_i^{p^l}\in
C(D)$ iff $[y_j^{p^{k_j}},x_i^{p^l}]=0$ for $0\leq j\leq n'-1$. By
Corollary 4.5 this is equivalent to $l>j$ or $k_j>i$ for every $0\leq
j\leq n'-1$, i.e. $l>j$ for every $j$ with $k_j\leq i$. Since
$k_0\leq\cdots\leq k_{n'-1}$ the smallest $l$ with this property is
$l'_i=\min\{j\mid k_j>i\}$, if this minimum is defined, i.e. if
$k_{n'-1}>i$, and $l'_i=n'$ otherwise. In the particular case when
$n'=0$ we have $D=R[x_0^{p^{l_0}},\cdots,x_{m'-1}^{p^{l_{m'-1}}}]$ so
$x_0,\ldots,x_{m-1}\in C(D)$. So in this case we take $l'_0=\cdots
=l'_{m-1}=0$. 

Similarly, for every $0\leq j\leq n-1$ the smallest $k$ with the
property that $y_j^{p^k}\in C(D)$ is $k'_j=\min\{ i\mid l_i>j\}$ if
$l_{m'-1}>j$ and $k'_j=m'$ otherwise. Again, if $m'=0$ then
$k'_0=\cdots =k'_{n-1}=0$. 

Note that $0\leq l'_0\leq\cdots\leq l'_{m-1}\leq n'$ and $0\leq
k'_0\leq\cdots\leq k'_{n-1}\leq m'$. Therefore Corollary
3.4 applies. In conclusion:

\blm Let $C=B_{m,n}(R)/I$ where $I\sb B_{m,n}(R)$ is an ideal. Let $0\leq m'\leq m$, $0\leq n'\leq n$, $0\leq l_0\leq\cdots\leq l_{m'-1}$ and $0\leq
k_0\cdots\leq k_{n'-1}$. 

We consider the subalgebra $D\sbq C$,
$$D=\langle x_0^{p^{l_0}},\ldots,x_{m'-1}^{p^{l_{m'-1}}},
y_0^{p^{k_0}},\ldots,y_{n'-1}^{p^{k_{n'-1}}}\rangle
=R[x_0^{p^{l_0}},\ldots,x_{m'-1}^{p^{l_{m'-1}}}]
[y_0^{p^{k_0}},\ldots,y_{n'-1}^{p^{k_{n'-1}}}].$$
For $0\leq i\leq m-1$, $0\leq j\leq n-1$  we define
$$l_i=\begin{cases}\min\{j\mid k_j>i\}&\text{if }k_{n'-1}>i\\
n'&\text{if }k_{n'-1}\leq i\end{cases}\quad\quad
k_j=\begin{cases}\min\{i\mid l_i>j\}&\text{if }l_{m'-1}>i\\
m'&\text{if }l_{n'-1}\leq j\end{cases}.$$

(If $n'=0$ then $l'_0=\cdots =l'_{m-1}=0$; if $m'=0$ then $k'_0=\cdots
=k'_{n-1}=0$.)

Then $0\leq l'_0\leq\cdots\leq l'_{m-1}\leq n'$, $0\leq
k'_0\leq\cdots\leq k'_{n-1}\leq m'$ and $C(D)\spq D'$, where
$D'=\langle x_0^{p^{l'_0}},\ldots,x_{m-1}^{p^{l'_{m-1}}},
y_0^{p^{k'_0}},\ldots,y_{n-1}^{p^{k'_{n-1}}}\rangle
=R[x_0^{p^{l'_0}},\ldots,x_{m-1}^{p^{l'_{m-1}}}]
[y_0^{p^{k'_0}},\ldots,y_{n-1}^{p^{k'_{n-1}}}]$. \elm

Note that $C(D)$ also writes as
$C(x_0^{p^{l_0}},\ldots,x_{m'-1}^{p^{l_{m'-1}}},
y_0^{p^{k_0}},\ldots,y_{n'-1}^{p^{k_{n'-1}}})$.

In the particular case $m'=0$, when
$D=R[y_0^{p^{k_0}},\ldots,y_{n'-1}^{p^{k_{n'-1}}}]$, we have
$k'_0=\cdots =k'_{n-1}=0$ so $D'=\langle
x_0^{p^{l'_0}},\ldots,x_{m-1}^{p^{l'_{m-1}}},
y_0,\ldots,y_{n-1}\rangle
=R[x_0^{p^{l'_0}},\ldots,x_{m-1}^{p^{l'_{m-1}}}][y]$.

Similarly, if $n'=0$, then
$D=R[x_0^{p^{l_0}},\ldots,x_{m'-1}^{p^{l_{m'-1}}}]$ and $l'_0=\cdots
=l'_{n-1}=0$, so $D'=\langle x_0,\ldots,x_{m-1},
y_0^{p^{k'_0}},\ldots,y_{n-1}^{p^{k'_{n-1}}}\rangle
=R[x][y_0^{p^{k'_0}},\ldots,y_{n-1}^{p^{k'_{n-1}}}]$.

We will prove that in Lemma 3.5 we have equality, i.e. $C(D)=D'$. But
first we need a preliminary result.

\blm If $[a,b]=1$ then $\ad (a)^k(b^n)=\frac{n!}{(n-k)!}b^{n-k}$ if
$k\leq n$ and $=0$ if $k>n$.

Similarly, $(-\ad (b))^k(a^n)=\frac{n!}{(n-k)!}a^{n-k}$ if $k\leq n$
and $=0$ if $k>n$.

In particular, $\ad (a)^n(b^n)=(-\ad (b))^n(a^n)=n!$.
\elm
\pf We have $\ad (a)(b^n)=[a,b^n]=nb^{n-1}$. By induction, if $k\leq
n$ then  $\ad (a)^k(b^n)=n(n-1)\cdots
(n-k+1)b^{n-k}=\frac{n!}{(n-k)!}b^{n-k}$. When $n=k$ we get $\ad
(a)^n(b^n)=n!$. It follows that $\ad (a)^{n+1}(b^n)=\ad
(a)(n!)=[a,n!]=0$. it follows that $\ad (a)^k(b^n)=0$ if $k\geq n+1$.

Since $-\ad (b)$ is given by $x\mapsto [x,b]$ we get the similar
results for $(-\ad (b))^k(a^n)$. \qed

\blm We use the notations from Lemma 3.5. If $M\sbq R[x]$, $M'\sbq
R(y]$ are $R$-submodules and $\mu :R[x]\otimes_RR[y]\to
B_{m,n}(R)$ is the multiplication map
$\alpha\otimes\beta\mapsto\alpha\beta$ then in $C$ we have: 

(i) $C(y_0^{p^{k_0}},\ldots,y_{n'-1}^{p^{k_{n'-1}}})\cap\mu
(R[x]\otimes_RM)=\mu
(R[x_0^{p^{l'_0}},\ldots,x_{m-1}^{p^{l'_{m-1}}}]\otimes_RM)$

(ii) $C(x_0^{p^{l_0}},\ldots,x_{m'-1}^{p^{l_{m'-1}}})\cap\mu 
(M'\otimes_RR[y])=\mu
(M'\otimes_RR[y_0^{p^{k'_0}},\ldots,y_{n-1}^{p^{k'_{n-1}}}])$.
\elm
\pf (i) By the case $m'=0$ of Lemma 3.5 we have $\mu
(R[x_0^{p^{l'_0}},\ldots,x_{m-1}^{p^{l'_{m-1}}}]\otimes_RM)\sbq
R[x_0^{p^{l'_0}},\ldots,x_{m-1}^{p^{l'_{m-1}}}][y]\sbq
C(y_0^{p^{k_0}},\ldots,y_{n'-1}^{p^{k_{n'-1}}})$. So we have the
$\spq$ inclusion.

For the reverse inclusion let $\alpha\in
C(y_0^{p^{k_0}},\ldots,y_{n'-1}^{p^{k_{n'-1}}})\cap\mu
(R[x]\otimes_RM)$. Assume that $\alpha\notin\mu
(R[x_0^{p^{l'_0}},\ldots,x_{m-1}^{p^{l'_{m-1}}}]\otimes_RM)$. 

We have the filtration
$R[x_0^{p^{l'_0}},\ldots,x_{m-1}^{p^{l'_{m-1}}}]=C_0\sbq\cdots\sbq
C_m=R[X]$, where $C_i=R[x_0,\ldots
x_{i-1},x_i^{p^{l'_i}},\ldots,x_{m-1}^{p^{l'_{m-1}}}]$. Then
$\alpha\in\mu (C_m\otimes_RM)\setminus\mu (C_0\otimes_RM)$ so there is
$0\leq i\leq m-1$ such that $\alpha\in\mu
(C_{i+1}\otimes_RM)\setminus\mu (C_i\otimes_RM)$. 

Next, we have a filtartion $C_{i+1}=C_{i,0}\spq\cdots\spq
C_{i,l'_i}=C_i$, where $C_{i,l}=R[x_0,\ldots,x_{i-1},
x_i^{p^l},x_{i+1}^{p^{l'_{i+1}}},\ldots,x_{m-1}^{p^{l'_{m-1}}}]$. Since
$\alpha\in\mu (C_{i,0}\otimes_RM)\setminus\mu (C_{i,l'_i}\otimes_RM)$
we have $\alpha\in\mu (C_{i,l}\otimes_RM)\setminus\mu
(C_{i,l+1}\otimes_RM)$ for some $0\leq l<l'_i$. 

Note that every power of $x_i^{p^l}$ writes as a power of
$x_i^{p^{l+1}}$ multiplied by one of the factors
$1,x_i^{p^l},\ldots,x_i^{p^l(p-1)}$. (If $a=bp+r$ with $0\leq r\leq
p-1$ then $(x_i^{p^l})^a=x_i^{p^lr}(x_i^{p^{l+1}})^b$.) It follows
that $C_{i,l}=\sum_{r=0}^{p-1}x_i^{p^lr}C_{i,l+1}$. Hence
$\alpha\in\mu
(C_{i,l}\otimes_RM)=\sum_{r=0}^{p-1}x_i^{p^lr}\mu
(C_{i,l+1}\otimes_RM)$. We write
$\alpha=\sum_{r=0}^{p-1}x_i^{p^lr}\alpha_r$, with $\alpha_r\in\mu
(C_{i,l+1}\otimes_RM)$. Since $\alpha\notin\mu (C_{i,l+1}\otimes_RM)$
we cannot have $\alpha =\alpha_0$ so $\alpha_1,\ldots,a_{p-1}$ are not
all zero. Hence
$\alpha=\alpha_0+\cdots +x_i^{p^lu}\alpha_u$ for some $1\leq u\leq
p-1$, with $\alpha_u\neq 0$. 

Since $l<l'_i=\min\{ j\mid k_j>i\}$ we have $k_l\leq i$. (Same happens
if $k_{n'-1}\leq i$, when $l'_i=n'$. In this case $k_l\leq
k_{n'-1}\leq i$.) Then, since $y_l^{p^{k_l}}$ commutes with $\alpha$,
so will $y_l^{p^i}$. For $h<i$ we have $[y_l^{p^i},x_h]=0$; also
$[y_l^{p^i},x_i^{p^{l+1}}]=0$; and for $h>i$ we have $l'_h\geq l'_i>l$
so $[y_l^{p^i},x_h^{p^{l'_h}}]=0$. So $y_l^{p^i}$ commutes with every
element in $C_{i,l+1}$. Since $y_l^{p^i}$ also commutes with every
element of $M\sbq R[y]$, it will commute with the elements of $\mu
(C_{i,l+1}\otimes_RM)$, in particular, with
$\alpha_0,\ldots,\alpha_u$. 

We have $\ad (y_l^{p^i})(\alpha )=[y_l^{p^i},\alpha]=0$ so $(\ad
(y_l^{p^i}))^u(\alpha )=0$. But
$\alpha=\sum_{r=0}^ux_i^{p^kr}\alpha_u$ and $y_l^{p^i}$ commutes with
$\alpha_r$ $\forall r$ so 
$$0=(\ad (y_l^{p^i}))^u(\alpha )=\sum_{r=0}^u(\ad
(y_l^{p^i}))^u(x_i^{p^ku})\alpha_r=u!\alpha_u.$$
(We have $[y_k^{p^i},x_i^{p^k}]=1$ so, by Lemma 3.6, $(\ad
(y_k^{p^i}))^u(x_i^{p^kr})=0$ if $r<u$ and $(\ad
(y_k^{p^i}))^u(x_i^{p^ku})=u!$.) But $u<p$ so $u!\in R^\times$. Hence
$\alpha_u=0$. Contradiction.

(ii) is similar. \qed

\bpr In Lemma 3.5 we have equality, $C(D)=D'$.
\epr
\pf We apply Lemma 3.7(i) with $M=R[y]$, when $\mu
(R[x]\otimes_RM)=R[x][y]=C$, and we get
$C(y_0^{p^{k_0}},\ldots,y_{n'-1}^{p^{k_{n'-1}}})=\mu
(R[x_0^{p^{l'_0}},\ldots,x_{m-1}^{p^{l'_{m-1}}}]\otimes_RR[y])$. 

Next we apply Lemma 3.7(ii) with
$M'=R[x_0^{p^{l'_0}},\ldots,x_{m-1}^{p^{l'_{m-1}}}]$. We have
$C(D)=C(x_0^{p^{l_0}},\ldots,x_{m'-1}^{p^{l_{m'-1}}})\cap
C(y_0^{p^{k_0}},\ldots,y_{n'-1}^{p^{k_{n'-1}}})=
C(x_0^{p^{l_0}},\ldots,x_{m'-1}^{p^{l_{m'-1}}})\cap\mu (M'\otimes
R[y])=\mu
(M'\otimes_RR[y_0^{p^{k'_0}},\ldots,y_{n-1}^{p^{k'_{n-1}}}])$. But
$M'=R[x_0^{p^{l'_0}},\ldots,x_{m-1}^{p^{l'_{m-1}}}]$ so $\mu
(M'\otimes_RR[y_0^{p^{k'_0}},\ldots,y_{n-1}^{p^{k'_{n-1}}}])=
R[x_0^{p^{l'_0}},\ldots,x_{m-1}^{p^{l'_{m-1}}}]
[y_0^{p^{k'_0}},\ldots,y_{n-1}^{p^{k'_{n-1}}}]$. \qed

\bco If $C=B_{m,n}(R)/I$ for some ideal $I\sb B_{m,n}(R)$ and
$m'\leq m$, $n'\leq n$ then in $C$ we have: 

(i) $C(y_0,\ldots,y_{n'-1})=R[x_0^{p^{n'}},\ldots,x_{m-1}^{p^{n'}}][y]$.

(ii) $C(x_0,\ldots,x_{m'-1})=R[x][y_0^{p^{m'}},\ldots,y_{n-1}^{p^{m'}}]$.

(iii) $C(x_0,\ldots,x_{m'-1},y_0,\ldots,y_{n'-1})=
R[x_0^{p^{n'}},\ldots,x_{m-1}^{p^{n'}}][y_0^{p^{m'}},\ldots,y_{n-1}^{p^{m'}}]$.

In particular,
$Z(C)=R[x_0^{p^{n}},\ldots,x_{m-1}^{p^{n}},
y_0^{p^{m}},\ldots,y_{n-1}^{p^{m}}]$.
\eco
\pf (iii) follows from Proposition 3.8 with
$(l_0,\ldots,l_{m'-1})=(0,\ldots,0)$ and
$(k_0,\ldots,k_{n'-1})=(0,\ldots,0)$ so that
$(k'_0,\ldots,k'_{n-1})=(m',\ldots,m')$ and
$(l'_0,\ldots,l'_{n-1})=(n',\ldots,n')$. (i) and (ii) are particular
cases, $m'=0$ and $n'=0$, respectively, of (iii). (See also the
remarks following Lemma 3.5.)

The formula for $Z(C)$ follows by taking $m'=m$ and $n'=n$ in
(iii). \qed

\section{The symbols $((\cdot,\cdot))_{p^m,p^n}$}

Throughout this section $K$ is a field of characteristic $p$ and
$m,n\in\NN$ are fixed. 

\bdf If $a=(a_0,\ldots,a_{m-1})\in W_m(K)$, $b=(b_0,\ldots,b_{n-1})\in
W_n(K)$ we define the algebra
$A_{((a,b))_{p^m,p^n}}=B_{m,n}(K)/(F^nx-a,F^my-b)$. 
\edf

Since in this section $m,n$ are fixed, for convenience we will write
$A_{((a,b))}$ instead of $A_{((a,b))_{p^m,p^n}}$. If the field $K$
needs to be specified we use the notation $A_{((a,b))}(K)$. 

If instead of $x,y$ we use other multivariables, say, $z,t$, then we
use the notation $A_{((a,b))}(z,t)$, which means
$B_{m,n}(K)(z,t)/(F^nz-a,F^mt-b)$.

Note that the relation $F^nx-a=0$ is equivalent to $F^nx=a$,
i.e. to $x_i^{p^n}=a_i$ for $0\leq i\leq m-1$. Hence
$(F^nx-a)=(x_0^{p^n}-a_0,\ldots,x_{m-1}^{p^n}-a_{m-1})$. Similarly,
$(F^my-b)=(y_0^{p^m}-b_0,\ldots,y_{n-1}^{p^m}-b_{n-1})$.

As a consequence,
$A_{((a,b)}=B_{m,n}(K)/(x_i^{p^n}-a_i,y_j^{p^m}-b_j\mid i<m,j<n)$. In
terms of generators and relations $A_{((a,b))}$ is generated by $x$
and $y$ with the relations $[x_i,x_j]=0$, $[y_i,y_j]=0$,
$[y_j,x_i]=c_{i,j}(x_0,\ldots,x_{i-1},y_0,\ldots,y_{j-1})$,
$x_i^{p^n}=a_i$ and $y_j^{p^m}=b_j$.

In particular, if $n=0$ then $B_{m,0}(K)=K[x]$ so
$A_{((a,b))}=K[x]/(x_0-a_0,\ldots,x_{m-1}-a_{m-1})=K$. Similarly, if
$m=0$ then $A_{((a,b))}=K[y]/(y_0-b_0,\ldots,y_{n-1}-b_{n-1})=K$.

\blm $A_{((a,b))}$ is central.
\elm 
\pf We apply Corollary 3.9 to $C=A_{((a,b))}$, which is a
quotient of $B_{m,n}(K)$. We have
$Z(A_{((a,b))})=K[x_0^{p^n},\ldots,x_{m-1}^{p^n},
y_0^{p^m},\ldots,y_{n-1}^{p^m}]$. But in $A_{((a,b))}$ we have
$x_i^{p^n}=a_i\in K$ and $y_j^{p^m}=b_j\in K$. So
$Z(A_{((a,b))})=K$. \qed 

\blm The multiplication map $\mu :K[x]\otimes_KK[y]\to
B_{m,n}(K)$, $\alpha\otimes\beta\mapsto\alpha\beta$, induces an
isomorphism $\overline\mu :K[x]/I_x\otimes_KK[y]/I_y\to A_{((a,b))}$,
where $I_x$ is the ideal
$(F^nx-a)=(x_0^{p^n}-a_0,\ldots,x_{m-1}^{p^n}-a_{m-1})$ of $K[x]$ and
$I_y$ is the ideal
$(F^my-b)=(y_0^{p^m}-b_0,\ldots,y_{n-1}^{p^m}-b_{n-1})$ of $K[y]$.
\elm
\pf Recall that $\mu$ is a $K$-linear isomorphism. We want to
identify the preimage under $\mu$ of the ideal
$(F^nx-a,F^my-b)$ from the definition of $A_{((a,b))}$. 

By Corollary 3.9 $x_0^{p^n}-a_0,\ldots,x_{m-1}^{p^n}-a_{m-1}\in
Z(B_{m,n}(K))$ so the ideal generated by them coincides with the
left ideal they generate. Since $B_{m,n}(K)$ is spanned by
products $PQ$, with $P\in K[x]$, $Q\in K[y]$, the ideal
$(F^nx-a)=(x_0^{p^n}-a_0,\ldots,x_{m-1}^{p^n}-a_{m-1})$ will be
spanned by $(x_i^{p^n}-a_i)PQ=\mu ((x_i^{p^n}-a_i)P\otimes Q)$, with
$P\in K[x]$, $Q\in K[y]$ and $0\leq i\leq m-1$. But the products
$(x_i^{p^n}-a_i)P$ with $0\leq i\leq m-1$ and $P\in K[x]$ span the
ideal $I_x$ of $K[x]$. Hence $(F^nx-a)$ is spanned by $PQ=\mu
(P\otimes Q)$ with $P\in I_x$, $Q\in K[y]$. Thus it is equal to $\mu
(I_x\otimes_KK[y])$. 

Similarly, since $y_0^{p^m}-b_0,\ldots,y_{n-1}^{p^m}-b_{n-1}\in
Z(B_{m,n}(K))$ the ideal
$(F^my-b)=(y_0^{p^m}-b_0,\ldots,y_{n-1}^{p^m}-b_{n-1})$ coincides with
the right ideal generated by
$y_0^{p^m}-b_0,\ldots,y_{n-1}^{p^m}-b_{n-1}$ in $B_{m,n}(K)$. Then
$(F^my-b)$ is spanned by the products $PQ(y_j^{p^m}-b_j)=\mu (P\otimes
Q(y_j^{p^m}-b_j))$, where $P\in K[x]$, $Q\in K[y]$ and $0\leq j\leq
n-1$. We get $(F^m(y)-b)=\mu (K[x]\otimes_K I_y)$.

It follows that $(F^nx-a,F^my-b)=\mu
(I_x\otimes_KK[y]+K[x]\otimes_KI_y)$. Hence $\mu$ induces an
isomorphism $\overline\mu$ between
$K[x]\otimes_KK[y]/(I_x\otimes_KK[y]+K[x]\otimes_KI_y)$ and
$B_{m,n}(K)/(F^nx-a,F^my-b)=A_{((a,b))}$. But
$K[x]\otimes_KK[y]/(I_x\otimes_KK[y]+K[x]\otimes_KI_y)=
K[x]/I_x\otimes_KK[y]/I_y$. \qed 

\bco We have $\dim_KA_{((a,b))}=p^{2mn}$ and $x_0^{i_0}\cdots
x_{m-1}^{i_{m-1}}y_0^{j_0}\cdots y_{n-1}^{j_{n-1}}$ with $0\leq
i_q\leq p^n-1$, $0\leq j_r\leq p^m-1$ are a basis. 
\eco
\pf The products $x_0^{i_0}\cdots x_{m-1}^{i_{m-1}}$ with $0\leq
i_q\leq p^n-1$ are a basis for $K[x]/I_x$ and the products
$y_0^{j_0}\cdots y_{n-1}^{j_{n-1}}$ with $0\leq j_r\leq p^m-1$ are a
basis for $K[y]/I_y$. It follows that $\overline\mu (x_0^{i_0}\cdots
x_{m-1}^{i_{m-1}}\otimes y_0^{j_0}\cdots y_{n-1}^{j_{n-1}})=x_0^{i_0}\cdots
x_{m-1}^{i_{m-1}}y_0^{j_0}\cdots y_{n-1}^{j_{n-1}}$ with $0\leq
i_q\leq p^n-1$, $0\leq j_r\leq p^m-1$ are a basis for
$A_{((a,b))}$. Since for every $i_q$ there are $p^n$ possible values
and for every $j_r$ there are $p^m$ possible values this basis has
$(p^n)^m(p^m)^n=p^{2mn}$ elements. \qed

\blm $A_{((a,b))}$ is simple.
\elm
\pf Let $0\neq I$ be an ideal of $A_{((a,b))}$.  We take $\alpha\in I$,
$\alpha\neq 0$, arbitrary and we prove that after a succesion of
transformations $\alpha\to [y_l^{p^k},\alpha ]$ we end up with an element
$\alpha\in I\cap K[y]$, $\alpha\neq 0$. Then we prove that after a
succesion of transformations $\alpha\to [\alpha,x_k^{p^l}]$ we end up
with an element $\alpha\in I\cap K$, $\alpha\neq 0$. Since $\alpha\in
I$ is invertible we have $I=A_{((a,b))}$. Hence $A_{((a,b))}$ is simple.

In the proof we use the particular cases of Proposition 3.2.

There are two steps.

Step 1: We prove that $J:=I\cap K[y]\neq\{ 0\}$.

By Corollary 4.3 every $\alpha\in I\setminus\{ 0\}$ writes uniquely in
the "standard form" as $\alpha =\sum_{i\in A}x^i\alpha_i$, where
$\emptyset\neq A\sbq\{ 0,\ldots,p^n-1\}^m$ and for every $i\in A$
$\alpha_i\in K[y]\setminus\{ 0\}$, with $\deg_{y_j}\alpha_i\leq p^m-1$
for $0\leq j\leq n-1$. We have $\alpha\in K[y]\setminus\{ 0\}$ iff
$A=\{0\}$, i.e. iff $\alpha =\alpha_0$. If $\alpha\notin K[y]$,
i.e. if $A\neq\emptyset,\{ 0\}$, then we define
$(k_\alpha,l_\alpha,q_\alpha )$ as follows: $k_\alpha$ is the largest
$k$ such that there is $i=(i_0,\ldots,i_{m-1})\in A$ with $i_k\neq 0$;
$l_\alpha$ is the smallest $l$ such that there is
$i=(i_0,\ldots,i_{m-1})\in A$ with $p^l\| i_{k_\alpha}$; $q_\alpha$ is
the largest $q$ with $p\nmid q$ such that there is
$i=(i_0,\ldots,i_{m-1})\in A$ with
$i_{k_\alpha}=p^{l_\alpha}q$. Obviously $k_\alpha,l_\alpha,q_\alpha$
are well defined and we have $0\leq k_\alpha\leq m-1$, $0\leq
l_\alpha\leq n-1$ and $1\leq q_\alpha\leq p^{n-l_\alpha}-1$ with
$p\nmid q_\alpha$. (Recall, if $i=(i_0,\ldots,i_{m-1})\in A$ then
$i_k<p^n$ $\forall k$.)

On triplets we define the order relation $\leq$, with $(k',l',q')\leq
(k,l,q)$ if $(k',-l',q')\leq (k,-l,q)$ in the lexicographic order. We
have $(k',l',q')<(k,l,q)$ if $k'<k$ or if $k'=k$ and $l'>l$ or if
$k'=k$, $l'=l$ and $q'<q$. 

Let $\alpha\in I\setminus\{ 0\}$. If $\alpha\in K[y]$ then we are
done. Otherwise let $(k,l,q)=(k_\alpha,l_\alpha,q_\alpha)$ and let
$\alpha'=[y_l^{p^k},\alpha ]\in I$. We prove that $\alpha'\in
I\setminus\{ 0\}$ and we have either $\alpha'\in K[y]$ or
$(k',l',q'):=(k_{\alpha'},l_{\alpha'},q_{\alpha'})<(k,l,q)$. If
$\alpha'\notin K[y]$ then we repeat the procedure and we define
$\alpha''=[y_{l'}^{p^{k'}},\alpha']$ and so on. At each step the
triplet $(k_\alpha,l_\alpha,q_\alpha)$ decreases. But
$(k_\alpha,l_\alpha,q_\alpha)$ belongs to a finite set so this proces
cannot go indefinitely. Eventually we get an element of $I\setminus\{
0\}$ belonging to $K[y]$.

We write $\alpha$ as above, $\alpha =\sum_{i\in A}x^i\alpha_i$. Then
$[y_l^{p^k},\alpha_i]=0$ so $\alpha'=\sum_{i\in
A}[y_l^{p^k},x^i]\alpha_i$. Let $A'=\{ i=(i_0,\ldots,i_{m-1})\in A\mid
p^l\| i_k\}$. By the definition of $l=l_\alpha$, $A'\neq\emptyset$. If
$i=(i_0,\ldots,i_{m-1})\in A'$ then $i_k=p^lq_i$ for some $q_i$ not
divisible by $p$. By definition $q=q_\alpha =\max\{ q_i\mid i\in
A'\}$. By the construction of $k=k_\alpha$ for every
$i=(i_0,\ldots,i_{m-1})\in A$ we have $i_h=0$ for $h>k$ so
$x^i=x_0^{i_0}\cdots x_k^{i_k}$. If $h<k$ then $[y_l^{p^k},x_h]=0$ so
$[y_l^{p^k},x^i]=x_0^{i_0}\cdots
x_{k-1}^{i_{k-1}}[y_l^{p^k},x_k^{i_k}]$. But by the construction of
$l=l_\alpha$ we have either $p^l\| i_k$, when $i\in A'$, or
$p^{l+1}\mid i_k$, when $i\notin A'$. If $i\notin A'$ then
$[y_l^{p^k},x_k^{p^{l+1}}]=0$ so $p^{l+1}\mid i_k$ 
implies $[y_l^{p^k},x_k^{i_k}]=0$ so $[y_l^{p^k},x^i]=0$. If $i\in A'$
then $[y_l^{p^k},x_k^{p^l}]=1$ so
$[y_l^{p^k},x_k^{i_k}]=[y_l^{p^k},x_k^{p^lq_i}]=
q_ix_k^{p^l(q_i-1)}=q_ix_k^{i_k-p^l}$. Hence
$[y_l^{p^k},x^i]=q_ix_0^{i_0}\cdots
x_{k-1}^{i_{k-1}}x_k^{i_k-p^l}=q_ix^{i-p^l{\rm e}_k}$. (Here
${\rm e}_0,\ldots,{\rm e}_{m-1}$ is the cannonical base of $\ZZ^m$.)
In conclusion, $\alpha'=\sum_{i\in A'}x^{i-p^l{\rm
e}_k}q_i\alpha_i$. For every $i\in A'$ we have $p\nmid q_i$ so
$\alpha_i\neq 0$ implies $q_i\alpha_i\neq 0$. Hence $\alpha'\in
I\setminus\{ 0\}$ and the set $A$ of indices corresponding to
$\alpha'$ is $A'-p^k{\rm e}_l=\{ i-p^l{\rm   e}_k\mid i\in
A'\}\neq\emptyset$ so $\alpha'\neq 0$. If $\alpha'\notin K[y]$ then
let $(k',l',q'):=(k_{\alpha'},l_{\alpha'},q_{\alpha'})$. Now for every
$i\in A'$ the entries of $i-p^l{\rm e}_k$ on the positions
$k+1,\ldots,m-1$ are $0$ so $k'\leq k$. The $k$th entry of
$i-p^l{\rm e}_k$ is $i_k-p^l=p^l(q_i-1)$. If $q=1$ then $p^l(q_i-1)=0$
$\forall i\in A'$ so $k'<k$, so $(k',l',q')<(k,l,q)$. Suppose that
$q>1$ so $p^l(q_i-1)\neq 0$ for some $i\in A'$. Then $k'=k$. If $p\mid
q_i-1$, so $p^{l+1}\mid p^l(q_i-1)$, $\forall i\in A'$ then $l'>l$ and
again $(k',l',q')<(k,l,q)$. Finally, if $p\nmid q_i-1$ so $p^l\|
p^l(q_i-1)$ for some $i\in A'$ then $l'=l$, but $q'=\max\{ q_i-1\mid
i\in A',\, p\nmid q_i-1\} \leq q-1<q$ and so $(k',l',q')<(k,l,q)$.

Step 2: We prove that $I\cap K=J\cap K\neq\{ 0\}$. 

Same as in Step 1, every $\alpha\in J\setminus\{ 0\}$ writes uniquely
in the "standard form" as $\alpha=\sum_{j\in B}a_jy^j$, where
$\emptyset\neq B\sbq\{ 0,\ldots,p^m-1\}^n$ and $a_j\in K\setminus\{
0\}$ $\forall j\in B$. We have $\alpha\in K\setminus\{ 0\}$ iff $B=\{
0\}$, i.e. iff $\alpha =\alpha_0$. If $\alpha\notin K$, i.e. if
$B\neq\emptyset,\{ 0\}$, then we define $(l_\alpha,k_\alpha,r_\alpha
)$ as follows: $l_\alpha$ is the largest $l$ such that there is
$j=(j_0,\ldots,j_{n-1})\in B$ with $j_l\neq 0$; $k_\alpha$ is the
samllest $k$ such that there is $j=(j_0,\ldots,j_{n-1})\in B$ with
$p^k\| j_{l_\alpha}$; $r_\alpha$ is the largest $r$ with $p\nmid r$
such that there is $j=(j_0,\ldots,j_{n-1})\in B$ with
$j_{l_\alpha}=p^{k_\alpha}r$. We have $0\leq l_\alpha\leq n-1$, $0\leq
k_\alpha\leq m-1$ and $1\leq r_\alpha\leq p^{m-k_\alpha}-1$, with
$p\nmid r_\alpha$. We use the same order relation on triplets as in
Step 1.

We prove that if $\alpha\in J\setminus K$ and
$(l,k,r)=(l_\alpha,k_\alpha,r_\alpha )$ then $\alpha'=[\alpha,
x_k^{p^l}]\in J\setminus\{ 0\}$ and we have either $\alpha'\in K$ or
$(l',k',r'):=(l_{\alpha'},k_{\alpha'},r_{\alpha'})<(l,k,r)$. Then by
the same induction argument from Step 1 we obtain an element
$\alpha\in J\setminus\{ 0\}$ with $\alpha\in K$. 

We write $\alpha =\sum_{j\in B}a_jy^j$ with $a_j\in K\setminus\{
0\}$. Let $B'=\{ j=(j_0,\ldots,j_{n-1})\in B\mid p^k\| j_l\}$. If
$j=(j_0,\ldots,j_{n-1})\in B'$ then $j_l=p^kr_j$ with $p\nmid r_j$ and
we have $r=\max\{ r_j\mid j\in B'\}$. By a similar proof as in Step 1,
we have $\alpha'=\sum_{j\in B'}r_ja_jy^{j-p^k{\rm e}_l}$. (This time
we use the fact that $[y_h,x_k^{p^l}]=0$ if $h<l$,
$[y_l^{p^{k+1}},x_k^{p^l}]=0$ and $[y_l^{p^k},x_k^{p^l}]=1$.) Then
$\alpha'\in I\cap K[y]=J$ and the set $B$ of indices corresponding to
$\beta$ is $B'-p^k{\rm e}_l\neq\emptyset$ so $\alpha'\neq 0$. Then, by
the same reasoning as in Step 1, we get that $\alpha'\in K$ or
$(l',k',r'):=(l_{\alpha'},k_{\alpha'},r_{\alpha'})<(l,k,r)$. \qed

As a consequence of  Lemma 4.1, Corollary 4.3 and Lemma 4.4 we have:

\btm $A_{((a,b))}$ is a central simple algebra of degree $p^{mn}$.
\etm

\bdf We define $((\cdot,\cdot ))=((\cdot,\cdot
))_{K,p^m,p^n}:W_m(K)\times W_n(K)\to\Br (K)$ by
$((a,b))=[A_{((a,b))}(K)]$. 

In particular, if $m=n$ we denote $((\cdot,\cdot
))_{K,p^n,p^n}=((\cdot,\cdot ))_{K,p^n}$. If the field $K$ is fixed we
drop the $K$ from the index.
\edf

Since $m,n$ are fixed we write $((\cdot,\cdot ))$ instead of
$((\cdot,\cdot ))_{p^m,p^n}$.

\btm (i) $((a,b))=((a+F^nc,b+F^md))$ $\forall a,c\in
W_m(K),\, b,d\in W_n(K)$.  

(ii) $((\cdot,\cdot ))$ is bilinear.

(iii) $((a,b))_{p^m,p^n}=-((b,a))_{p^n,p^m}$. In particular,
$((\cdot,\cdot))_{p^n}$ is skew-symmetric.

(iv) If $m=n$ then $((a,bc))+((b,ac))+((c,ab))=0$.
\etm
\pf (i) By Lemma 2.20 we have
$$B_{m,n}(K)(x,y)=B_{m,n}(K)(x+c,y+d).$$
Then
\begin{multline*}
A_{((a,b))}(x,y)=B_{m,n}(K)(x+c,y+d)/(F^nx-a,F^my-b)\\
=B_{m,n}(K)(x+c,y+d)/(F^n(x+c)-(a+F^nc),
F^m(y+d)-(b+F^md))\\
=A_{((a+F^nc,b+F^md))}(x+c,y+d).
\end{multline*}
By taking the classes in the Brauer group we get
$((a,b))=((a+F^nc,b+F^md))$.

(iii) By Lemma 2.19 $B_{m,n}(K)(x,y)^{op}=B_{n,m}(K)(y,x)$. Hence
\begin{multline*}
A_{((a,b))_{p^m,p^n}}(K)(x,y)^{op}=
(B_{m,n}(K)(x,y)/(F^nx-a,F^my-b))^{op}\\
=B_{n,m}(K)(y,x)/(F^nx-a,F^my-b)=
A_{((b,a))_{p^n,p^m}}(K)(y,x).
\end{multline*}
By taking the classes in the Brauer group we get
$-((a,b))_{p^m,p^n}=((b,a))_{p^n,p^m}$. 

(ii) We prove first that $((0,0))=0$. To do this we note that
$A_{((0,0))}(K)=A_{((0,0))}(\FF_p)\otimes_{\FF_p}K$. But $\Br
(\FF_p)=0$ so $A_{((0,0))}(\FF_p)\cong M_{p^{mn}}(\FF_p)$. It follows
that $A_{((0,0))}(K)\cong M_{p^{mn}}(K)$ so $((0,0))=0$. 

We have 
$$A_{((a,b))}(x,y)\otimes_KA_{((c,d))}(z,t)=
C/(F^nx-a,F^my-b,F^nz-c,F^mt-d),$$
where $C:=B_{m,n}(K)(x,y)\otimes_K B_{m,n}(K)(z,t)$. By
Lemma 2.21 $C$ also writes as $C=B_{m,n}(K)(x+z,y)\otimes_K
B_{m,n}(K)(z,t-y)$. Also the relations $F^nx=a$,
$F^my=b$, $F^nz=c$, $F^mt=d$ are equivalent to
$F^n(x+z)=a+c$, $F^my=b$, $F^nz=c$, $F^m(t-y)=d-b$. Hence
\begin{multline*}
(F^nx-a,F^my-b,F^nz-c,F^mt-d)\\
=(F^n(x+z)-(a+c),F^my-b,F^nz-c,F^m(t-y)-(d-b)).
\end{multline*}
It follows that
$$A_{((a,b))}(x,y)\otimes_KA_{((c,d))}(z,t)=
A_{((a+c,b))}(x+z,y)\otimes_KA_{((c,d-b))}(z,t-y).$$
By taking classes in the Brauer group we get
$((a,b))+((c,d))=((a+c,b))+((c,d-b))$. In particular,
$((a,b))+((c,b))=((a+c,b))+((c,0))$. Similarly, we also have
$((c,b))+((a,b))=((c+a,b))+((a,0))$ so $((a,0))=((c,0))$. Hence
$((a,0))$ is independent off $a$. But $((0,0))=0$ so $((a,0))=0$
$\forall a\in W_m(K)$. Hence
$((a,b))+((c,b))=((a+c,b))+((c,0))=((a+c,b))$.

For the linearity in the second variable we use the skew-symmetry from
(iii). 

(iv) If $D=A_{((a,bc))}(x_1,y_1)\otimes_K
A_{((b,ac))}(x_2,y_2)\otimes_K A_{((c,ab))}(x_3,y_3)$ then $D=C/I$,
where
$C=B_n(K)(x_1,y_1)\otimes_KB_n(K)(x_2,y_2)\otimes_KB_n(K)(x_3,y_3)$
and $I$ is the ideal $(F^nx_1-a,F^ny_1-bc, F^nx_2-b,F^ny_2-ac,
F^nx_3-c,F^ny_3-ab)$.

But by Lemma 2.22 $C$ also writes as
$$C=B_n(K)(x_1,y_1-x_2x_3)\otimes_K
B_n(K)(x_2,y_2-x_1x_3)\otimes_K B_n(K)(x_3,y_3-x_1x_2).$$
Also the relations $F^nx_1=a$, $F^ny_1=bc$, $F^nx_2=b$, $F^ny_2=ac$,
$F^nx_3=c$ and $F^ny_3=ab$ are equivalent to $F^nx_1=a$,
$F^n(y_1-x_2x_3)=0$, $F^nx_2=b$, $F^n(y_2-x_1x_3)=0$, $F^nx_3=c$ and
$F^n(y_3-x_1x_2)=0$. Hence $I=(F^nx_1-a,F^n(y_1-x_2x_3),
F^nx_2-b,F^n(y_2-x_1x_3), F^nx_3-c,F^n(y_3-x_1x_2))$. It follows that
$D=C/I=A_{((a,0))}(x_1,y_1-x_2x_3)\otimes_K
A_{((b,0))}(x_2,y_2-x_1x_3)\otimes_K A_{((c,0))}(x_3,y_3-x_1x_2)$.

Thus $((a,bc))+((b,ac))+((c,ab))=[D]=((a,0))+((b,0))+((c,0))=0$. \qed

{\bf Remark.} The proof of Theorem 4.6(iv), using Lemma 2.22, follows
the idea in the case $n=1$ from [BK, 8.1.1], where it is refered as
``the most complicated fifth isomorphism''.
\medskip

As a consequence of Theorem 4.6(i) and (ii) and the fact that
$W_m(K)/F^n(W_m(K))$ and $W_n(K)/F^m(W_n(K))$ are $p^l$-torsion, we
get

\bco $((\cdot,\cdot ))_{p^m,p^n}$ is a bilinear defined as
$$((\cdot,\cdot ))_{p^m,p^n}:W_m(K)/F^n(W_m(K))\times
W_n(K)/F^m(W_n(K))\to\Brr{p^l}(K),$$
where $l=\min\{m,n \}$. 

In particular, if $m=n$ then $((\cdot,\cdot ))=((\cdot,\cdot ))_{p^n}$
is defined
$$((\cdot,\cdot )):W_n(K)/F^n(W_n(K))\times
W_n(K)/F^n(W_n(K))\to\Brr{p^n}(K).$$
\eco

See [B, Corollary 3.12].

Theorem 4.6(ii), (iii) and (iv) in the case $m=n$ are equivalent to:

\bpr There is a group morphism $\alpha_{p^n}:\Omega^1(W_n(K))/\diff
W_n(K)\to\Brr{p^n}(K)$ given by $a\diff b\mapsto ((a,b))_{p^n}$.
\epr

See [B, Proposition 3.6] and the following Remark.

\section{The adjoint property of Frobenius and Verschiebung}

In this section we prove that the operators $F$ and $V$ are adjoint
with respect to the symbols $((\cdot,\cdot ))_{p^m,p^n}$.

\blm Let $a\in W_m(K)$, $b\in W_n(K)$. Let $0\leq k\leq m$, $0\leq
l\leq n$. Same as in Proposition 3.1, we denote the generators of
$A_{((a,b))_{p^m,p^n}}$ as $x=(x',x'')$, $y=(y',y'')$, with
$x'=(x_0,\ldots,x_{k-1})$, $x''=(x_k,\ldots,x_{m-1})$,
$y'=(y_0,\ldots,y_{l-1})$, $y''=(y_l,\ldots,y_{n-1})$.

Similarly, we denote $a=(a',a'')$ and $y=(b',b'')$, with
$a'=(a_0,\ldots,a_{k-1})$, $a''=(a_k,\ldots,a_{m-1})$,
$b'=(b_0,\ldots,b_{l-1})$, $b''=(b_l,\ldots,b_{n-1})$.

Then in $A_{((a,b))_{p^m,p^n}}$ we have
$$\langle F^lx'',F^ky''\rangle
=A_{((a'',b''))_{p^{m-k},p^{n-l}}}(F^lx'',F^ky'').$$
\elm
\pf Since $x=(x',x'')$ and $a=(a',a'')$ the relation $F^nx=a$ from
$A_{((a,b))_{p^m,p^n}}$ also writes as $F^nx'=a'$ and
$F^nx''=a''$. Similarly $F^my=b$ writes as $F^my'=b'$ and
$F^my''=b''$. 

Let $C$ and $D$ be the subalgebras generated by $F^lx''$ and $F^ky''$
in $B_{m,n}(K)$ and $A_{((a,b))_{p^m,p^n}}$, respectively. Then the
projection $B_{m,n}(K)\to
A_{((a,b))_{p^m,p^n}}=B_{m,n}(K)/(F^nx-a,F^my-b)$ sends $C$ to $D$. So
we have a canonical surjective morphism $h:C\to D$. The relations
$F^nx''=a''$ and $F^my''=b''$ from $A_{((a,b))_{p^m,p^n}}$ also hold
in the sublagebra $D$.

By Proposition 3.1 in $B_{m,n}(K)$ we have
$C=B_{m-k,n-l}(K)(F^lx'',F^ky'')$. So if $z=(z_0,\ldots,z_{m-k-1})$
and $t=(t_0,\ldots,t_{n-l-1})$ are multivariables then we have an
isomorphism $f:B_{m-k,n-l}(K)(z,t)\to C$ given by $z\mapsto F^lx''$,
$t\mapsto F^ky''$. Then $g:=hf:B_{m-k,n-l}(K)(z,t)\to D$, given by
$z\mapsto F^lx''$, $t\mapsto F^ky''$, is a surjective morphism of
algebras. Since $g(F^{n-l}z-a'')=F^{n-l}(F^lx'')-a''=F^nx''-a''=0$ and
$g(F^{m-k}t-b'')=F^{m-k}(F^ky'')-b''=F^my''-b''=0$ we have
$F^{m-k}z-a'',F^{n-l}t-b''\in\ker g$. It follows that $g$ induces a
morphism $\overline g:B_{m-k,n-l}(K)(z,t)/(F^{n-l}z-a'',F^{m-k}t-b'')=
A_{((a'',b''))_{p^{m-k},p^{n-l}}}(z,t)\to D$. Since $g$ is surjective,
so is $\overline g$. But $A_{((a'',b''))_{p^{m-k},p^{n-l}}}(z,t)$ is a
simple algebra so in fact $\overline g$ is an isomorphism. Since
$\overline g$ is given by $z\mapsto F^nx''$, $t\mapsto F^my''$ we have
$D=A_{((a'',b''))_{p^{m-k},p^{n-l}}}(F^lx'',F^ky'')$, as claimed. \qed

From now on we regard truncated $p$-typical Witt vectors as classes of
full $p$-typical Witt vectors, i.e. we identify
$W_n(K)=W(K)/V^n(K)$. This has the advantage that we can switch
between truncations of different lenghts. Then in short notation the
group $W_m(K)/F^n(W_m(K))$ writes as $W(K)/(V^m,F^n)$, where by
$(V^m,F^n)$ we mean the group fenerated by the images of $V^m$ and
$F^n$. Similarly for $W_n(K)/F^m(W_n(K))$. Then $((\cdot,\cdot
))_{p^m,p^n}$ is defined as
$$((\cdot,\cdot ))_{p^m,p^n}:W(K)/(V^m,F^n)\times
W(K)/(V^n,F^m)\to\Brr{p^l}(K),\quad l=\min\{ m,n\}.$$
(See [B, 3.14].)

\btm Let $a,b\in W(K)$.

(i) If $m\geq 0$, $n\geq 1$ then
$((Fa,b))_{p^m,p^n}=((a,Vb))_{p^m,p^n}=((a,b))_{p^m,p^{n-1}}$.

(ii) If $m\geq 0$, $n\geq 1$ then
$((a,Fb))_{p^m,p^n}=((Va,b))_{p^m,p^n}=((a,b))_{p^{m-1},p^n}$.

Recall that if $m$ or $n=0$ then $((a,b))_{p^m,p^n}=0$.
\etm

Note that (ii) follows from (i) by using the skew-symmetry from
Theorem 4.6(iii). So we only have to prove (i).
\medskip

{\bf Idea of the proof}

If $C=A_{((Fa,b))_{p^m,p^n}}$ or $A_{((a,Vb))_{p^m,p^n}}$ then we find
the subalgebras $A,B\sbq C$ with $[A,B]=0$, $A\cong
A_{((a,b))_{p^m,p^{n-1}}}$ and $B\equiv M_{p^m}(K)$. Since $[A,B]=0$,
by the universal property of the tensor product there is a morphism
$f:A\otimes_KB\to C$ given by
$\alpha\otimes\beta\mapsto\alpha\beta$. Since $A$ and $B$ are
c.s.a. so is $A\otimes_KB$. Therefore $f$ is injective. But by
Corollary 4.3 we have $\dim_KA=p^{2m(n-1)}$ and $\dim_KB=p^{2m}$ so
$\dim_KA\dim_KB=p^{2mn}=\dim_KC$. Hence $f$ is an isomorphism so
$C\cong A\otimes_KB$. It follows that
$[C]=[A]+[B]=((a,b))_{p^m,p^{n-1}}+0=((a,b))_{p^m,p^{n-1}}$,
i.e. $((Fa,b))_{p^m,p^n}$ or
$((a,Vb))_{p^m,p^n}=((a,b))_{p^m,p^{n-1}}$, accordingly.
\medskip

{\bf Proof of $((Fa,b))_{p^m,p^n}=((a,b))_{p^m,p^{n-1}}$}

We write $a=(a_0,\ldots,a_{m-1})$ and
$b=(b_0,\ldots,b_{n-1})=(b',b_{n-1})$, where
$b'=(b_0,\ldots,b_{n-2})$. Since $b'$ is the truncation of $b$ in
$W_{n-1}(K)$, by $((a,b))_{p^m,p^{n-1}}$ we mean
$((a,b'))_{p^m,p^{n-1}}$. 

If $x,y$ are the generators of $A_{((Fa,b))_{p^m,p^n}}$ then
$y=(y',y_{n-1})$, where $y'=(y_0,\ldots,y_{n-2})$. We have
$C=B_{m,n}(K)/(F^nx-Fa,F^my-b)$. 

We define the following ideals:
\begin{center}
$I=$ the ideal $(F^nx-Fa,F^my-b)$ of $B_{m,n}(K)$

$I'=$ the ideal $(F^nx-Fa,F^my'-b')$ of $B_{m,n-1}(K)$

$J=$ the ideal $(F^{n-1}x-a,F^my'-b')$ of $B_{m,n-1}(K)$
\end{center}
Note that $A_{((Fa,b))_{p^m,p^n}}=B_{m,n}(K)/I$ and
$A_{((a,b'))_{p^m,p^{n-1}}}=B_{m,n-1}(K)/J$. Also note that $F^my=b$
implies $F^my'=b'$ so $I'\sbq I$ and $F^{n-1}x=a$ implies $F^nx=Fa$ so
$I'\sbq J$.

\blm We have $I'=I\cap B_{m,n-1}(K)$.
\elm
\pf Since $I'\sbq I$ the inclusion map $B_{m,n-1}(K)\sbq B_{m,n}(K)$
induces a morphism $f:B_{m,n-1}(K)/I'\to B_{m,n}(K)/I$. Then $I'=I\cap
B_{m,n-1}(K)$ is equivalent to the injectivity of $f$. 

By Corollary 4.3 we have that $S=\{ x_0^{i_0}\cdots
x_{m-1}^{i_{m-1}}y_0^{j_0}\cdots y_{n-1}^{j_{n-1}}\mid 0\leq i_q\leq
p^n-1,\, 0\leq j_r\leq p^m-1\}$ is a basis of
$B_{m,n}(K)/I=A_{((Fa,b))_{p^m,p^n}}$. We prove a similar result for
$B_{m,n-1}(K)/I'$. To do this we replicate the proofs of Lemma 4.2 and
Corollary 4.3.

By Corollary 3.9
$Z(B_{m,n-1}(K))=K[x_0^{p^{n-1}},\ldots,x_{m-1}^{p^{n-1}},
y_0^{p^m},\ldots,y_{n-2}^{p^m}]$. Let $I'_x$ be the ideal $(F^nx-Fa)$
of $K[x]$ and let $I'_{y'}$ be the ideal $(F^my'-b')$ of
$K[y']$. Since the generators
$x_0^{p^n}-a_0^p,\ldots,x_{m-1}^{p^n}-a_{m-1}^p$ of $I'_x$ and the
generators $y_0^{p^m}-b_0,\ldots,y_{n-2}^{p^m}-b_{n-2}$ of $I'_{y'}$
belong to $Z(B_{m,n-1}(K))$ and $I'=(F^nx-Fa,F^my'-b')$, by the same
proof as for the Lemma 4.2, we get that the isomorphism $\mu
:K[x]\otimes_KK[y']\to B_{m,n-1}(K)$ induces an isomorphism $\mu
:K[x]/I'_x\otimes_KK[y']/I'_{y'}\to B_{m,n-1}(K)/I'$. Then we proceed
as  for Corollary 4.3. The products $x_0^{i_0}\cdots
x_{m-1}^{i_{m-1}}$ with $0\leq i_q\leq p^n-1$ are a basis for
$K[x]/I'_x$ and the products $y_0^{j_0}\cdots y_{n-2}^{j_{n-2}}$ with
$0\leq j_r\leq p^m-1$ are a basis for $K[y']/I'_{y'}$. Hence $S'=\{
x_0^{i_0}\cdots x_{m-1}^{i_{m-1}}y_0^{j_0}\cdots y_{n-2}^{j_{n-2}}\mid
0\leq i_q\leq p^n-1,\, 0\leq j_r\leq p^m-1\}$ is a basis for
$B_{m,n-1}/I'$.

To conclude the proof, note that the elements of the basis $S'$ are
sent by $f$ to simlar elements of $B_{m,n}(K)/I$, which are part of
the basis $S$ so they are linearly independent. It follows that $f$ is
injective. \qed

\blm If $B$ is the subalgebra $\langle F^{n-1}x,y_{n-1}\rangle$ of
$A_{((Fa,b))_{p^m,p^n}}$ then $B\cong M_{p^m}(K)$.
\elm
\pf We use Lemma 5.1 with $k=0$, $l=n-1$ and we get
$B=A_{((Fa,b_{n-1}))_{p^m,p}}(F^{n-1}x,y_{n-1})$. By Theorem 4.5 $B$
is a c.s.a. of degree $p^m$. By Theorem 4.6(i) we have
$[B]=((Fa,b_{n-1}))_{p^m,p}=((0,b_{n-1}))_{p^m,p}=0$. Hence $B\cong
M_{p^m}(K)$. \qed

\blm If $A$ is the centralizer $C(B)$ of $B$ in
$A_{((Fa,b))_{p^m,p^n}}$ then $A\cong A_{((a,b'))_{p^m,p^{n-1}}}$.
\elm
\pf Since $B\cong M_{p^m}(K)$ is simple, by the double centralizer
theorem, $A$ is also simple and
$\dim_KB\dim_KA=\dim_KA_{((Fa,b))_{p^m,p^n}}$. But $\dim_KB=p^{2m}$
and, by Corollary 4.3, $\dim_KA_{((Fa,b))_{p^m,p^n}}=p^{2mn}$. Thus
$\dim_KA=p^{2m(n-1)}$.

Since $B\spq\langle F^{n-1}x\rangle=\langle
x_0^{p^{n-1}},\ldots,x_{m-1}^{p^{n-1}}\rangle$ we have $A\sbq
A_0:=C(x_0^{p^{n-1}},\ldots,x_{m-1}^{p^{n-1}})$. But by the case
$n'=0$ of Proposition 3.8 we have
$A_0=K[x][y_0,\ldots,y_{n-2},y_{n-1}^{p^m}]$. (See also the remarks
following Lemma 3.5.) But in $A_{((Fa,b))_{p^m,p^n}}$ we have
$y_{n-1}^{p^m}=b_{n-1}\in K$ so
$A_0=K[x][y_0,\ldots,y_{n-2}]=K[x][y']$. But $x$ and $y'$ are the
generators of $B_{m,n-1}(K)$ so $A_0$ is the image of
$B_{m,n-1}(K)\sbq B_{m,n}(K)$ in
$A_{((Fa,b))_{p^m,p^n}}=B_{m,n}(K)/I$. Since $I\cap B_{m,n-1}(K)=I'$
we have $A_0=B_{m,n-1}(K)/I'$.

Since $I'\sbq J$ we have a surjective morphism $f:B_{m,n-1}(K)/I'\to
B_{m,n-1}(K)/J$ i.e. $f:A_0\to A_{((a,b'))_{p^m,p^{n-1}}}$. We denote
by $g:A\to A_{((a,b'))_{p^m,p^{n-1}}}$ the restriction $f_{|A}$. Since
$A$ is simple $g$ is injective. But by Corollary 4.3
$\dim_KA_{((a,b))_{p^m,p^{n-1}}}=p^{2m(n-1)}=\dim_KA$. It follows that
$g$ is an isomorphism so $A\cong A_{((a,b'))_{p^m,p^{n-1}}}$. \qed

By Lemmas 5.4 and 5.5 we have the subalgebras $A,B$ of
$A_{((Fa,b))_{p^m,p^n}}$ with $A\cong A_{((a,b'))_{p^m,p^{n-1}}}$,
$B\cong M_{p^m}(K)$ and , since $A=C(B)$, $[A,B]=0$. As seen in the
preamble to our proof, this implies
$((Fa,b))_{p^m,p^n}=((a,b'))_{p^m,p^{n-1}}$.

\medskip
{\bf Proof of $((a,Vb))_{p^m,p^n}=((a,b))_{p^m,p^{n-1}}$}

If $b=(b_0,\ldots,b_{n-2})\in W_{n-1}(K)$ then
$Vb=(0,b)=(0,b_0,\ldots,b_{n-2})\in W_n(K)$. The generators of
$A_{((a,Vb))}$ are $x=(x_0,\ldots,x_{m-1})$ and $y=(y_0,y'')$, where
$y''=(y_1,\ldots,y_{n-1})$. Then the relation $F^my=Vb$ means
$y_0^{p^m}=0$ and $F^my''=b$, i.e. $y_j^{p^m}=b_{j-1}$ when $j\geq
1$.

\blm Let $D$ be a division algebra and let $r\geq 0$. If $y\in M_r(D)$
is nilpotent then $y^r=0$.
\elm
\pf By [K, chapter II, \S 2, Example (1)] all simple left
$M_r(D)$-modules are $\cong D^r$. Since $M_r(D)$ is semismple every
left $M_r(D)$-module is a direct sum of simple modules, i.e. of copies
of $D^r$. In particular, $M_r(D)$ as a left $M_r(D)$-module is the
direct sum of $r$ copies of $D^r$. If $M_r(D)=I_0\sp
I_1\space\cdots\sp I_l$ is a strictly descending sequence of left
ideals then each $I_k$ writes as a direct sum of $r_k$ copies of
$D^r$, with $r=r_0>r_1>\cdots >r_l\geq 0$. It follows that $l\leq r$.

We consider the descending sequence of ideals $M_r(D)\spq M_r(D)y\spq
M_r(D)y^2\spq\cdots$. As seen above this sequence cannot be strictly
decreasing and, moreover, if $l$ is the smallest index such that
$M_r(D)y^l=M_r(D)y^{l+1}$ then $l\leq r$. Then for every $k\geq l$ we
have $M_r(D)y^ly^{k-l}=M_r(D)y^{l+1}y^{k-l}$,
i.e. $M_r(D)y^k=M_r(D)y^{k+1}$. Hence
$M_r(D)y^l=M_r(D)y^{l+1}=M_r(D)y^{l+2}=\cdots$. But $y$ is nilpotent
so $y^N=0$ for $N$ large enough. It follows that
$M_r(D)y^l=M_r(D)y^N=0$, which implies $y^l=0$. But $r\geq l$ so
$y^r=0$. \qed

\blm If $A$ is the subalgebra $\langle Fx,y''\rangle$ of
$A_{((a,Vb))_{p^m,p^n}}$ then $A\cong A_{((a,b))_{p^m,p^{n-1}}}$.
\elm
\pf We use Lemma 5.1 with $k=0$, $l=1$. Since $Vb=(0,b)$ we have
$A=A_{((a,b))_{p^m,p^{n-1}}}(Fx,y'')$. \qed

\blm If $B$ is the centralizer $C(A)$ of $A$ in
$A_{((a,Vb))_{p^m,p^n}}$ then $B\cong M_{p^m}(K)$.
\elm
\pf Since $A\cong A_{((a,b))_{p^m,p^{n-1}}}$ is simple, by the double
centralizer theorem, $B$ ia also simple and
$\dim_KA\dim_KB=\dim_KA_{((a,Vb))_{p^m,p^n}}$. But by Corollary 4.3 we
have $\dim_KA=p^{2m(n-1)}$ and
$\dim_KA_{((a,Vb))_{p^m,p^n}}=p^{2mn}$. Thus $\dim_KB=p^{2m}$.

We have $A=\langle Fx,y''\rangle =\langle
x_0^p,\ldots,x_{m-1}^p,y_1,\ldots,y_{n-1}\rangle$. By Proposition 3.2
we have $[y_0,x_i^p]=0$ $\forall i$ so $y_0\in C(A)=B$. If $\alpha\in
Z(B)$ then $\alpha$ commutes with $y_0\in B$. But we also have
$\alpha\in B=C(A)$ so it commutes with the generators
$x_0^p,\ldots,x_{m-1}^p,y_1,\ldots,y_{n-1}$ of $A$. Therefore
$\alpha\in C(x_0^p,\ldots,x_{m-1}^p,y_0,y_1,\ldots,y_{n-1})=
K[x_0^{p^n},\ldots,x_{m-1}^{p^n}][y_0,y_1^{p^m},\ldots,y_{n-1}^{p^m}]$. (See
Proposition 3.8.) But in $A_{((a,Vb))_{p^m,p^n}}$ we have
$x_i^{p^n}=a_i\in K$ and, if $j\geq 1$, $y_j^{p^m}=b_{j-1}\in K$. So
in fact $\alpha\in K[y_0]$. In conclusion $K\sbq Z(B)\sbq K[y_0]$. But
$B$ is simple so $Z(B)$ is a field. Suppose that $Z(B)\neq K$ and let
$\alpha\in Z(B)\setminus K\sbq K[y_0]\setminus K$. Then $\alpha
=\alpha_0+y_0P(y_0)$ for some $\alpha_0\in K$ and $P\in K[X]$. Since
$\alpha\not\in K$ we have $y_0P(y_0)\neq 0$. Since $\alpha_0\in K\sbq
Z(B)$ we have $y_0P(y_0)=\alpha -\alpha_0\in Z(B)$. But $y_0^{p^m}=0$
so in $Z(B)$ we have $y_0P(y_0)\neq 0$ but $(y_0P(y_0))^{p^m}=0$. Thus
$Z(B)$ is not a field. Contradiction. So $Z(B)=K$. Since $B$ is also
simple, it is a c.s.a. 

Now $B$ is a c.s.a. with $\dim_KB=p^{2m}$ so $\deg B=p^m$. It follows
that $B\cong M_r(D)$, where $D$ is a central division algebra with
$r\deg D=p^m$. Assume that $B\not\cong M_{p^m}(K)$. It follows that
$r<p^m$. Since $y_0\in B\cong M_r(D)$ is nilpotent, by Lemma 5.6 we
have $y_0^r=0$. But this is impossible since $r<p^m$ so $y_0^r$ is an
element in the basis of $A_{((a,Vb))_{p^m,p^n}}$ from Corollary
4.3. Hence $B\cong M_{p^m}(K)$. \qed 

By Lemmas 5.7 and 5.8 we have the subalgebras $A,B$ of
$A_{((a,Vb))_{p^m,p^n}}$ with $A\cong A_{((a,b))_{p^m,p^{n-1}}}$,
$B\cong M_{p^m}(K)$ and , since $B=C(A)$, $[A,B]=0$. As seen in the
preamble to our proof, this implies
$((a,Vb))_{p^m,p^n}=((a,b))_{p^m,p^{n-1}}$. 

We are now able to recover all the properties of the symbols
$((\cdot,\cdot ))_{p^m,p^n}$ defined in [B], except those involving
$[\cdot,\cdot)_{p^n}$, such as [B, Proposition 3.8, Corollary 3.9,
Definition 2, Proposition 3.13]. 

\bco (i) If $a,b\in W(K)$ then for every $m,n,i,j,k,l\geq 0$ we have
$$((F^iV^ja,F^kV^lb))_{p^m,p^n}=\begin{cases}
((a,b))_{p^{m-j-k},p^{n-i-l}}&\text{if }m>j+k,\, n>i+l\\
0&\text{otherwise}
\end{cases}.$$

(ii) If $l\geq m,n$ then $((a,b))_{p^m,p^n}=((V^{l-m}a,V^{l-n}b))_{p^l}$. 

(iii) If $m\geq n$ then
$((a,b))_{p^n}=((V^{m-n}a,V^{m-n}b))_{p^m}=p^{m-n}((a,b))_{p^m}$.
\eco
\pf (i) Since $F$ and $V$ are adjoint with respect to $((\cdot,\cdot
))_{p^m,p^n}$ we have
$((F^iV^ja,F^kV^lb))_{p^m,p^n}=((V^{j+k}a,V^{i+l}b))_{p^m,p^n}$. If
$j+k\geq m$ then $V^{j+k}a=0$ in $W_m(K)$. If $i+l\geq n$ then
$V^{i+l}b=0$ in $W_n(K)$. In both cases
$((V^{j+k}a,V^{i+l}b))_{p^m,p^n}=0$. If $m>j+k$ and $n>i+l$ then
$((V^{j+k}a,V^{i+l}b))_{p^m,p^n}=((a,b))_{p^{m-j-k},p^{n-i-l}}$
follows by repeated use of the relations
$((Va,b))_{p^m,p^n}=((a,b))_{p^{m-1},p^n}$ and
$((a,Vb))_{p^m,p^n}=((a,b))_{p^m,p^{n-1}}$. 

(ii) follows directly from (i) since $((\cdot,\cdot
))_{p^l}=((\cdot,\cdot ))_{p^l,p^l}$. Similarly for the same equality
from (iii). For the second equality we use the adjoint property of $F$
and $V$ and we get
$((V^{m-n}a,V^{m-n}b))_{p^m}=((F^{m-n}V^{m-n}a,b))_{p^m}
=((p^{m-n}a,b))_{p^m}=p^{m-n}((a,b))_{p^m}$. \qed

We now state the representation theorem for $\Brr{p^n}(K)$ from [B,
Theorem 4.10] in terms of the new symbols $((\cdot,\cdot ))_{p^n}$ we
introduced here.

\btm We have an isomorphism $\alpha_{p^n}:G_n\to\Brr{p^n}(K)$, where
$$G_n=\Omega^1(W_n(K))/(Fa\diff b-a\diff Vb\mid a,b\in W_n(K),~\wp
([a])\dlog b\mid a,b\in K,\, b\neq 0),$$
given by $a\diff b\mapsto ((a,b))_{p^n}$.
\etm
\pf As seen in the proof of [B, Theorem 4.10] the ingredients we need
are the following:

The fact that there is a linear map
$\alpha_{p^n}:\Omega^1(W_n(K))\to\Brr{p^n}(K)$ given by $a\diff
b\mapsto ((a,b))_{p^n}$. This follows from Porposition 4.8.

The fact that $F$ and $V$ are adjoint with respect to $((\cdot,\cdot
))_{p^n}$, which follows from Theorem 5.2.

The relation $\alpha_{p^n}(\wp([a])\dlog b)=0$ $\forall a,b\in K$,
$b\neq 0$. This will follow from Lemma 5.11 bellow.

The relation $((a,b))_{p^n}=((Va,Vb))_{p^{n+1}}$ $\forall a,b\in
W_n(K)$, which follows from Lemma 5.9(iii).

The induction step $n=1$. This follows from the fact that the symbol
$((\cdot,\cdot ))_p$ introduced here coincides with the symbol defined
in [B], which we proved in the introduction. Then the map
$\alpha_1:G_1\to\Brr p(K)$ we introduced here coincides with the one
from [B], which we know it is an isomorphism. Note that the induction
step $n=1$ of [B, Theorem 4.10] is just [GS, Theorem 9.2.4].\qed

\blm If $a\in W_n(K)$ and $b\in\kk$ then $\alpha_{p^n}(\wp (a)\dlog
[b])=0$.
\elm
\pf We have $\wp (a)\dlog [b]=\wp (a)[b]^{-1}\diff [b]=Fa[b]^{-1}\diff
[b]-a[b]^{-1}\diff [b]$ so $\alpha_{p^n}(\wp (a)\dlog
[b])=((Fa[b]^{-1},[b]))_{p^n}-((a[b]^{-1},[b]))_{p^n}$. Hence we must
prove that $((Fa[b]^{-1},[b]))_{p^n}=((a[b]^{-1},[b]))_{p^n}$.

We use Corollary 5.9(iii) and the adjoint property of $F$ and $V$ and
we get 
\begin{multline*}
((a[b]^{-1},[b]))_{p^n}=((V(a[b]^{-1}),V[b]))_{p^{n+1}}=
((F(a[b]^{-1}),F[b]))_{p^{n+1}}\\
=\alpha_{p^{n+1}}(F(a[b]^{-1})\diff
F[b])=\alpha_{p^{n+1}}(Fa[b]^{-p}\diff
[b]^p)=\alpha_{p^{n+1}}(Fa[b]^{-p}p[b]^{p-1}\diff [b])\\
=p\alpha_{p^{n+1}}(Fa[b]^{-1}\diff [b])=p((Fa[b]^{-1}\diff
[b]))_{p^{n+1}}=((Fa[b]^{-1}\diff [b]))_{p^n}.
\end{multline*}
\qed




\appendix

\section{$B_{m,n}(R)$ as an Azumaya algebra over its center}

Let $R$ be a ring of characteristic $p$. We prove that
$B_{m,n}(R)$ is an Azumaya algebra over it's center,
$Z(B_{m,n}(R))=
R[x_0^{p^n},\ldots,x_{m-1}^{p^n},y_0^{p^m},\ldots,y_{n-1}^{p^m}]$. This
is an analogue of the similar result involving the usual Weyl algebras
in positive characteristic, proved in [R]. 

We proceed like in [BK2, \S 3.2] and we use the following alternative
definition for Azumaya algebras. An $S$-algebra $C$ is called an
Azumaya algebra of degree $k$ if there is a flat $S$-ring $S'$ such
that $C\otimes_SS'\cong M_k(S')$ for some $k\geq 1$.

We consider the polynomial algebra $S=R[\alpha,\beta ]$, where $\alpha
=(\alpha_0,\ldots,\alpha_{m-1})$, $\beta
=(\beta_0,\ldots,\beta_{n-1})$, and the algebra
$C=B_{m,n}(S)/(F^nx-\alpha,F^my-\beta)$. Note that the
relations $F^nx=\alpha$, $F^my=\beta$ can be written as
$x_i^{p^n}=\alpha_i$ $\forall i$, $y_j^{p^m}=\beta_j$ $\forall j$.

\blm There is an isomorphism of $R$-algebras between $B_{m,n}(R)$
and $C$. Also $S\sbq C$ and $S=Z(C)$.
\elm
\pf As an $R$ algebra $C$ is generated by $\alpha$, $\beta$, $x$, $y$,
with the relations $[x_i,x_j]=0$, $[y_i,y_j]=0$
$[y_j,x_i]=c_{i,j}(x,y)$, $x_i^{p^n}=\alpha_i$, $y_j^{p^m}=\beta_j$
and the commutativity relations between each entry of $\alpha$ and
$\beta$ and all the other generators. For $B_{m,n}(R)$ we have the
generators $x$, $y$ and the relations $[x_i,x_j]=0$, $[y_i,y_j]=0$
$[y_j,x_i]=c_{i,j}(x,y)$. Since the relations among generators in
$B_{m,n}(R)$ also hold in $C$ there is a morphism of $R$-algebras
$f:B_{m,n}(R)\to C$ given by $x\mapsto x$, $y\mapsto
y$. Conversely, we have a morphism $g:C\to B_{m,n}(R)$ given by
$x\mapsto x$, $y\mapsto y$, $\alpha\mapsto F^n(x)$,
$\beta\mapsto F^m(y)$. Such morphism exists because the relations
among the generators from $C$ also hold in $B_{m,n}(R)$. Indeed, the
relations $[x_i,x_j]=0$, $[y_i,y_j]=0$ and $[y_j,x_i]=c_{i,j}(x,y)$
from $C$ are the same in $B_{m,n}(R)$. The relations
$x_i^{p^n}=\alpha_i$ and $y_j^{p^m}=\beta_j$ correspond to
$x_i^{p^n}=x_i^{p^n}$ and $y_j^{p^m}=y_j^{p^m}$ and the commutativity
relations involving the generators $\alpha_i$ and $\beta_j$ correspond
to commutativity relations involving $x_i^{p^n}$ and $y_j^{p^m}$. But
these follow from $x_i^{p^n},y_j^{p^m}\in Z(B_{m,n}(R))$. Obviously
$g\circ f=1_{B_{m,n}(R)}$, as it is given by $x_i\mapsto x_i$,
$y_j\mapsto y_j$ and $f\circ g=1_C$, as it is given by $x_i\mapsto
x_i$, $y_j\mapsto y_j$, $\alpha_i\mapsto x_i^{p^n}=\alpha_i$ and
$\beta_j\mapsto y_j^{p^m}=\beta_j$. Hence $f$ and $g$ are isomorphisms
inverse to each other.

Since the products $x^iy^j$ are linealy independent in
$B_{m,n}(R)$ we have that $Z(B_{m,n}(R))=
R[x_0^{p^n},\ldots,x_{m-1}^{p^n},y_0^{p^m},\ldots,y_{n-1}^{p^m}]$
holds strictly. Since $f(x_i^{p^n})=x_i^{p^n}=\alpha_i$ and
$f(y_j^{p^m})=y_j^{p^m}=\beta_j$ this implies that
$Z[C]=R[\alpha_0,\ldots,\alpha_{m-1},\beta_0,\ldots,\beta_{n-1}]$
strictly. Hence $Z(C)=S$. More precisely, $Z(C)$ is the image of $S$
in $C$. But the strictness property means that the monomials
$\alpha^i\beta^j$ with $i\in\II_m$, $j\in\II_n$ are linearly
independent in $C$. So the map
$S=R[\alpha_0,\ldots,\alpha_{m-1},\beta_0,\ldots,\beta_{n-1}]\to C$ is
an embedding, i.e. $S\sbq C$. \qed

\btm $B_{m,n}(R)$ is an Azumaya algebra of degre $p^{mn}$ over its
center.
\etm
\pf In the view of Lemma A.1, we must prove that $C$ is an Azumaya
algebra over $S$. We consider the multi-radical extension $S\sbq S'$,
with $S'=R[\eta,\theta]$, where $\eta =(\eta_0,\ldots,\eta_{m-1})$ and
$\theta =(\theta_0,\ldots,\theta_{n-1})$ satisfy
$\eta_i^{p^n}=\alpha_i$, $\theta_j^{p^m}=\beta_j$, i.e. $F^n\eta
=\alpha$, $F^m\theta =\beta$. Now $S'$ is a free $S$-module with
the basis $\eta_0^{i_0}\cdots\eta_{m-1}^{i_{m-1}}
\theta_0^{j_0}\cdots\theta_{n-1}^{j_{n-1}}$, with $0\leq i_k\leq
p^n-1$, $0\leq j_l\leq p^m-1$. Hence $S'$ is a faithfully flat
extension of $S$. Thus it suffices to prove that
$C':=C\otimes_SS'\cong M_{p^{mn}}(S')$. 

We use techniques that are similar to those from the proof of Theorem
4.6(i). By Lemma 2.20 $B_{m,n}(S')(x,y)=B_{m,n}(S')(x-\eta,y-\theta )$.
The relations $F^nx=\alpha =F^n\eta$ and $F^my=\beta
=F^m\theta$ are equivalent to $F^n(x-\eta )=0$ and
$F^m(y-\theta )=0$. So we have the equality of ideals
$(F^nx-\alpha,F^my-\beta )=(F^n(x-\eta ),F^m(y-\theta
))$. Then $C'=B_{m,n}(S')(x,y)/(F^nx-\alpha,F^my-\beta
)$ also writes as $C'=B_{m,n}(S')(x-\eta,y-\theta )/(F^n(x-\eta
),F^m(y-\theta ))$. It follows that $C'\cong
B_{m,n}(S')/(F^nx,F^my)$. Then we have
$C'=C_0\otimes_{\FF_p}S'$, where
$C_0=B_{m,n}(\FF_p)/(F^nx,F^my)=A_{((0,0))_{p^m,p^n}}(\FF_p)$. Hence
$C_0$ is a c.s.a. of degree $p^{mn}$ over $\FF_p$. But $\Br (\FF_p)=0$
so $C_0\cong M_{p^{mn}}(\FF_p)$. It follows that
$C'=C_0\otimes_{\FF_p}S'\cong M_{p^{mn}}(S')$. \qed


{\bf References}
\smallskip







[B] C.N. Beli, "A representation theorem for the $p^n$ torsion of the
Brauer group in characteristic $p$", preprint. Available on arXiv at:
\begin{verbatim*}
https://arxiv.org/abs/1711.00980
\end{verbatim*}

[BK1] A. Belov-Kanel, M. Kontsevich, "Automorphisms of the Weyl
Algebra", Letters in Mathematical Physics 74, No. 2, 181-199 (2005).

[BK2] A. Belov-Kanel, M. Kontsevich, "Jacobian Conjecture is stably
equivalent to Dixmier Conjecture", Moscow Mathematical Journal 7,
No. 2, 209-218 (2007).

[GS] P. Gille, T. Szamuely, "Central simple algebras and Galois
cohomology", Cambridge Studies in Advanced Mathematics 101, Cambridge
University Press (2006).

[H] Michiel Hazewinkel, "Formal  groups  and  applications", Pure and
Applied Mathematics,vol. 78, Academic Press Inc. [Harcourt Brace
Jovanovich Publishers], New York, (1978). 

[K] A. Knapp, "Advanced Algebra", Springer Verlag (2007).

[KMRT] M.-A. Knus, A. Merkurjev, M. Rost, and J.-P. Tignol, ``The Book
of Involutions'', AMS Colloquium Publications, Vol. 44 (1998).

[R] P. Revoy, "Alg\`ebres de Weyl en caract\'eristique p", Comptes
Rendus Hebdomadaires des S\'eances de l’Acad\'emie des Sciences,
S\'erie A 276, 225-228 (1973).
\bigskip

Institute of Mathematics Simion Stoilow of the Romanian Academy, Calea
Grivitei 21, RO-010702 Bucharest, Romania.

E-mail address: Constantin.Beli@imar.ro

\end{document}